\documentclass[a4paper,english]{amsart}

\usepackage[T1]{fontenc}
\usepackage[latin9]{inputenc}

\usepackage[pdftex]{graphicx}

\usepackage{a4wide}

\usepackage{babel}

\usepackage{amsmath}
\usepackage{amsthm}
\usepackage{amssymb}

\usepackage{enumitem} 

\newcommand{\R}{\mathbb {R}}
\newcommand{\N}{\mathbb {N}}
\def\pref#1{(\text{\ref{#1}})}
\newcommand{\plap}{\Delta_p}
\DeclareMathOperator{\dvg}{div}
\newcommand{\sobspace}{W^{1,p}_0(\Omega)}
\newcommand{\sobdual}{W^{-1,q}(\Omega)}
\newcommand{\dualpairing}[2]{\left\langle #1,#2\right\rangle}
\DeclareMathOperator{\per}{Per}
\newcommand{\cset}{\mathcal{C}}

\theoremstyle{plain}
\newtheorem{plapnum_proposition}{Proposition}[section]
\newtheorem{plapnum_lemma}[plapnum_proposition]{Lemma}
\newtheorem*{plapnum_aux_lemma}{Auxiliary Lemma}
\theoremstyle{remark}
\newtheorem*{plapnum_remark}{Remark}

\newcommand{\pdfpngextension}{png}
\newcommand{\dsuonegraphs}[2]{%
  \begin{picture}(30,57)
    \put(0,0){\makebox(30,4){$p=#1$}}
    \put(2,26){\includegraphics[width=26mm, clip=true, viewport= 700 300 1800 1580]{figs/D_disk/#2.jpg}}
    \put(3,7){\includegraphics[width=24mm]{figs/D_disk/#2.\pdfpngextension}}
  \end{picture}
  }
\newcommand{\dsutwographs}[2]{%
  \begin{picture}(30,62)
    \put(0,0){\makebox(30,4){$p=#1$}}
    \put(2,31){\includegraphics[width=26mm, clip=true, viewport= 700 300 1800 1580]{figs/D_disk/#2.jpg}}
    \put(3,7){\includegraphics[width=24mm]{figs/D_disk/#2.\pdfpngextension}}
  \end{picture}
  }
\newcommand{\dsemgraphs}[2]{%
  \begin{picture}(30,54)
    \put(2,24){\includegraphics[width=26mm, clip=true, viewport= 700 300 1800 1580]{figs/D_disk/#2.jpg}}
    \put(3,0){\includegraphics[width=24mm]{figs/D_disk/#2.\pdfpngextension}}
    \put(0,46){#1}
  \end{picture}
  }
\newcommand{\squonegraphs}[2]{%
  \begin{picture}(30,56)
    \put(0,0){\makebox(30,4){$p=#1$}}
    \put(1,29){\includegraphics[width=28mm, clip=true, viewport= 600 360 1890 1560]{figs/D_square2/#2.jpg}}
    \put(3,7){\includegraphics[width=24mm]{figs/D_square2/#2.\pdfpngextension}}
  \end{picture}
  }
\newcommand{\squtwographs}[2]{%
  \begin{picture}(30,60)
    \put(0,0){\makebox(30,4){$p=#1$}}
    \put(1,34){\includegraphics[width=28mm, clip=true, viewport= 600 360 1890 1560]{figs/D_square2/#2.jpg}}
    \put(3,7){\includegraphics[width=24mm]{figs/D_square2/#2.\pdfpngextension}}
  \end{picture}
  }
\newcommand{\reutwographs}[2]{%
  \begin{picture}(37.5,33)
    \put(0,0){\makebox(37.5,4){$p=#1$}}
    \put(3.75,6){\scalebox{-1}[1]{\includegraphics[width=30mm]{figs/D_rectangle_2_1.75/#2.\pdfpngextension}}}
  \end{picture}
  }
\newcommand{\troneuonegraphs}[2]{%
  \begin{picture}(30,56)
    \put(0,30){\includegraphics[width=30mm, clip=true, viewport= 650 475 1680 1355]{figs/D_triangle1/#2.jpg}}
    \put(1.5,3){\includegraphics[width=27mm]{figs/D_triangle1/#2.\pdfpngextension}}
    \put(0,0){\makebox(30,4){$p=#1$}}
  \end{picture}
  }
\newcommand{\troneutwographs}[2]{%
  \begin{picture}(30,54)
    \put(0,30){\includegraphics[width=30mm, clip=true, viewport= 650 475 1680 1260]{figs/D_triangle1/#2.jpg}}
    \put(1.5,3){\includegraphics[width=27mm]{figs/D_triangle1/#2.\pdfpngextension}}
    \put(0,0){\makebox(30,4){$p=#1$}}
  \end{picture}
  }
\newcommand{\troneutwoextragraphs}[2]{%
  \begin{picture}(30,57)
    \put(0,30){\includegraphics[width=30mm, clip=true, viewport= 650 475 1680 1370]{figs/D_triangle1_evensym/#2.jpg}}
    \put(1.5,3){\includegraphics[width=27mm]{figs/D_triangle1_evensym/#2.\pdfpngextension}}
    \put(0,0){\makebox(30,4){$p=#1$}}
  \end{picture}
  }
\newcommand{\trthreeuonegraphs}[2]{%
  \begin{picture}(30,56)
    \put(1,23){\includegraphics[width=28mm, clip=true, viewport= 660 340 1510 1340]{figs/D_triangle3/#2.jpg}}
    \put(3,2){\includegraphics[width=24mm]{figs/D_triangle3/#2.\pdfpngextension}}
    \put(0,0){\makebox(30,4){$p=#1$}}
  \end{picture}
  }
\newcommand{\trthreeutwographs}[2]{%
  \begin{picture}(30,56)
    \put(1,23){\includegraphics[width=28mm, clip=true, viewport= 660 340 1510 1340]{figs/D_triangle3/#2.jpg}}
    \put(3,2){\includegraphics[width=24mm]{figs/D_triangle3/#2.\pdfpngextension}}
    \put(0,0){\makebox(30,4){$p=#1$}}
  \end{picture}
  }
\newcommand{\trthreeutwoextragraphs}[2]{%
  \begin{picture}(30,59)
    \put(1,26){\includegraphics[width=28mm, clip=true, viewport= 660 340 1510 1340]{figs/D_triangle3_p_1-5/#2.jpg}}
    \put(3,2){\includegraphics[width=24mm]{figs/D_triangle3_p_1-5/#2.\pdfpngextension}}
    \put(0,0){\makebox(30,4){#1}}
  \end{picture}
  }
\newlength{\tableheadsep}
\newlength{\tablecolwidthtreq}
\begin{document}
\title[Numerical investigation of eigenvalues of the $p$-Laplace operator]{Numerical investigation of the smallest eigenvalues of the \boldmath $p$-Laplace operator on planar domains}
\author{Ji\v{r}\'{\i} Hor\'ak}
\email{jhorak@math.uni-koeln.de}
\date{\today}
\subjclass[2010]{Primary: 35J92, 49M30; Secondary: 49R05}

\begin{abstract}
  The eigenvalue problem for the $p$-Laplace operator with $p>1$ on
  planar domains with the zero Dirichlet boundary condition is
  considered. The Constrained Descent Method and the Constrained
  Mountain Pass Algorithm are used in the Sobolev space setting to
  numerically investigate the dependence of the two smallest
  eigenvalues on $p$. Computations are conducted for values of $p$
  between 1.1 and 10. Symmetry properties of the second eigenfunction
  are also examined numerically. While for the disk an odd symmetry
  about the nodal line dividing the disk in halves is maintained for
  all the considered values of $p$, for rectangles and triangles
  symmetry changes as $p$ varies. Based on the numerical evidence the
  change of symmetry in this case occurs at a certain value $p_0$
  which depends on the domain.
\end{abstract}

\maketitle

\section{Introduction}
For a bounded domain $\Omega\subset\R^N$, $N\in\N$ and a parameter
$p\in(1,\infty)$ consider the nonlinear eigenvalue problem
\begin{equation}
  \label{eq:evproblem}
  \begin{aligned}
    -\plap u &= \lambda |u|^{p-2} u &&\text{in }\Omega, \\
    u &= 0 &&\text{on }\partial\Omega
  \end{aligned}
\end{equation}
to be solved for a real function $u:\Omega\to\R$ and a parameter
$\lambda\in\R$. The operator $\plap u := \dvg \left(|\nabla
  u|^{p-2}\nabla u\right)$ is called the $p$-Laplace operator. If for
a certain $\lambda$ a nontrivial weak solution $u\in\sobspace$
of~\pref{eq:evproblem} exists, we call $\lambda$ and $u$ Dirichlet
eigenvalue and eigenfunction of the $p$-Laplace operator,
respectively. Problem \pref{eq:evproblem} is homogeneous but in
general not additive in $u$.

From \cite{Anane,Lindqvist1,Lindqvist1ad,BelloniKawohl1} and others it
is well known that there exists the smallest eigenvalue $\lambda_1$
and that it is positive, isolated and simple (i.e., the corresponding
eigenfunction $u_1$ is unique up to multiplication by a
constant). Moreover, for any eigenfunction $u$ it holds: $u$
corresponds to $\lambda_1$ if and only if it does not change its sign
on $\Omega$. In \cite{GarciaPeral} using a variational approach the
authors constructed a nondecreasing sequence of eigenvalues
accumulating at infinity. Since between $\lambda_1$ and the next
member of this sequence there are no other eigenvalues, as it was
shown in \cite{AnaneTsouli1}, we call this second smallest eigenvalue
$\lambda_2$ and a corresponding eigenfunction $u_2$. In general,
however, it is not known yet whether this sequence contains all the
eigenvalues. Nodal domains of variational eigenfunctions were studied
in \cite{DrabekRobinson}. The regularity results of \cite{DiBenedetto}
imply that any eigenfunction (perhaps redefined on a class of measure
zero) is of class $C^{1,\alpha}(\Omega)$ for some $\alpha>0$.

An early attempt at computing several eigenpairs of the $p$-Laplace
operator on a planar domain ($N=2$) numerically is due to Brown and
Reichel \cite{BrownReichel}. Under the assumption of radial symmetry
they used a shooting method for the resulting ordinary differential
equation. The first genuinely two-dimensional approach was taken by
Yao and Zhou \cite{YaoZhou1} using their local minimax method based on
a variational formulation. For a square $\Omega=\{(x_1,x_2)\ |\
x_1,x_2\in(0,2)\}$ and $p\in\{1.75, 2.5, 3.0\}$ the authors computed
approximations to seven eigenvalues and corresponding
eigenfunctions. They observed that the found eigenfunction $u_2$ has
an odd symmetry about $x_1=1$ for $p<2$ and about $x_1=x_2$ for $p>2$.

The goal of the current work is to apply the numerical variational
methods of \cite{Ho1} to compute approximations of the two smallest
eigenvalues and to visualize the corresponding eigenfunctions on a
planar domain. In particular the focus is
\begin{itemize}
\item to extend the Constrained Mountain Pass Algorithm from the
  Hilbert space setting (as described in~\cite{ChMcK1} and \cite{Ho1})
  to the Banach space $\sobspace$; to verify that this algorithm is
  suitable even for computations with $p$ ``far'' from 2;
\item to observe the behavior of the eigenpairs for a large range of
  $p$ and compare it with the known theoretical results about the
  asymptotics for $p\to 1$ and $p\to\infty$;
\item to observe changes in symmetry of $u_2$ on various domains.
\end{itemize}

In Section~\ref{sec:background} we review known results about the
variational properties of $\lambda_1$ and $\lambda_2$ and their
asymptotic behavior. The variational numerical methods applied to
compute the eigenpairs are summarized in
Sec.~\ref{sec:num_method}. The choice of a descent direction in the
Banach space $\sobspace$ is discussed in detail here, too. In
Sec.~\ref{sec:numerical_results} we present the numerical results for
several planar domains. We pay a particular attention to the
dependence of the eigenvalues on $p$ and changes of symmetry of the
second eigenfunction. Several issues concerning the application of the
numerical methods (like mesh refinement, choice of parameters, etc.)
are addressed in Sec.~\ref{sec:numerics_remarks}. Finally,
Sec.~\ref{sec:conclusion} summarizes our numerical observations and
the Appendix provides proofs of several claims used in
Sec.~\ref{sec:num_method}.

\section{Background material}\label{sec:background}

In the Introduction we mentioned the existence of the first two
eigenvalues $\lambda_1$ and $\lambda_2$. Now we review some known
results about their variational characterization based on the above
references and their asymptotic behavior for $p$ close to 1 and $p$
large.

\subsection{Variational characterization of $\lambda_1$ and $\lambda_2$}
Define two continuously Fr\'echet differentiable functionals $I,J\in
C^1\big(\sobspace,\R\big)$:
\begin{equation}
  \label{eq:IJ}
  I(u) := \int_\Omega |\nabla u|^p\ dx, \qquad
  J(u) := \int_\Omega |u|^p\ dx.
\end{equation}
Their Fr\'echet derivatives $I'(u),J'(u)$ are members of the dual
space of $\sobspace$ which we denote by $\sobdual$, where
$\frac{1}{p}+\frac{1}{q}=1$, and are given by
\begin{equation}
  \label{eq:IJprime}
  \dualpairing{I'(u)}{\phi} = p\int_\Omega |\nabla u|^{p-2}\nabla u\nabla\phi\ dx ,\qquad
  \dualpairing{J'(u)}{\phi} = p\int_\Omega |u|^{p-2} u\phi\ dx .
\end{equation}

Two observation can be made: 1.~After testing \pref{eq:evproblem} with
$\phi\in\sobspace$ and integrating by parts it becomes clear that
\pref{eq:evproblem} is the Euler-Lagrange equation $I'(u)-\lambda
J'(u)=0$ (up to the factor $p$) which all critical points of $I$ with
respect to the constraint
\begin{equation}
  \label{eq:S}
  S:=\big\{u\in\sobspace\ \big|\ J(u)=1\big\}
\end{equation}
must satisfy for some value of the Lagrange multiplier $\lambda$.

2.~If $(\lambda,u)$ is an eigenpair and we test \pref{eq:evproblem}
with $u$, we obtain the Rayleigh quotient
\begin{equation}
  \label{eq:rayleigh}
  \lambda=\frac{\int_\Omega|\nabla u|^p\ dx}{\int_\Omega|u|^p\ dx}.
\end{equation}
Since both its numerator and denominator are homogeneous of the same
degree in $u$, finding the smallest eigenvalue $\lambda$ is the same
as minimizing $I$ on $S$:
\begin{equation}
  \label{eq:lambda1}
  \lambda_1=\min_{u\in S} I(u) .
\end{equation}

A variational minimax characterization of the second eigenvalue
$\lambda_2$ based on the Krasnoselskii genus was given in
\cite{GarciaPeral}. Alternatively, since for $u_1\in S$ both $u_1$ and
$-u_1$ are local minimizers of $I$ on $S$, a mountain pass
characterization of $\lambda_2$ is also possible
\cite{CuestadeFigueiredoGossez}:
\begin{equation}
  \label{eq:lambda2}
  \lambda_2=\inf_{\gamma\in\Gamma}\max_{u\in\gamma([0,1])}I(u),
\end{equation}
where $\Gamma=\{\gamma\in C([0,1],S)\ |\ \gamma(0)=u_1,\
\gamma(1)=-u_1\}$ is the family of all paths in $S$ connecting the two
local minimizers. Hence for the numerical computations we have the
setting required by the Constrained Mountain Pass Algorithm of
\cite{Ho1}.

\subsection{Asymptotic behavior of $\lambda_1$ and $\lambda_2$ as $p\to 1$}
To make the dependence of an eigenvalue on the domain $\Omega$ and the
parameter $p$ explicit in our notation we will write
$\lambda(\Omega;p)$ if necessary (and similarly for
eigenfunctions). The main result of \cite{KawohlFridman} implies that
for $\Omega$ with a Lipschitz boundary
\begin{equation}
  \label{eq:lambda1to1}
  \lim_{p\to 1}\lambda_1(\Omega;p) = h_1(\Omega), \qquad\text{where }
  h_1(\Omega):=\min_{D\subset\Omega}\frac{\per(D)}{|D|}
\end{equation}
is called Cheeger constant, $\per(D)$ denotes the perimeter of $D$
measured with respect to $\R^N$ and $|D|$ its $N$-dimensional Lebesgue
measure. A minimizer in the definition of $h_1(\Omega)$ is called a
Cheeger set of $\Omega$. Furthermore, any convex planar domain
$\Omega$ possesses a unique Cheeger set $\cset_\Omega$ and
\begin{equation}
  \label{eq:u1to1}
  \lim_{p\to 1}u_1(\Omega;p)=\chi_{\cset_\Omega} \quad\text{in $L^1$ along a subsequence.}
\end{equation}
Here the eigenfunctions $u_1(\Omega;p)$ have been normalized to 1 in
the $L^\infty$-norm, $\chi_{\cset_\Omega}$ is the indicator function
of $\cset_\Omega$. 

A detailed description of how to find the Cheeger set $\cset_\Omega$
for a convex planar domain $\Omega$ is given in
\cite{LachandRobertKawohl}. Its main property is
\begin{equation}
  \label{eq:Cheegerset}
  \cset_\Omega=\bigcup \left\{B\subset\Omega\ \left|\ B\text{ is a ball of radius }
    \textstyle\frac{1}{h_1(\Omega)}\right.\right\}.
\end{equation}

In \cite{Parini} it was shown that for $\Omega$ with a Lipschitz
boundary it holds:
\begin{equation}
  \label{eq:lambda2to1}
  \lim_{p\to 1}\lambda_2(\Omega;p) = h_2(\Omega),
\end{equation}
where
\begin{equation}
  \label{eq:h2}
  h_2(\Omega):=\min\left\{\mu\in\R\ \Bigg|\ \exists D_1,D_2\subset\Omega,
  D_1\cap D_2=\emptyset\text{ and }
    \max_{i=1,2}\frac{\per(D_i)}{|D_i|}\leq\mu\right\}
\end{equation}
is called the second Cheeger constant and the convention
$\per(D)/|D|=\infty$ is used if $|D|=0$. Any two sets $D_1, D_2$ for
which the minimum in the definition of $h_2(\Omega)$ is achieved are
called coupled Cheeger sets of $\Omega$. For a result about the
$L^1$-convergence of the second eigenfunctions we refer to
\cite[Thm.~5.11]{Parini}.

\subsection{Asymptotic behavior of $\lambda_1$ and $\lambda_2$ as $p\to\infty$}
For a bounded domain $\Omega$ of $\R^N$ a limit problem of
\pref{eq:evproblem} as $p\to\infty$ is studied in
\cite{JuutinenLindqvistManfredi,JuutinenLindqvist} for an unknown
function $u$ and an unknown real parameter $\Lambda$ (see
\cite[Definition 2.1]{JuutinenLindqvist}). The smallest $\Lambda$ for
which this limit problem admits a nontrivial viscosity solution is
called the first $\infty$-eigenvalue and denoted $\Lambda_1$. For
$\Lambda_1$ there exists a positive viscosity solution and it holds:
\begin{equation}
  \label{eq:lambda1toinfty}
\lim_{p\to\infty}\big(\lambda_1(\Omega;p)\big)^{1/p} = \Lambda_1(\Omega),
\end{equation}
\begin{align}
  \Lambda_1(\Omega) &=\frac{1}{r_1}, \qquad\text{where } r_1:=\sup\{r>0\ |\ \exists
  \text{ an open ball } B\subset\Omega \text{ of radius } r\}, \label{eq:Lambda1}\\
  \Lambda_1(\Omega) &=\min\left\{\left.\frac{\|\nabla u\|_{L^\infty(\Omega)}}{\|u\|_{L^\infty(\Omega)}}\ \right|\ 
    u\in W^{1,\infty}_0(\Omega)\setminus\{0\}\right\}.\label{eq:Lambda1var}
\end{align}
The characterization \pref{eq:Lambda1var} is an analogy of
\pref{eq:rayleigh} and \pref{eq:lambda1}. Furthermore, for any
sequence $\{u_1(\Omega;p_i)\}_{i=1}^\infty$ with $p_i\to\infty$ and
$\|u_1(\Omega;p_i)\|_{L^{p_i}(\Omega)}=1$ there exists a subsequence
converging uniformly to a viscosity solution of the limit problem for
$\Lambda_1(\Omega)$.

The smallest $\Lambda$ for which the limit problem admits a viscosity
solution with at least two nodal domains is called the second
$\infty$-eigenvalue and denoted $\Lambda_2$. From the definition it
follows that $\Lambda_1\leq\Lambda_2$. If $\Lambda_1<\Lambda_2$, then
for $\Lambda\in(\Lambda_1,\Lambda_2)$ zero is the only solution of the
limit problem. It holds:
\begin{equation}
  \label{eq:lambda2toinfty}
  \lim_{p\to\infty}\big(\lambda_2(\Omega;p)\big)^{1/p} = \Lambda_2(\Omega),
\end{equation}
\begin{align}
  \Lambda_2(\Omega) &=\frac{1}{r_2}, \quad\text{where } r_2:=\sup\{r>0\ |\ \exists
  \text{ disjoint open balls } B_1, B_2\subset\Omega \text{ of radius } r\},  \label{eq:Lambda2} \\
  \Lambda_2(\Omega) &=\inf_{\gamma\in\Gamma}\max_{u\in\gamma([0,1])}\|\nabla u\|_{L^\infty(\Omega)},
  \label{eq:Lambda2var}
\end{align}
where $\Gamma$ is defined as in \pref{eq:lambda2}, $u_1$ is any first
$\infty$-eigenfunction and $S:=\{u\in W^{1,\infty}_0(\Omega)\ |\
\|u\|_{L^\infty(\Omega)}=1\}$. Furthermore, for any sequence
$\{u_2(\Omega;p_i)\}_{i=1}^\infty$ with $p_i\to\infty$ and
$\|u_2(\Omega;p_i)\|_{L^{p_i}(\Omega)}=1$ there exists a subsequence
converging uniformly to a viscosity solution of the limit problem for
$\Lambda_2(\Omega)$ which has at least two nodal domains.

\section{Numerical methods}\label{sec:num_method}
\begin{figure}[h]
  \centering
  \setlength{\unitlength}{1mm}
  \begin{picture}(108,56)(0,14)
    \put(-1,10){\includegraphics{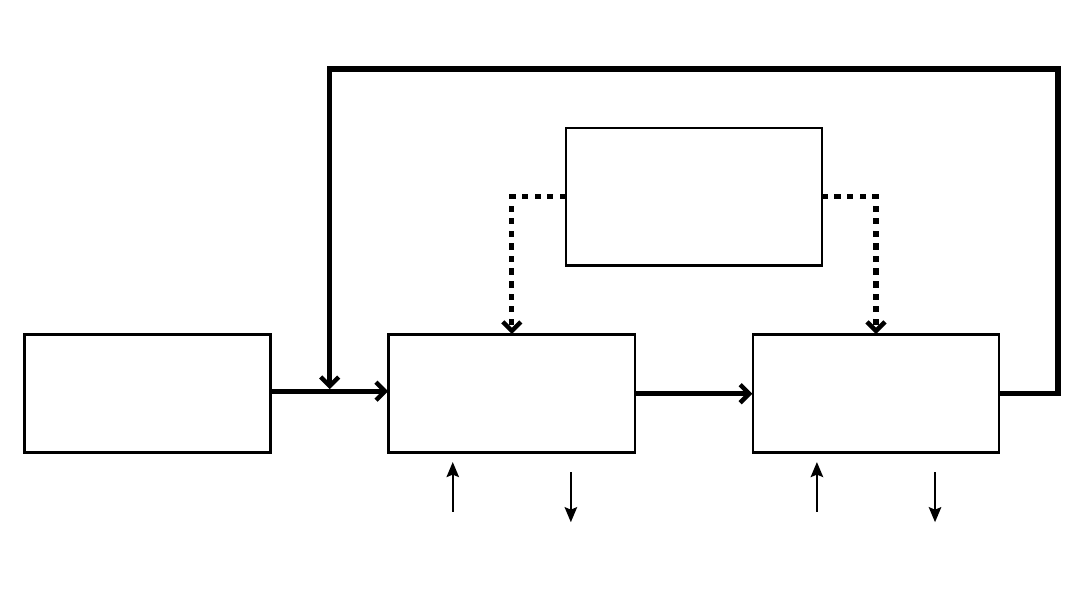}}
    \put(2.5,29.5){\parbox{23mm}{\centering\small Triangulation of~$\Omega$, FEM}}
    \put(39.5,29.5){\parbox{23mm}{\centering\small CDM (iterative)}}
    \put(76.5,29.5){\parbox{23mm}{\centering\small CMPA (iterative)}}
    \put(57.6,49){\parbox{24mm}{\centering\small Aug.\ Lagrangian for $(-\Delta_p)^{-1}$ (iterative)}}
    \put(34,14.5){\parbox{23mm}{\centering\small $e_0$}}
    \put(45.5,13.5){\parbox{23mm}{\centering\small $(\lambda_1,u_1)$}}
    \put(68,31.5){\small $u_1$}
    \put(71,14.5){\parbox{23mm}{\centering\small $e_\mathrm{M}$}}
    \put(85,13.5){\parbox{23mm}{\centering\small $(\lambda_2,u_2),\ldots$}}
    \put(28, 64.5){\parbox{83mm}{\centering\small loop over a range of values of $p\in(1,\infty)$}}
  \end{picture}
  \caption{Flowchart of the numerical computations.}
  \label{fig:flowchart}
\end{figure}
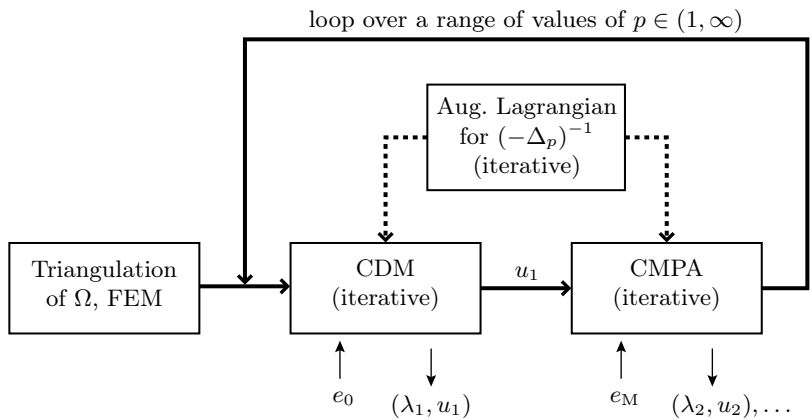

An overview of the numerical methods used to compute approximations of
the first and the second Dirichlet eigenpair of the $p$-Laplace
operator is given in Fig.~\ref{fig:flowchart}. In this section we will
describe these methods. Our goal is to find $u_1$ as a minimizer of
$I$ on $S$ according to \pref{eq:lambda1} and $u_2$ as a mountain pass
point of $I$ on $S$ according to \pref{eq:lambda2}. We first
discretize the planar domain $\Omega$ using a mesh of triangles and
apply the finite element method to approximate $\sobspace$ by a finite
dimensional subspace. Then we fix $p\in(1,\infty)$ and use a variant
of the \emph{Constrained Steepest Descent Method} (CDM) to find the
first eigenpair, and the \emph{Constrained Mountain Pass Algorithm}
(CMPA) to find the second eigenpair. We implement both methods based
on \cite{Ho1}. There are, however, several important issues arising
from the fact that we work in a Banach space and not a Hilbert space
as in \cite{Ho1}. How to deal with these issues will also be explained
in this section. For the computation of the descent direction the
Augmented Lagrangian Method of \cite{GlowinskiMarroco} is applied.

\subsection{Finite element method}\label{sec:fem}
A finite element approximation of the $p$-Laplacian was studied in
\cite{BarrettLiu}. We adopt this approach for our computations. The
planar domain $\Omega$ is approximated by a polygonal domain
$\Omega^h$ which is partitioned into a finite number of triangles of
diameter at most $h$. Let $\{a_i\}_{i=1}^k$ be the set of those
triangle vertices which lie in the interior of $\Omega^h$. Functions
$\{\phi_i\}_{i=1}^k$ forming a basis of the $k$-dimensional subspace
$V^h_0$ of $W^{1,p}_0(\Omega^h)$ are chosen linear on each triangle
with $\phi_i(a_j)=\delta_{ij}$, where $\delta_{ij}$ is the Kronecker
delta, and zero on $\partial\Omega^h$. The space $V^h_0$ is our finite
element approximation of the Sobolev space $\sobspace$.

In \cite{BarrettLiu} a detailed description of this method was given
for the boundary value problem
\begin{equation}
  \label{eq:basicproblem}
  \begin{aligned}
    -\plap u &= f &&\text{in }\Omega, \\
    u &= 0 &&\text{on }\partial\Omega
  \end{aligned}
\end{equation}
with the right-hand side $f\in L^2(\Omega)$. Since it is a
straightforward task to adapt it to our problem \pref{eq:evproblem}
with $\lambda |u|^{p-2}u$ on the right-hand side we will not show the
details here.

We will however mention one additional technical detail involved. The
evaluation of the functionals given in \pref{eq:IJ} and
\pref{eq:IJprime} for functions from $V^h_0$ amounts to adding up the
contributions of the individual triangles that make up $\Omega^h$. For
example, for $I(u)$ with $u\in V^h_0$ one merely needs to integrate a
constant on every triangle. The situation is different for $J(u)$. Let
$T\subset\Omega^h$ be a triangle with area $|T|$ and vertices $A$,
$B$, and $C$ and let $u$ be a linear function on $T$ with values
$u_A$, $u_B$, and $u_C$ at these vertices, respectively. If these
values are mutually different, then the following formula holds:
\begin{equation}
  \label{eq:Jontriangle}
  \int_T |u|^p\,dx = \frac{2\,|T|}{(p+1)(p+2)(u_C-u_A)}\left(
    \frac{|u_C|^{p+2}-|u_B|^{p+2}}{u_C-u_B} - \frac{|u_A|^{p+2}-|u_B|^{p+2}}{u_A-u_B}
  \right).
\end{equation}
By inspecting this formula we see that great care must be taken when
implementing it to avoid numerical cancellations. This is crucial for
the success of our method. A similar situation occurs when evaluating
$\dualpairing{J'(u)}{\phi}$ (or the right-hand side of equation in
\pref{eq:evproblem} in the weak formulation).

\subsection{Direction of descent}\label{sec:descent}
An important ingredient of the variational numerical methods CDM and
CMPA is finding a descent direction of the functional $I$ on the
constraint set $S$. How this is accomplished in the Hilbert space
setting was shown in~\cite{Ho1}: Let $\nabla I(u)$ be the Riesz
representation of $I'(u)$ (i.e., the gradient) and $P_u$ the
orthogonal projection on the tangent space of $S$ at $u\in S$. Then
\begin{equation}\label{eq:orthogonalprojection}
  w_u=-P_u\nabla I(u),\qquad u\in S
\end{equation}
gives the steepest descent direction of $I$ at $u$ with respect to
$S$.

Because of the lack of orthogonality in the Banach space
$\sobspace$ we need to take a different approach. Let
\begin{equation}
  \label{eq:tgspace}
  T_uS:=\left\{\left.v\in \sobspace\ \right|\ \dualpairing{J'(u)}{v}=0\right\},\qquad
  \|v\|:=\left(\int_\Omega|\nabla v|^p\,dx\right)^{1/p}
\end{equation}
denote the tangent space of $S$ at $u\in S$ and the norm of $v\in
\sobspace$, respectively. The problem of finding the steepest descent
direction of $I$ with respect to $S$ can be written as follows: for a
given $u\in S$ which is not a critical point of $I$ with respect to
$S$
\begin{equation}
  \label{eq:steepestdescentBanach}
  \text{minimize } \dualpairing{I'(u)}{w} \quad \text{ subject to } \quad w\in\{v\in T_uS\ |\ \|v\|=1\}.
\end{equation}
It has a unique solution as Lemma~\ref{lemma1} in the Appendix shows.
The Euler-Lagrange equation that this solution must satisfy can be written in the form
\begin{equation}
  \label{eq:elg_steepest_descent}
  -\Delta_p w = \beta\left( -\Delta_p u - \alpha |u|^{p-2}u \right),
\end{equation}
where $\alpha,\beta\in\R$ are unknown. The coefficient $\beta$ comes
from the requirement $\|w\|=1$. After testing
\pref{eq:elg_steepest_descent} by $w$ it can be seen that $\beta<0$
since the minimum in \pref{eq:steepestdescentBanach} is negative. For
$p\neq 2$ finding the right $\alpha$ is not an easy problem.

We will try to find a different convenient descent direction instead,
not necessarily the steepest one. A simple calculation shows that
under no constraints the steepest descent direction of $I$ at $u\in B$
is given by $-u$. For $u\in S$ we consider the point $w_u\in T_u S$
closest to $-u$, i.e., the unique solution of the minimization problem
\begin{equation}
  \label{eq:descentBanach}
  \text{minimize } \|w+u\| \quad \text{ subject to } \quad w\in T_uS.
\end{equation}
The minimizer must satisfy the Euler-Lagrange equation
\begin{equation}
  \label{eq:elg_descent}
  -\Delta_p(w+u)=\alpha |u|^{p-2}u
\end{equation}
for some $\alpha\in\R$. Unlike \pref{eq:elg_steepest_descent}, this
equation can be solved easily for $w$:
\begin{equation}
  \label{eq:descent_direction}
  w_u=-u+\frac{1}{\int_\Omega |u|^{p-2}u \,v_u\,dx}\ v_u,\qquad
  \text{where }v_u:=(-\Delta_p)^{-1}\left(|u|^{p-2}u\right).
\end{equation}
The operator $(-\Delta_p)^{-1}$ is discussed later in
Sec.~\ref{sec:inverse}.  Lemma \ref{lemma2} in the Appendix shows
that $w_u$ is, indeed, a descent direction of $I$ with respect to
$S$. This descent direction is used in our implementation of the
variational numerical methods CDM and CMPA.

\begin{plapnum_remark}
  1.~Observe that if $-\Delta_p$ were linear, equations
  \pref{eq:elg_steepest_descent} and \pref{eq:elg_descent} would
  coincide (after setting $\beta=-1$). Hence in case $p=2$ they yield
  the same descent direction (which is the one given by
  \pref{eq:orthogonalprojection}).

  2.~With $\beta=-1$ in \pref{eq:elg_steepest_descent}, both equations
  \pref{eq:elg_steepest_descent} and \pref{eq:elg_descent} yield a
  zero solution if and only if $u$ is a critical point of $I$ with
  respect to $S$.
\end{plapnum_remark}

\subsection{Inverse of the $p$-Laplace operator}\label{sec:inverse}
A classical result (see, e.g., \cite[Theorem 1.3]{Struwe}) says that
for any $f\in\sobdual$, the dual of $\sobspace$ with
$\frac{1}{p}+\frac{1}{q}=1$, the problem
\begin{equation}
  \label{eq:constant_rhs}
  \begin{aligned}
    -\plap u &= f &&\text{in }\Omega, \\
    u &= 0 &&\text{on }\partial\Omega
  \end{aligned}
\end{equation}
has a unique weak solution in $\sobspace$. This means that the operator 
\begin{equation}
  \label{eq:operator}
  -\Delta_p: \sobspace\to\sobdual \quad\text{given by}\quad
  \dualpairing{-\Delta_p u}{v}=\int_\Omega|\nabla u|^{p-2}\nabla u\nabla v\,dx
\end{equation}
is invertible. We denote its inverse by $(-\Delta_p)^{-1}$.

In order to use the descent direction given by
\pref{eq:descent_direction} we need to compute $v_u$ first, i.e., we
need to solve problem \pref{eq:constant_rhs} numerically. For that we
apply the Augmented Lagrangian Method of \cite{GlowinskiMarroco}. Here
we give a brief description of this method. Let $V^h_0$ be again the
subspace of continuous functions of $W^{1,p}_0(\Omega^h)$ which are
linear on every triangle of a triangulation of $\Omega^h$, $D^h$ the
space of functions with values in $\R^2$ defined on $\Omega^h$ which
are constant on each triangle, and $r>0$ a parameter. For the
Augmented Lagrangian
\begin{equation}
  \label{eq:augmented_lagrangian}
  \mathcal{L}_r(v,t,\mu)=\frac{1}{p}\int_\Omega |t|^p\,dx - \dualpairing{f}{v}+
  \frac{r}{2}\int_\Omega|\nabla v-t|^2\,dx + \int_\Omega\mu\cdot(\nabla v-t)\,dx,
\end{equation}
where $v\in V^h_0$ and $t,\mu\in D^h$, a saddle point $(u,s,\eta)$ is
searched for such that
\begin{equation}
  \label{eq:al_saddle_point}
  \mathcal{L}_r(u,s,\mu)\leq\mathcal{L}_r(u,s,\eta)\leq\mathcal{L}_r(v,t,\eta)
  \qquad \forall (v,t,\mu)\in V^h_0\times D^h\times D^h.
\end{equation}
A sequence $(u^{(n)},s^{(n)},\eta^{(n)})$ approximating $(u,s,\eta)$ is
constructed as follows: choose $\big(s^{(0)},\eta^{(1)}\big)\in D^h\times D^h$ and for
$n\in\N$ solve
\begin{gather}
  \begin{aligned}
    -r\Delta u^{(n)} &= f+\eta^{(n)}\cdot\nabla - r\,s^{(n-1)}\cdot\nabla &&\text{in }\Omega, \\
    u^{(n)} &=0 &&\text{on }\partial\Omega, \\
  \end{aligned} \label{al_step1}\\
  \big|s^{(n)}\big|^{p-2}s^{(n)}+r\,s^{(n)} = r\nabla u^{(n)}+\eta^{(n)}, \label{al_step2}\\
  \eta^{(n+1)}=\eta^{(n)}+r\big(\nabla u^{(n)}-s^{(n)}\big).\label{al_step3}
\end{gather}

For given $s^{(n-1)}$ and $\eta^{(n)}$ the boundary value problem
\pref{al_step1} can be solved for $u^{(n)}$. For this one just needs some
standard algorithm for finding the inverse of the Laplace operator
with Dirichlet boundary conditions. The equation in \pref{al_step1} is
understood in the weak sense: for example the term $\eta^{(n)}\cdot\nabla$
is evaluated as $\int_\Omega\eta^{(n)}\cdot\nabla\phi\,dx$ for a test
function $\phi\in V^h_0$.

Next, equation \pref{al_step2} is used to find $s^{(n)}$. The $\R^2$-norm
of $s^{(n)}$ must satisfy
\begin{equation}
  \label{eq:al_step2a}
  \big|s^{(n)}\big|^{p-1}+r\,\big|s^{(n)}\big|=\big|r\,\nabla u^{(n)}+\eta^{(n)}\big|,
\end{equation}
which on each triangle is just a scalar nonlinear equation with one
unknown. For each triangle it can be solved, e.g., by Newton's
method. After $|s^{(n)}|$ has been obtained, $s^{(n)}$ can be computed
immediately from \pref{al_step2}.

At the end $\eta$ is updated according to \pref{al_step3} and a new
iteration step can be started.

The convergence of this method was studied in
\cite{GlowinskiMarroco}. We use the norm of $\nabla u^{(n)}-s^{(n)}$ to
measure the convergence.

\subsection{Constrained Descent Method}\label{sec:CDM}
The Constrained Descent Method (CDM) is applied to find the first
eigenpair of the $p$-Laplace operator: $u_1$ is found as the minimizer
of $I$ with respect to $S$, $\lambda_1=I(u_1)$. As mentioned above, it
differs from the Constrained Steepest Descent Methods of \cite{Ho1} in
the way the descent direction is chosen.

The method solves numerically the following initial value problem:
\begin{equation}
  \label{eq:cdm_ivp}
  \frac{d}{d t}u(t)=w_{u(t)},\qquad u(0)=e_0\in S,
\end{equation}
where $w_u$ is given by \pref{eq:descent_direction} for $u\in
S$. Proposition~\ref{proposition1} in the Appendix states that this
problem has a unique solution $u(t)\in S$ for $t\in(0,\infty)$ and
that $u(t)$ gets arbitrarily close to a critical point of $I$ with
respect to $S$ as $t\to\infty$.

After choosing the starting point $e_0\in S$ and setting
$u^{(0)}:=e_0$ the initial value problem is solved by repeating the
following two steps: First (Euler's step), given $u^{(n-1)}$ find
$\bar u^{(n)}=u^{(n-1)}+\Delta t^{(n)}\,w_{u^{(n-1)}}$ with some small
value $\Delta t^{(n)}>0$. Second (scaling), define $u^{(n)}=c\bar
u^{(n)}$, where the coefficient $c\in\R$ is chosen such that
$u^{(n)}\in S$. In case $I\big(u^{(n)}\big) > I\big(u^{(n-1)}\big)$,
halve the step $\Delta t^{(n)}$ and compute $\bar u^{(n)}$ and
$u^{(n)}$ again. If this halving has to be repeated and $\Delta
t^{(n)}$ becomes very small (smaller than a prescribed threshold
value), stop the algorithm. The norm of the descent direction
$\|w_{u^{(n-1)}}\|$ is used to measure convergence of $u^{(n)}$ to an
eigenfunction $u$. When computing $w_u$ according to
\pref{eq:descent_direction} the integral $\nu:=\int_\Omega |u|^{p-2}u
\,v_u\,dx$ has to be evaluated. If $\|w_{u^{(n-1)}}\|$ is small, then
$\big(1/\nu^{(n-1)}\big)^{p-1}$ approximates the eigenvalue.

We note that at every step of CDM the Augmented Lagrangian Method of
Sec.~\ref{sec:inverse} has to be applied to compute the descent
direction $w_{u^{(n-1)}}$.

\subsection{Constrained Mountain Pass Algorithm}\label{sec:CMPA}
Suppose that an approximation of the first eigenvalue $\lambda_1$ and
eigenfunction $u_1$ of the $p$-Laplace operator have been
computed. Constrained Mountain Pass Algorithm (CMPA) is applied to
find the second eigenpair: $u_2$ is found as a mountain pass point of
$I$ on $S$ lying ``between'' the two local minimizers $u_1$ and
$-u_1$, $\lambda_2=I(u_2)$. Again, it differs from CMPA described in
detail in \cite{Ho1} in the choice of the descent direction.

We give a short summary here based on the original description of the
Mountain Pass Algorithm by Choi and McKenna \cite{ChMcK1}: Take a
discretized path $\{z_j\}_{j=0}^{P}\subset S$ connecting $z_0:=u_1$
with $z_P:=-u_1$. After finding the path point $z_m=:z^\mathrm{max}$
at which $I$ is maximal along the path, move this point a small
distance in the tangent space to~$S$ at $z^\mathrm{max}$ in the
descent direction $w_{z^\mathrm{max}}$ and then scale it (as in CDM)
to come back to $S$. Thus the path has been deformed on~$S$ and the
maximum of $I$ lowered. Repeat this deforming of the path until the
maximum along the path cannot be lowered anymore: a mountain pass
point of $I$ with respect to~$S$ has been reached.

To construct the initial path connecting $u_1$ and $-u_1$ in $S$ we
choose an intermediate point $e_\mathrm{M}\in S\setminus\{\pm
u_1\}$, set $k:=[P/2]$ and define:
\begin{align}
  &\bar z_j:=u_1+\frac{j}{k}(e_\mathrm{M}-u_1) &&\text{for }j\in\{0,\ldots,k\}, \notag\\
  &\bar z_j:=e_\mathrm{M}+\frac{j-k}{P-k}(-u_1-e_\mathrm{M}) &&\text{for }j\in\{k,\ldots,P\}, \notag\\
  &z_j:=c_j\bar z_j\in S\text{ (scaling to $S$ as in Sec.~\ref{sec:CDM})}
  &&\text{for }j\in\{0,\ldots,P\}. \notag
\end{align}
Connecting $u_1$ and $-u_1$ by a line segment without the intermediate
point $e_\mathrm{M}$ would not work. Such a line segment passes
through $0$ and hence cannot be scaled to get to $S$.

Finally, as in CDM, $\|w_{z^\mathrm{max}}\|$ is used to measure
convergence to an eigenfunction $u$. The corresponding eigenvalue
$\lambda$ is computed as in CDM, too.  At every step of CMPA the
Augmented Lagrangian Method of Sec.~\ref{sec:inverse} has to be
applied to compute the descent direction $w_{z^\mathrm{max}}$.

\section{Numerical results}\label{sec:numerical_results}
In this section numerical results will be given for the following
planar domains: the unit disk, the square with side length 2, the
rectangle with sides 2 and 7/4, the isosceles triangle with base and
height 1, the isosceles triangle with base 1 and height 3/4, and the
equilateral triangle with side 1. Unless explicitly stated otherwise
the computed eigenfunctions will be plotted as a surface over the
domain with heights given by the function values and as a contour plot
of these values (like, e.g., in Fig.~\ref{fig:disk_u1}). In order to
better compare the shapes the eigenfunctions in these figures have
been scaled to have the same maximum value. We do not explicitly
differentiate between two eigenfunctions $u$ and $\tilde u$ if $\tilde
u(x)=cu(Tx)$, where $c\in\R$ is a scaling coefficient and
$T:\Omega\to\Omega$ is some symmetry transformation of $\Omega$ (e.g.,
for a square a rotation by $\pi/2$ about the center of the square).

\subsection{Unit Disk}
Let
\begin{displaymath}
  \Omega=\left\{(x_1,x_2)\in\R^2\,\big|\, x_1^2+x_2^2<1\right\}.
\end{displaymath}
Before presenting the numerical results we make a remark about the
radially symmetric case. It is known that the first eigenfunction for
the disk is radially symmetric. One important question about the
second eigenfunction for the disk has been whether it is radially
symmetric, too. In \cite{Parini,BenDrabGirg} the authors proved that for
$p$ close to 1 the answer is no. The eigenvalue problem
\pref{eq:evproblem} under the assumption of radial symmetry $u=u(r)$,
$r\in(0,1)$ becomes
\begin{equation}
  \label{eq:evproblem_radsym}
  \begin{gathered}
    -\left(r|u'|^{p-2}u'\right)^\prime = \lambda r |u|^{p-2}u, \\
    u'(0)=0,\quad u(1)=0.
  \end{gathered}
\end{equation}
This and a related problem are treated, for example, in
\cite{BrownReichel} and \cite{BenDrabGirg}, where numerical approaches
play an important role. For our numerical investigation we adapt the
\emph{genuine 2D method} of Sec.~\ref{sec:num_method} in the
following ways:
\begin{itemize}
\item all integrals are one-dimensional,
\item the weight $r$ is introduced,
\item the natural boundary condition is implemented at $r=0$ (the zero boundary condition stays at $r=1$).
\end{itemize}
Since these modifications are rather elementary, we will not describe
them in more detail. We will refer to this method as \emph{radial 2D
  method}.

For the computations carried out by the genuine 2D method the domain
$\Omega$ was approximated by a polygon and discretized using 68,608
triangles. For the computations carried out by the radial 2D method
the interval $(0,1)$ was divided into 1,000 subintervals of the same
length.

\begin{table}[h]
  \centering\vspace{1ex}
  \begin{tabular}[t]{c|c|c|c|c|cc|c|c|c|}
    \cline{2-5}\cline{8-10}
    (a) & $p$ & $\lambda_1$ & $\lambda_2$ & \multicolumn{1}{c|}{$\lambda_2^\mathrm{rad}$} & \qquad\qquad & (b) & $p$ & $\lambda_1$ & $\lambda_2^\mathrm{rad}$ \\\cline{2-5}\cline{8-10}
    \multicolumn{10}{c}{\vspace{\tableheadsep}}\\\cline{2-5}\cline{8-10}
    & 1.1 & 2.5690 & 4.2008 & 5.6809 &&& 1.1 & 2.5688 & 5.6762 \\\cline{2-5}\cline{8-10}
    & 1.2 & 2.9654 & 5.0707 & 7.2277 &&& 1.2 & 2.9653 & 7.2251 \\\cline{2-5}\cline{8-10}
    & 1.3 & 3.3263 & 5.9604 & 8.9302 &&& 1.3 & 3.3260 & 8.9279 \\\cline{2-5}\cline{8-10}
    & 1.4 & 3.6740 & 6.9072 & 10.861 &&& 1.4 & 3.6739 & 10.858 \\\cline{2-5}\cline{8-10}
    & 1.5 & 4.0179 & 7.9310 & \multicolumn{1}{c}{} &&& 1.5 & 4.0177 & 13.073 \\\cline{2-4}\cline{8-10}
    & 1.6 & 4.3623 & 9.0465 & \multicolumn{1}{c}{} &&& 1.6 & 4.3621 & 15.626 \\\cline{2-4}\cline{8-10}
    & 1.7 & 4.7097 & 10.266 & \multicolumn{1}{c}{} &&& 1.7 & 4.7095 & 18.574 \\\cline{2-4}\cline{8-10}
    & 1.8 & 5.0618 & 11.604 & \multicolumn{1}{c}{} &&& 1.8 & 5.0616 & 21.982 \\\cline{2-4}\cline{8-10}
    & 1.9 & 5.4194 & 13.072 & \multicolumn{1}{c}{} &&& 1.9 & 5.4192 & 25.921 \\\cline{2-4}\cline{8-10}
    & 2.0 & 5.7834 & 14.683 & \multicolumn{1}{c}{} &&& 2.0 & 5.7831 & 30.471 \\\cline{2-4}\cline{8-10}
    & 2.1 & 6.1542 & 16.452 & \multicolumn{1}{c}{} &&& 2.1 & 6.1539 & 35.725 \\\cline{2-4}\cline{8-10}
    & 2.2 & 6.5320 & 18.395 & \multicolumn{1}{c}{} &&& 2.2 & 6.5317 & 41.788 \\\cline{2-4}\cline{8-10}
    & 2.3 & 6.9173 & 20.527 & \multicolumn{1}{c}{} &&& 2.3 & 6.9169 & 48.780 \\\cline{2-4}\cline{8-10}
    & 2.4 & 7.3102 & 22.866 & \multicolumn{1}{c}{} &&& 2.4 & 7.3097 & 56.836 \\\cline{2-4}\cline{8-10}
    & 2.5 & 7.7107 & 25.432 & \multicolumn{1}{c}{} &&& 2.5 & 7.7102 & 66.112 \\\cline{2-4}\cline{8-10}
    & 3.0 & 9.8323 & 42.460 & \multicolumn{1}{c}{} &&& 3.0 & 9.8314 & 137.93 \\\cline{2-4}\cline{8-10}
    & 4.0 & 14.683 & 110.71 & \multicolumn{1}{c}{} &&& 4.0 & 14.681 & 559.02 \\\cline{2-4}\cline{8-10}
    & 5.0 & 20.351 & 273.00 & \multicolumn{1}{c}{} &&& 5.0 & 20.347 & 2,132.7 \\\cline{2-4}\cline{8-10}
    & 6.0 & 26.832 & 649.47 & \multicolumn{1}{c}{} &&& 6.0 & 26.823 & 7,822.6 \\\cline{2-4}\cline{8-10}
    & 8.0 & 42.210 & 3,430.1 & \multicolumn{1}{c}{} &&& 8.0 & 42.182 & 97,462 \\\cline{2-4}\cline{8-10}
    & 10.0 & 60.784 & 17,071 & \multicolumn{1}{c}{} &&& 10.0 & 60.715 & $1.1359\cdot 10^6$ \\\cline{2-4}\cline{8-10}
  \end{tabular}\vspace{2ex}
  \caption{Eigenvalues for the disk computed numerically by:
    (a) the genuine 2D method, (b) the radial 2D method.}
  \label{tab:disk_lambda}
\end{table}

\begin{figure}[p]
  \centering
  \setlength{\unitlength}{1mm}
  \begin{picture}(150,57)
    \put(0,0){\dsuonegraphs{1.1}{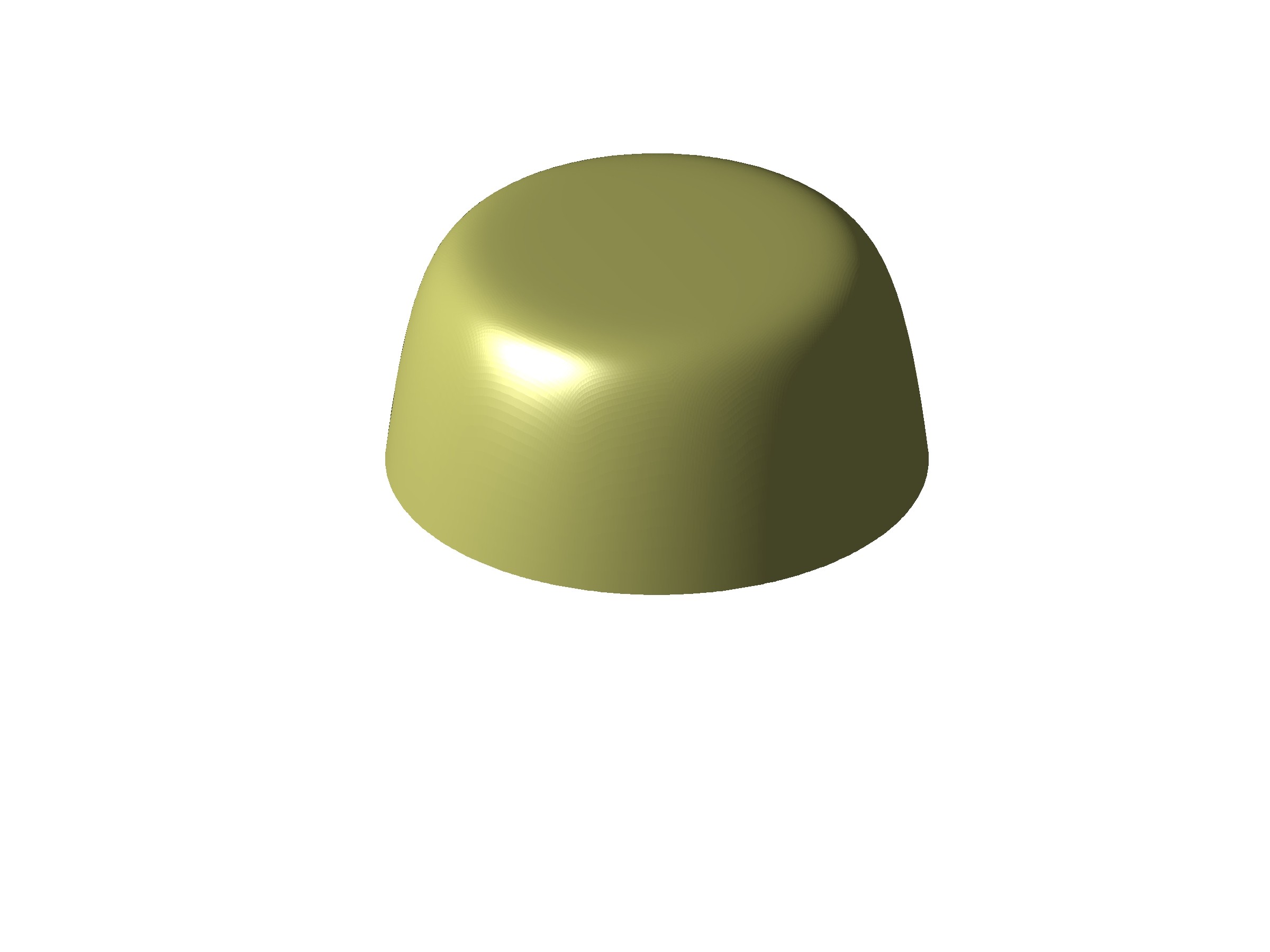}}
    \put(30,0){\dsuonegraphs{1.4}{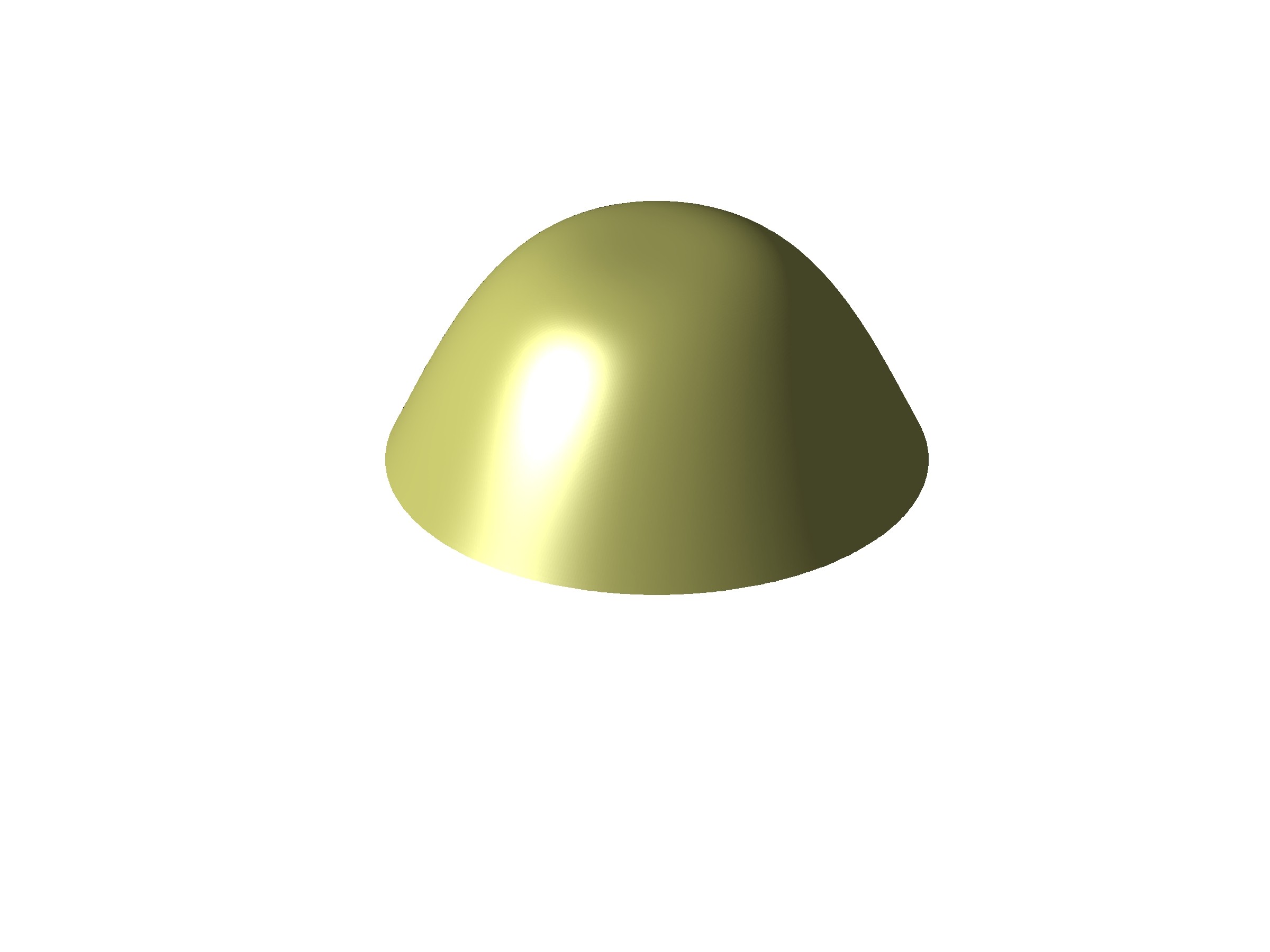}}
    \put(60,0){\dsuonegraphs{2.0}{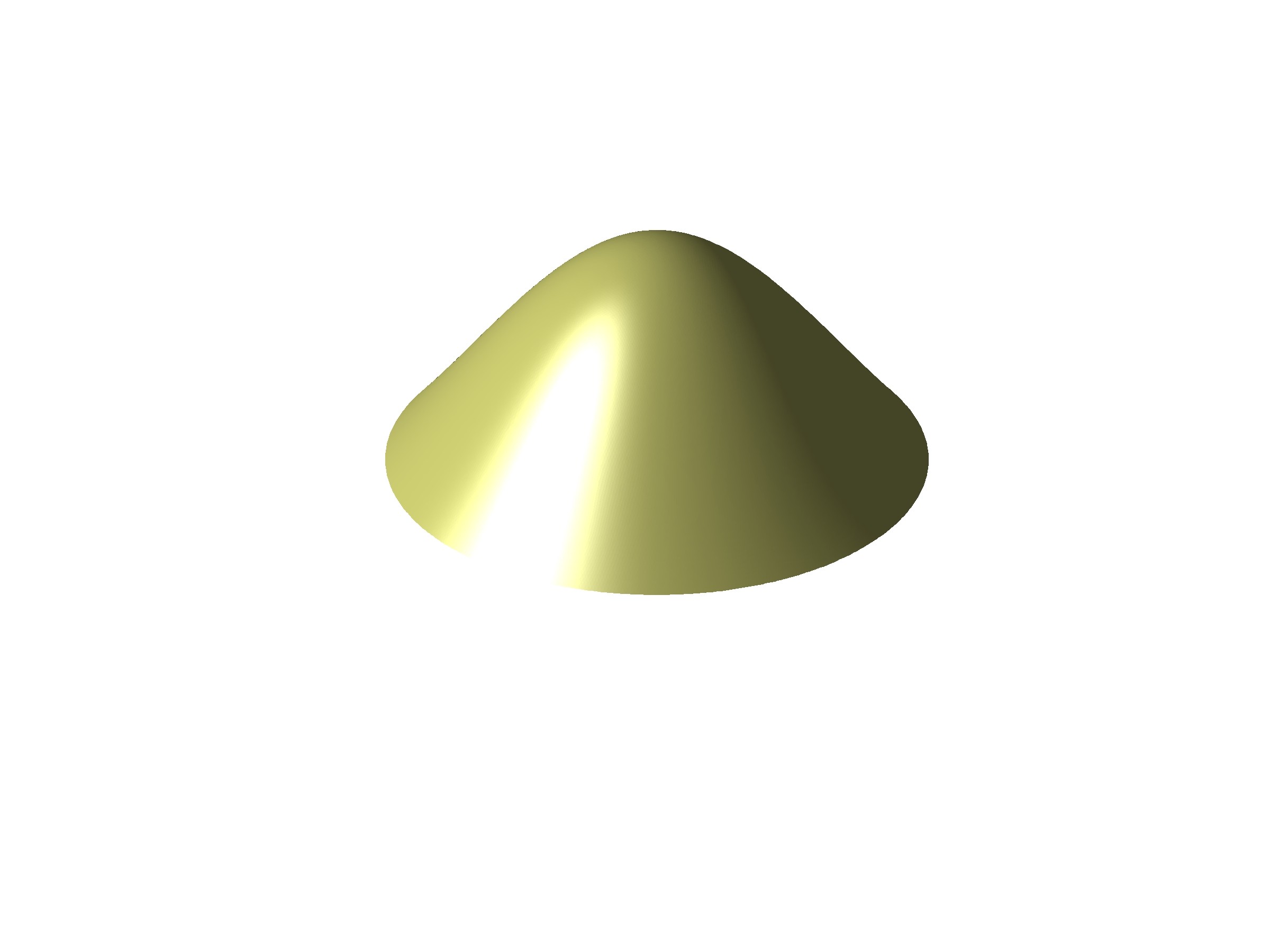}}
    \put(90,0){\dsuonegraphs{2.5}{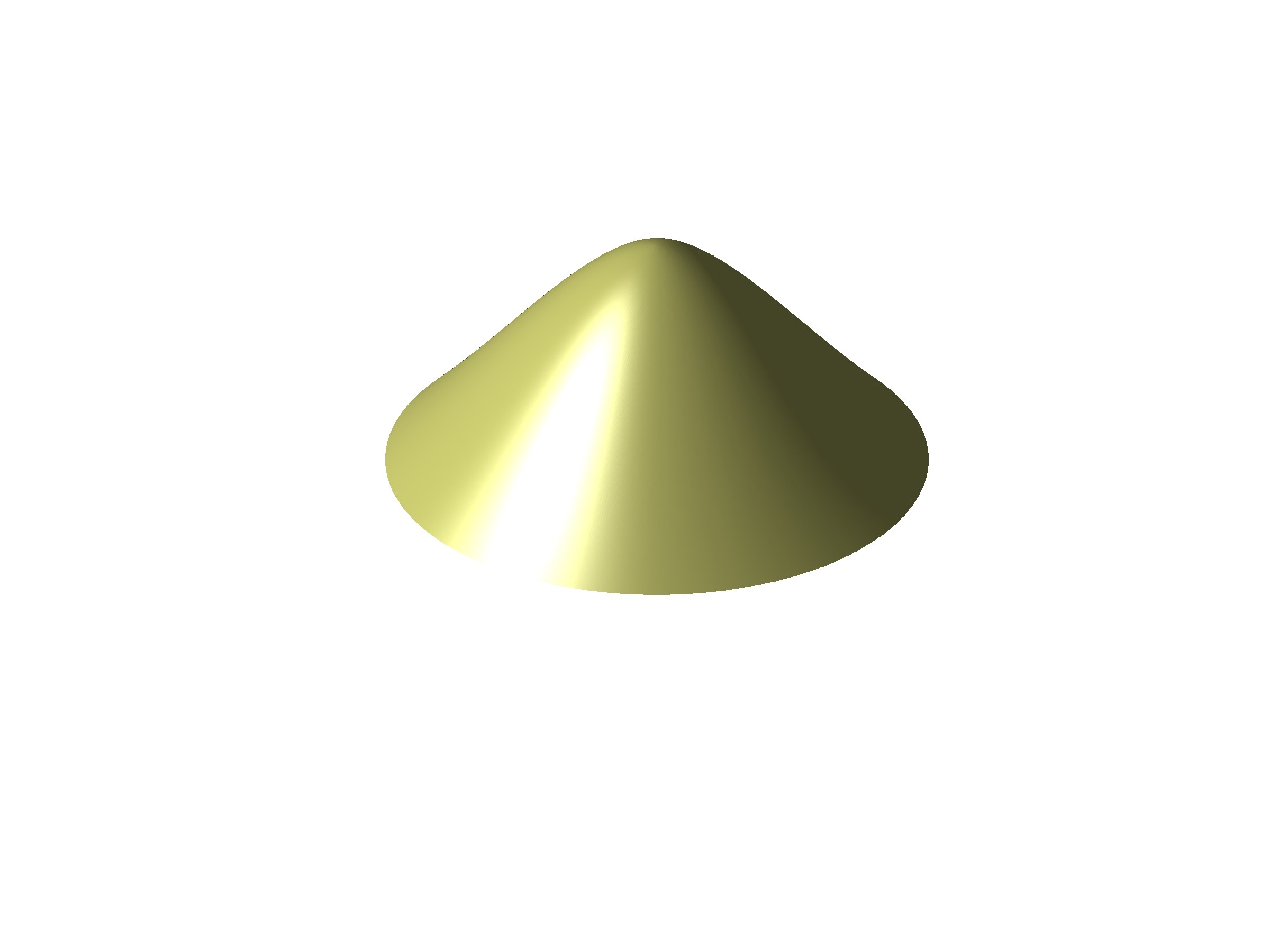}}
    \put(120,0){\dsuonegraphs{8.0}{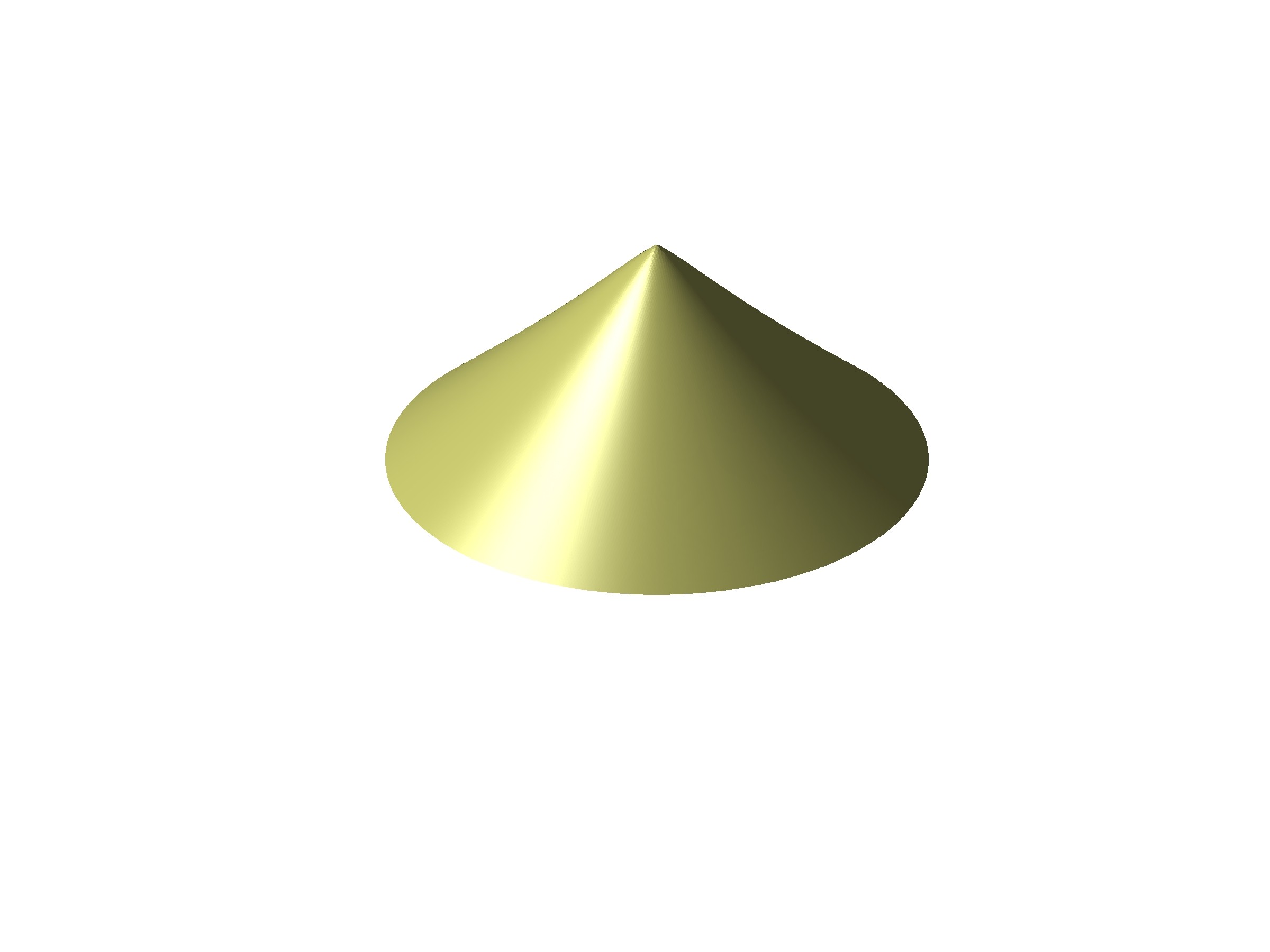}}
  \end{picture}
  \caption{The numerically computed first eigenfunction $u_1$ for the disk.}
  \label{fig:disk_u1}
\end{figure}

\begin{figure}[p]
  \centering
  \setlength{\unitlength}{1mm}
  \begin{picture}(150,62)
    \put(0,0){\dsutwographs{1.1}{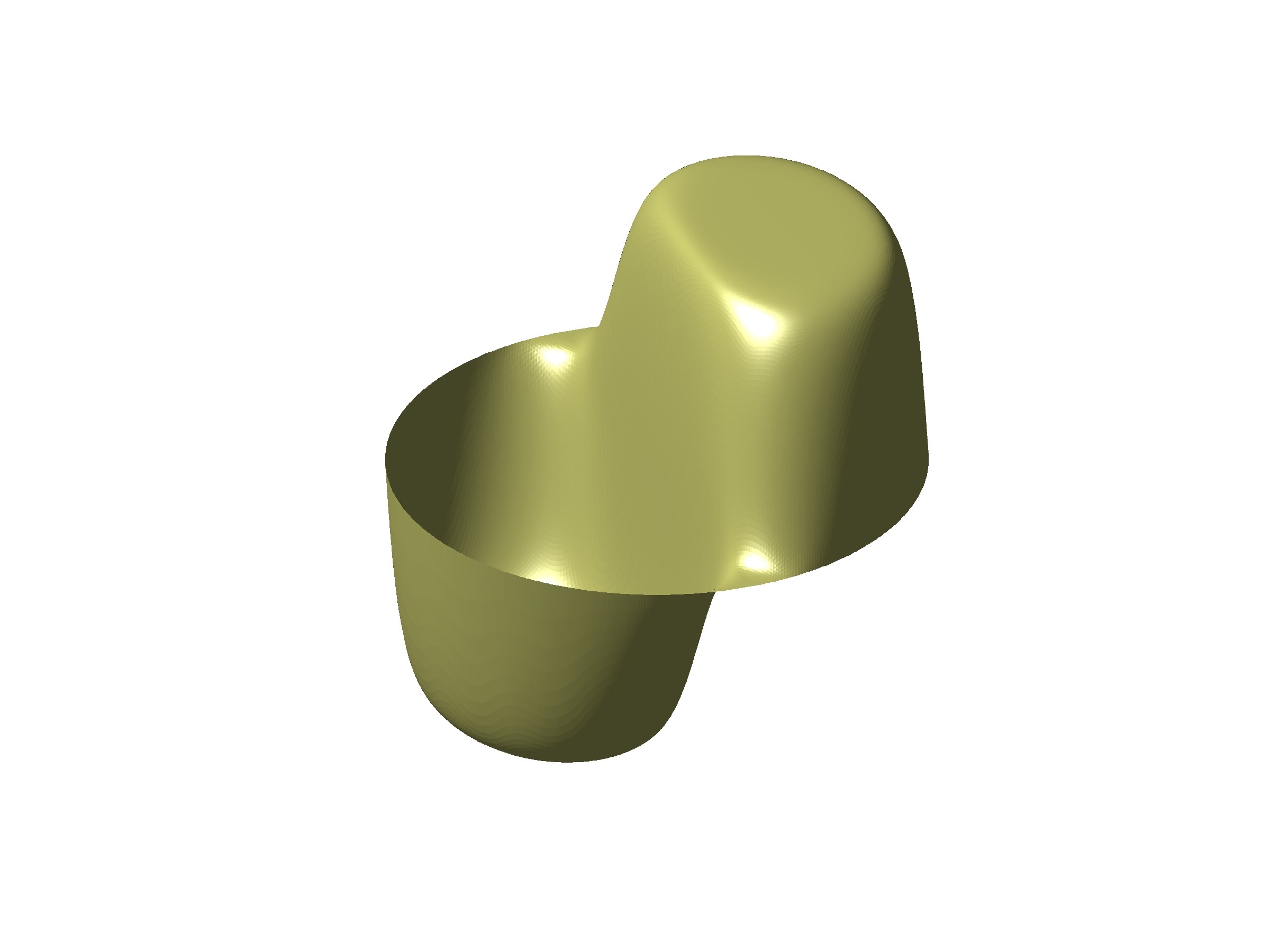}}
    \put(30,0){\dsutwographs{1.4}{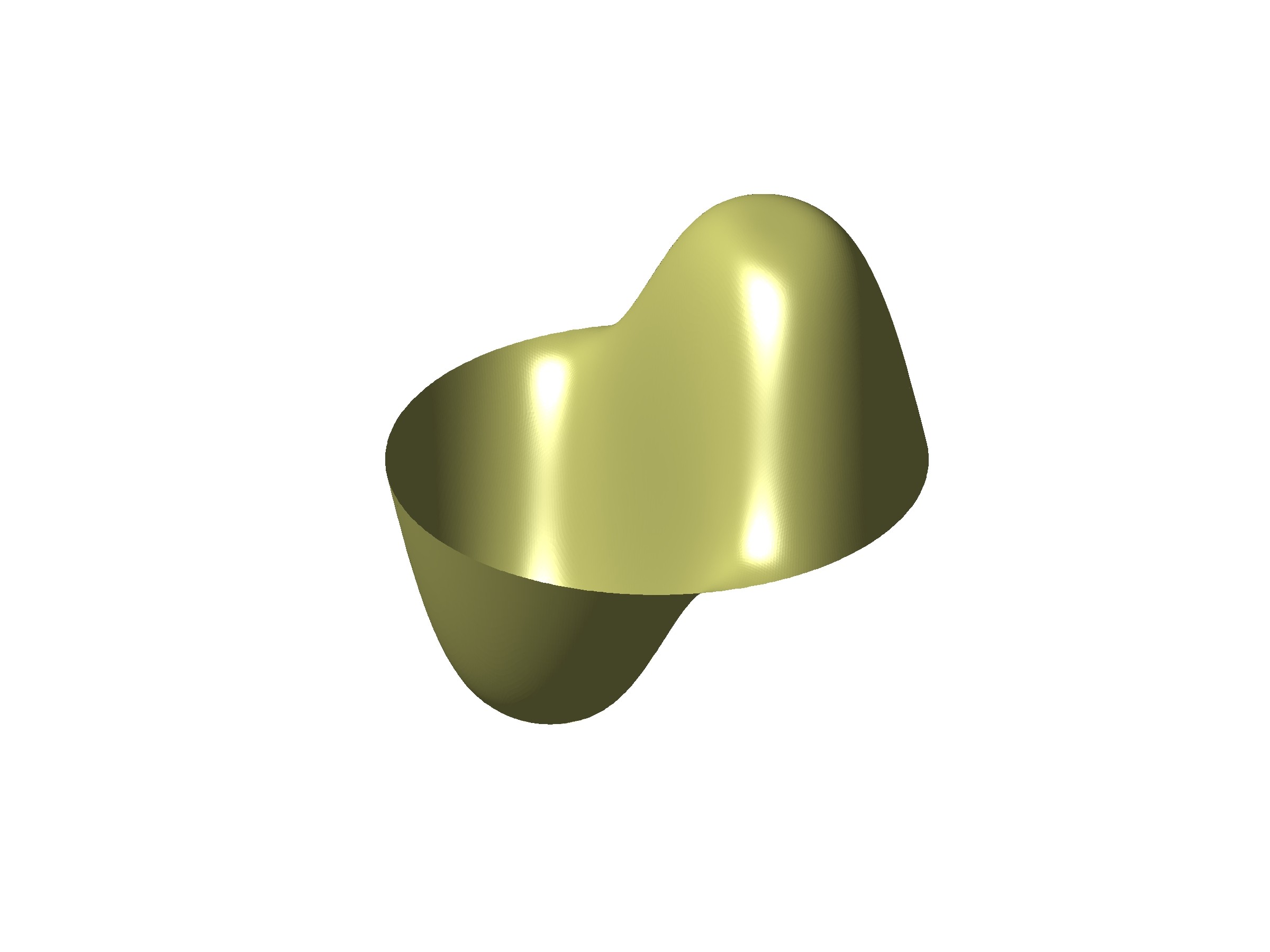}}
    \put(60,0){\dsutwographs{2.0}{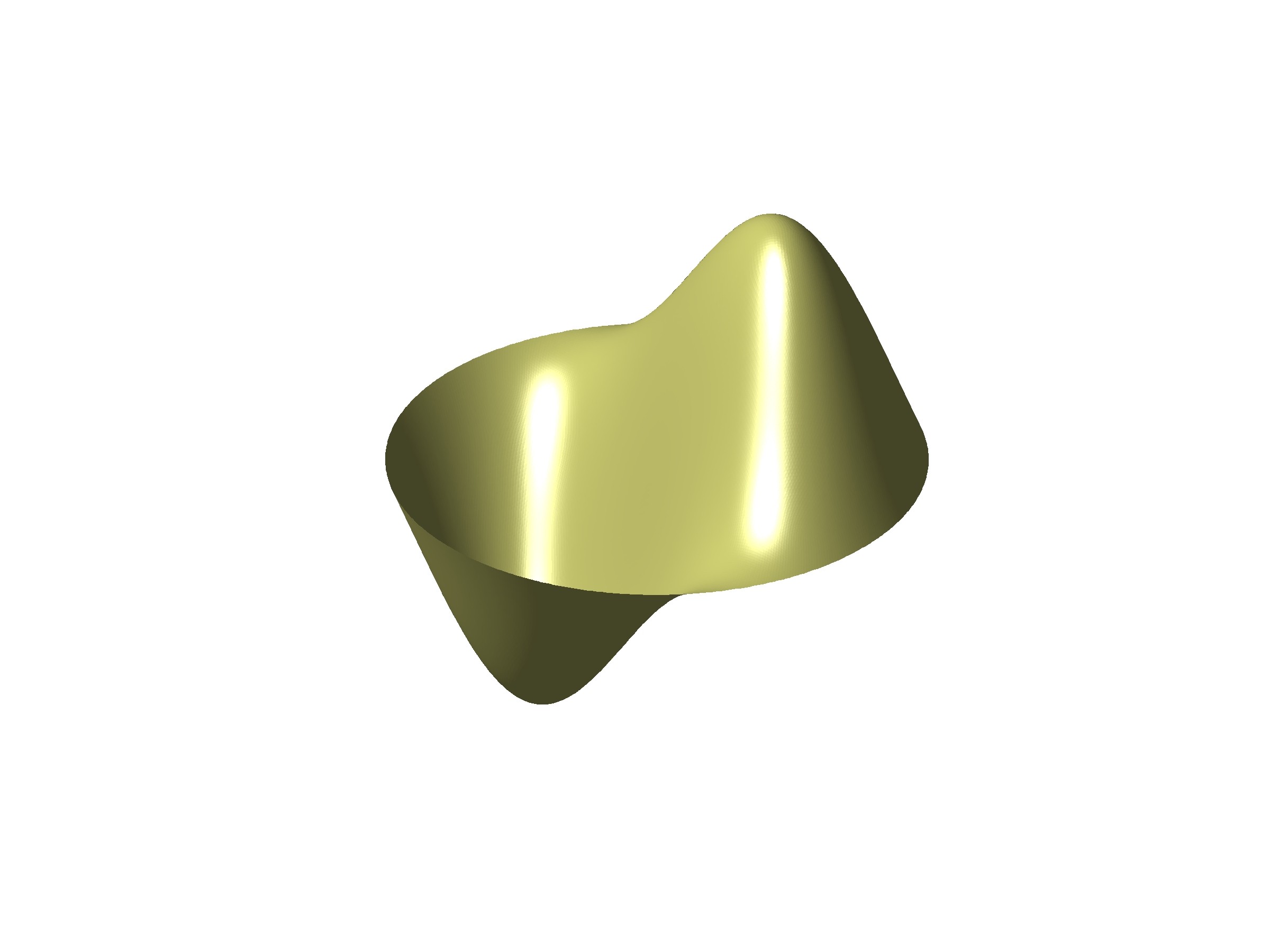}}
    \put(90,0){\dsutwographs{2.5}{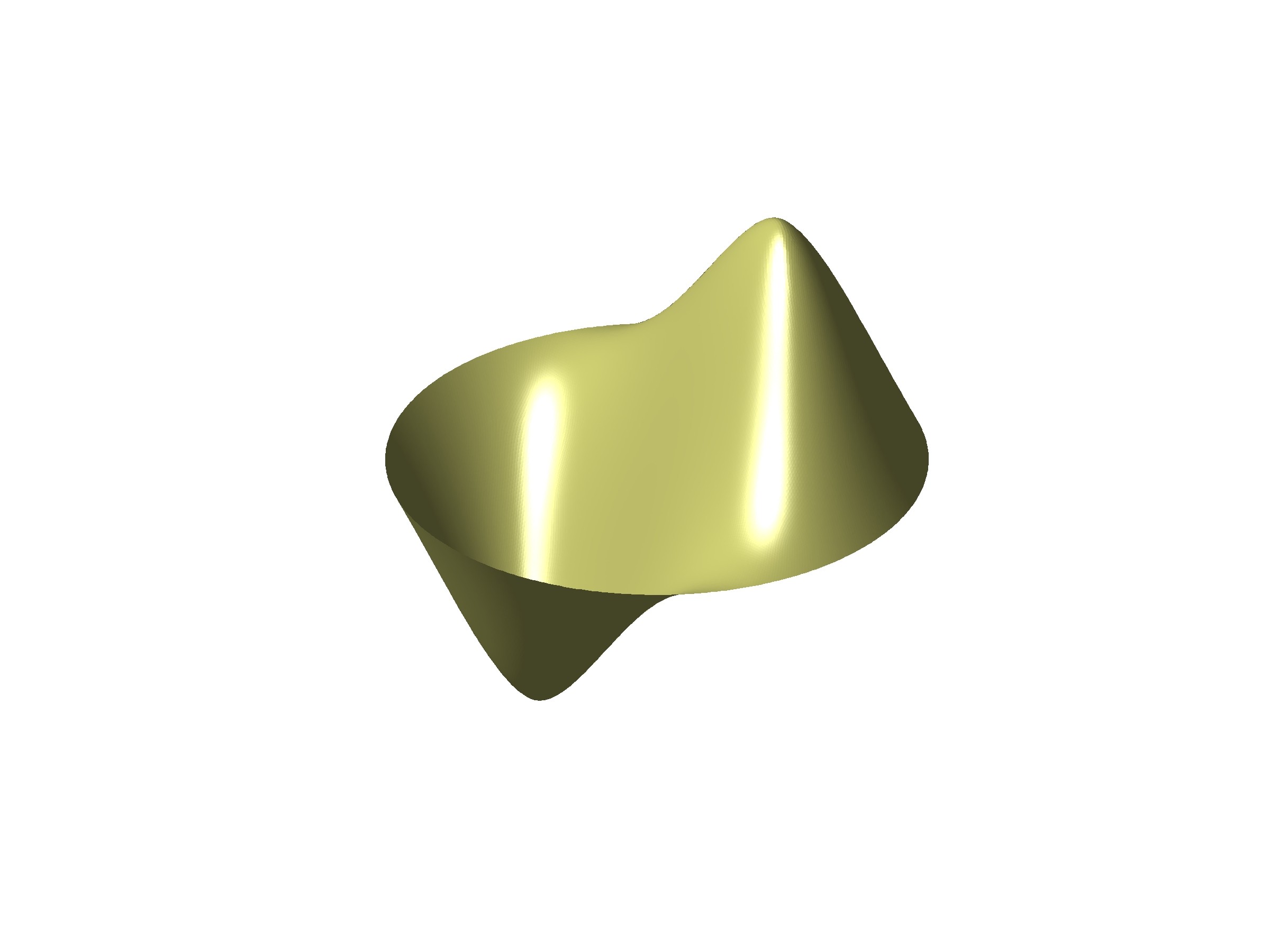}}
    \put(120,0){\dsutwographs{8.0}{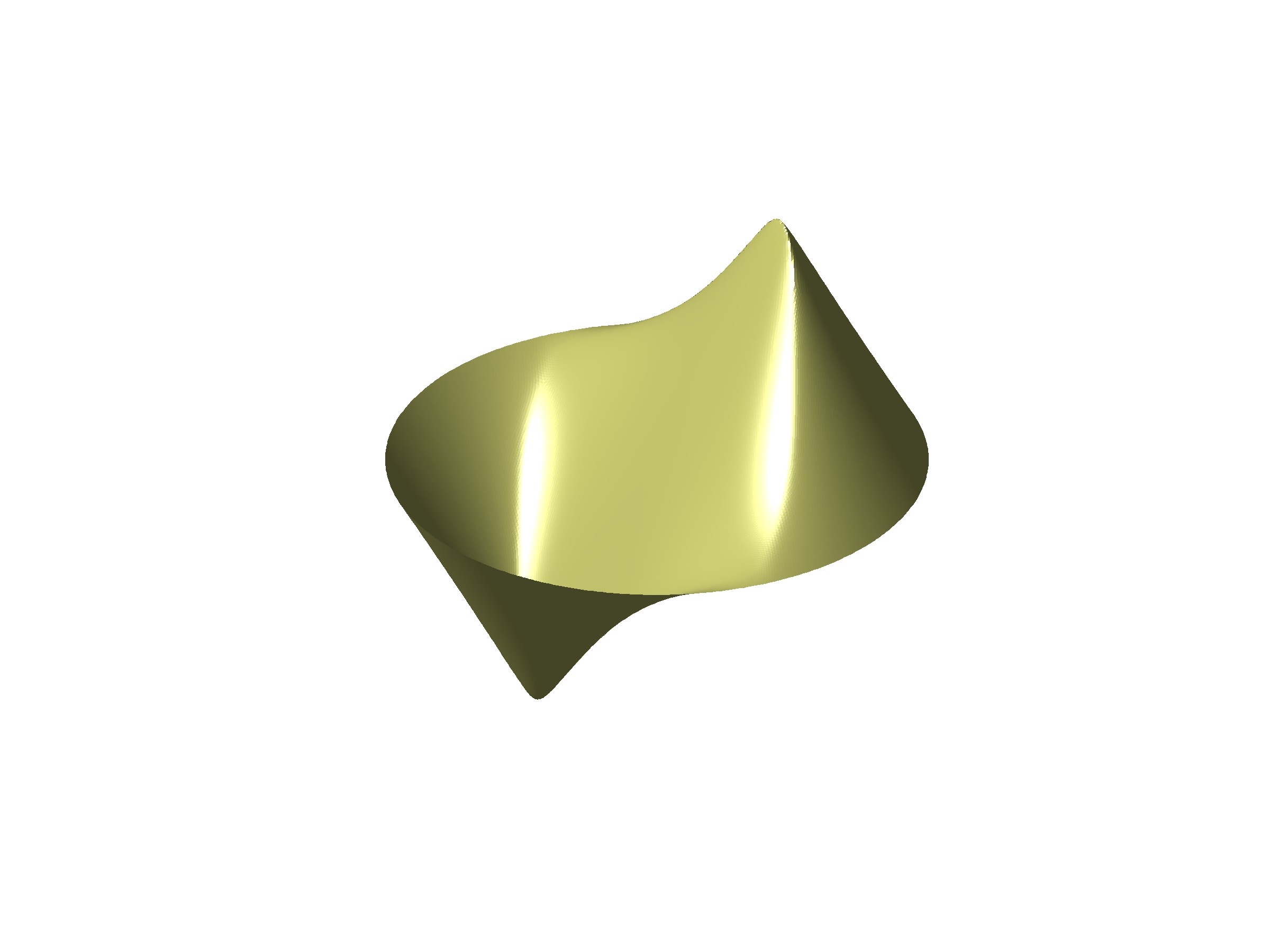}}
  \end{picture}
  \caption{The numerically computed second eigenfunction $u_2$ for the disk.}
  \label{fig:disk_u2}
\end{figure}

\begin{figure}[t]
  \centering
  \setlength{\unitlength}{1mm}
  \begin{picture}(150,54)
    \put(0,0){\dsemgraphs{(a)}{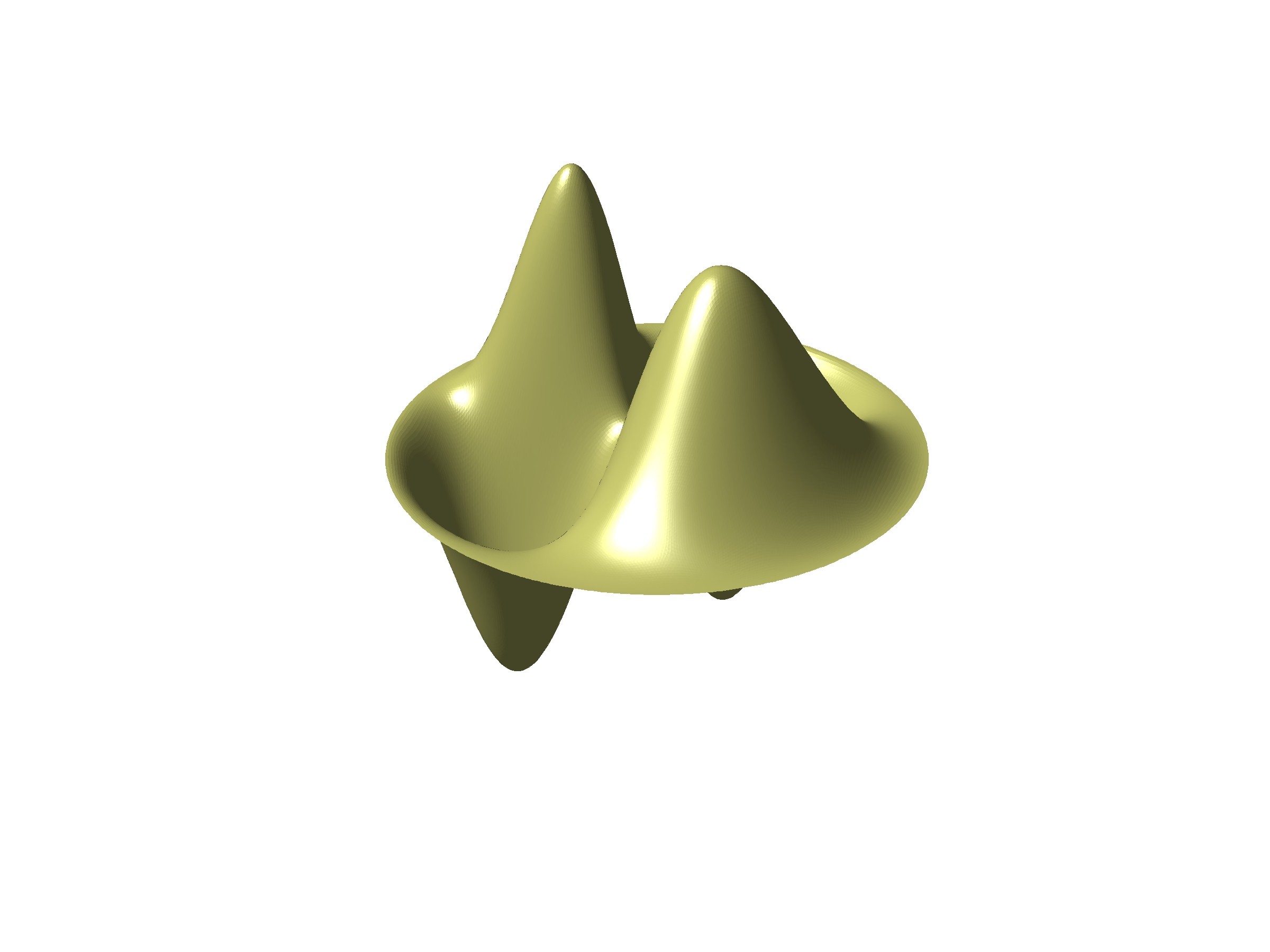}}
    \put(30,0){\dsemgraphs{(b)}{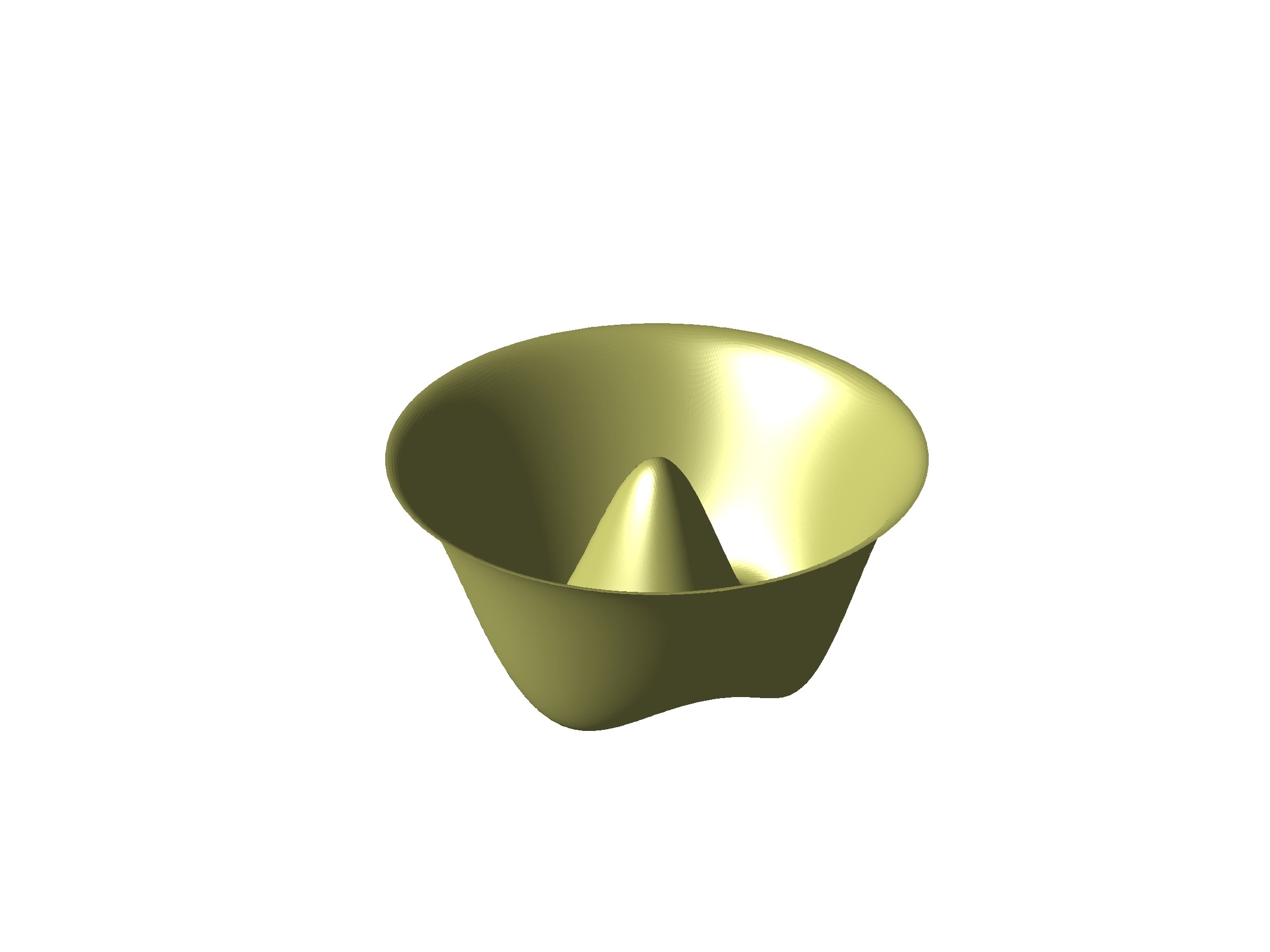}}
    \put(60,0){\dsemgraphs{(c)}{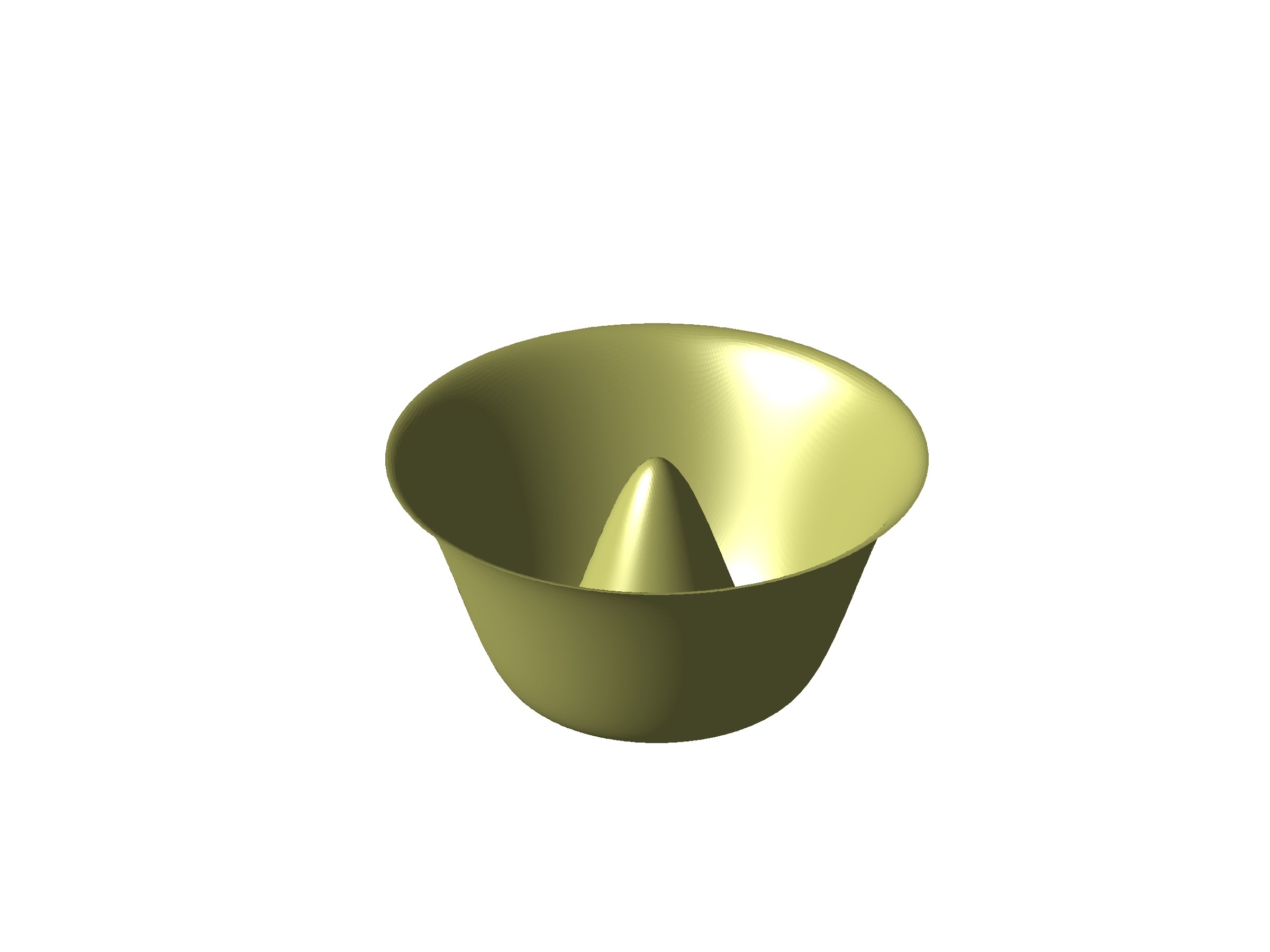}}
    \put(94,2){\includegraphics[width=60mm]{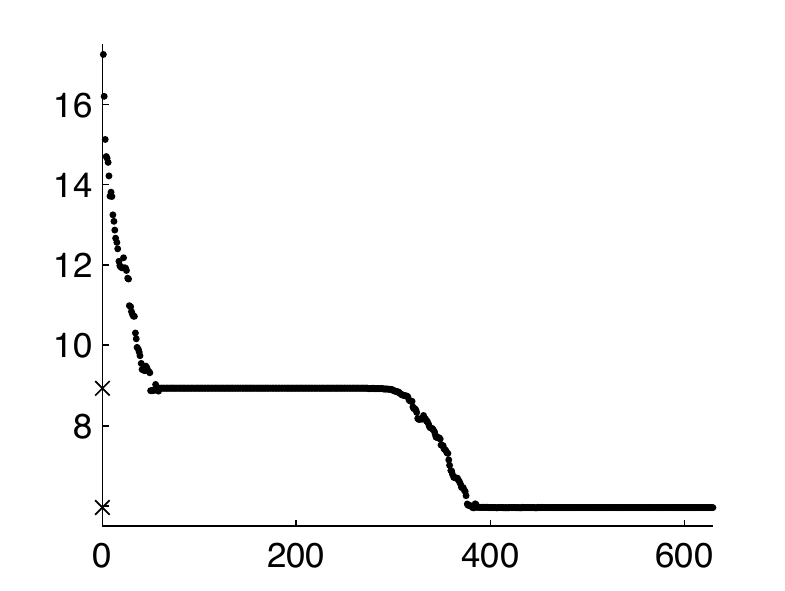}}
    \put(135,0){\small iterations}
    \put(100,45){$I(z^\mathrm{max})$}
    \put(93,7.5){\scriptsize 5.9604}
    \put(93,16.5){\scriptsize 8.9302}
    \put(117,19.5){$u_2^\mathrm{rad}$}
    \put(144,10){$u_2$}
    \put(130,32){$p=1.3$}
    \put(92,46){(d)}
  \end{picture}  
  \caption{(a)--(c) Intermediate point $e_\mathrm{M}$ of the initial
    path used in CMPA. (d) Maximum value of the Dirichlet
    functional $I$ along the path during the run of CMPA with
    $e_\mathrm{M}$ shown in figure (c) for $p=1.3$.}
\label{fig:disk_em}
\end{figure}

Figures \ref{fig:disk_u1} and \ref{fig:disk_u2} show the
eigenfunctions $u_1$ and $u_2$ computed by the genuine 2D method for
several values of $p$, respectively. The corresponding eigenvalues
$\lambda_1$ and $\lambda_2$ for these and other values of $p$ are
listed in Table~\ref{tab:disk_lambda}(a). Figure~\ref{fig:disk_em}(a)
shows the shape of the intermediate point $e_\mathrm{M}$ on the
initial path connecting $u_1$ and $-u_1$ we used for CMPA to find
$u_2$ for all the listed values of $p$. The function $u_2$ found this
way seems to posses an odd symmetry with respect to its nodal
line. The slope of this nodal line in the coordinate system
$(x_1,x_2)$ depends on the computation. For the depiction in
Fig.~\ref{fig:disk_u2} we rotated $\Omega$ in each case to make the
slope appear the same. CMPA needed between 120 and 600 iterations
to converge.

\begin{figure}[p]
  \centering
  \setlength{\unitlength}{1mm}
  \begin{picture}(150,62)
    \put(0,0){\dsutwographs{1.1}{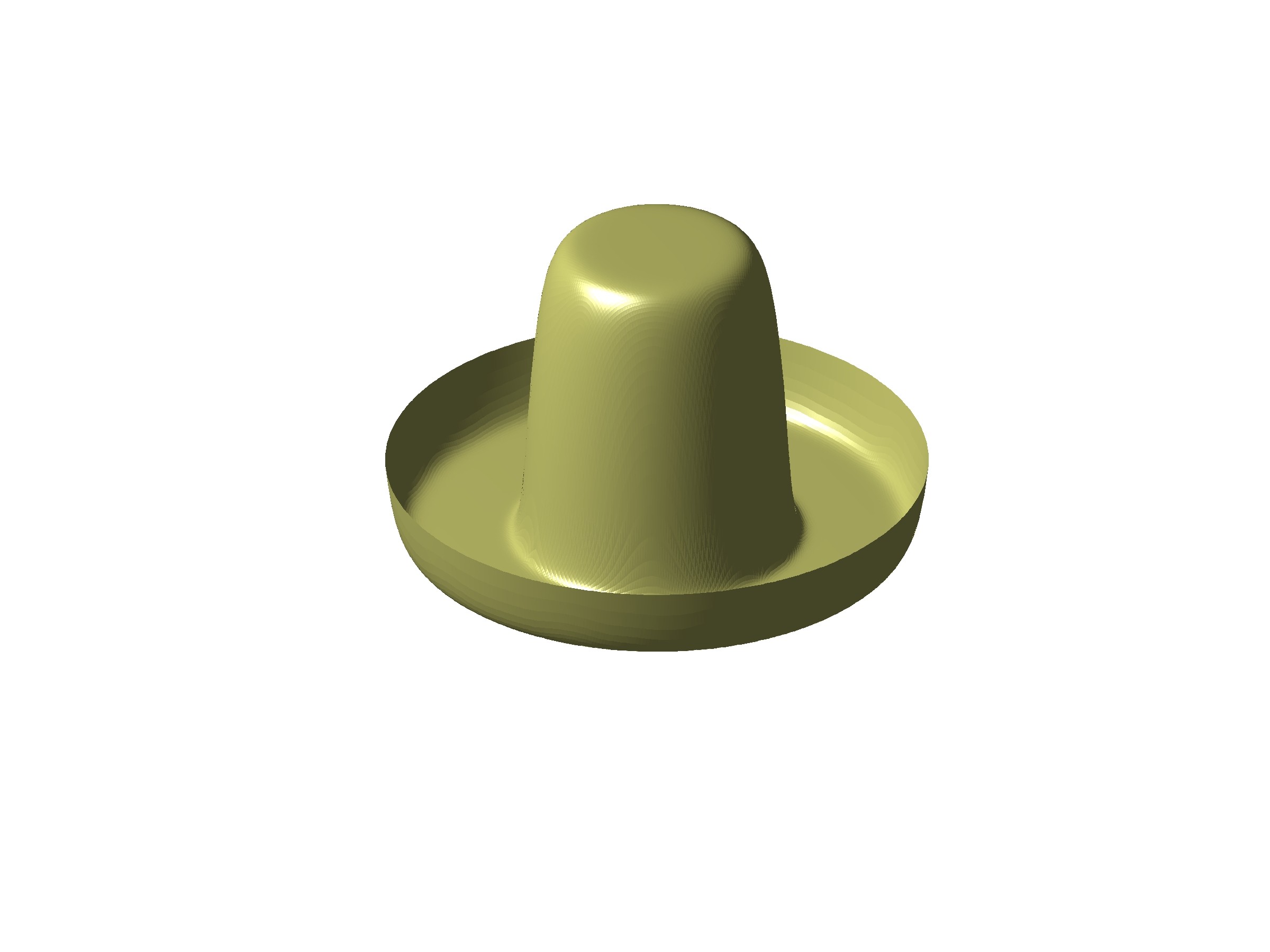}}
    \put(0,57){(a)}
    \put(33,31){\includegraphics[height=31mm]{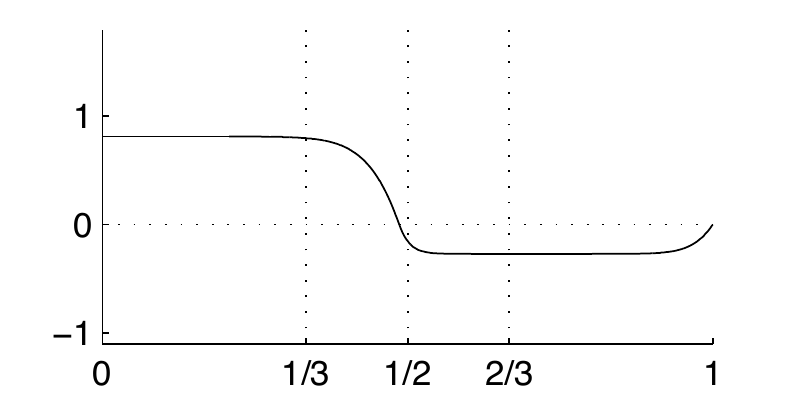}}
    \put(90,31){\includegraphics[height=31mm]{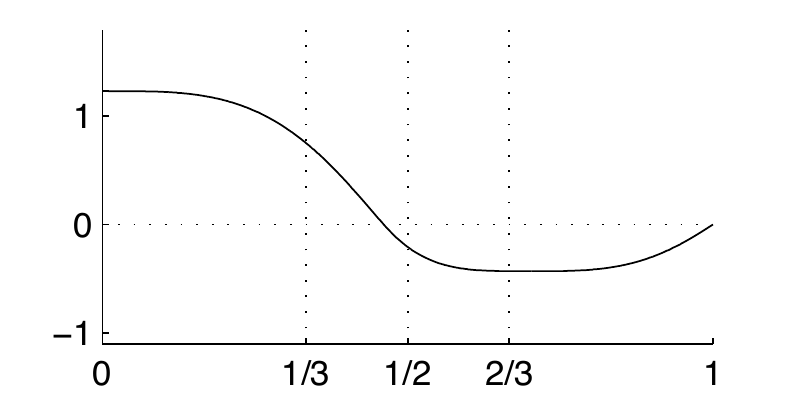}}
    \put(33,0){\includegraphics[height=31mm]{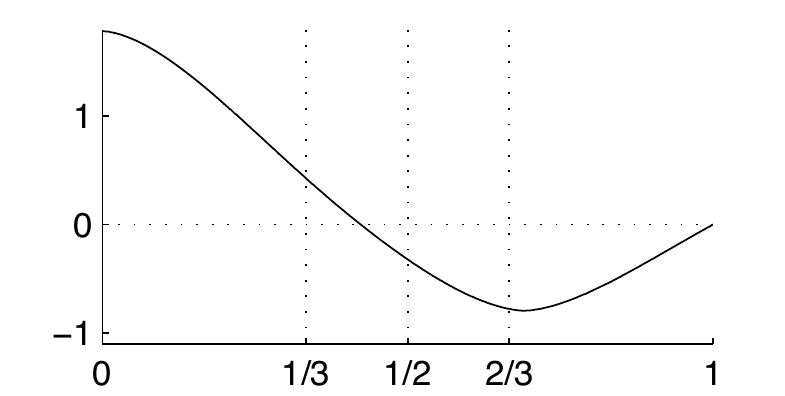}}
    \put(90,0){\includegraphics[height=31mm]{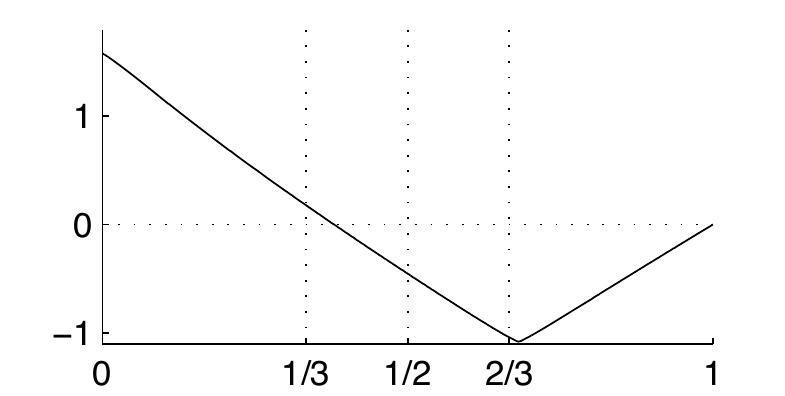}}
    \put(75,54){$p=1.1$}
    \put(132,54){$p=1.4$}
    \put(75,23){$p=2.5$}
    \put(132,23){$p=8.0$}
    \put(33,57){(b)}
  \end{picture}  
  \caption{The numerically computed radially symmetric second
    eigenfunction $u_2^\mathrm{rad}$: (a) using the genuine 2D method;
    (b) using the radial 2D method. The profile of $u_2^\mathrm{rad}$
    for the radial coordinate $r\in(0,1)$ is shown, the scaling along
    the vertical axis is chosen such that $\|u_2^\mathrm{rad}\|_p=1$.}
  \label{fig:disk_u2rad}
\end{figure}
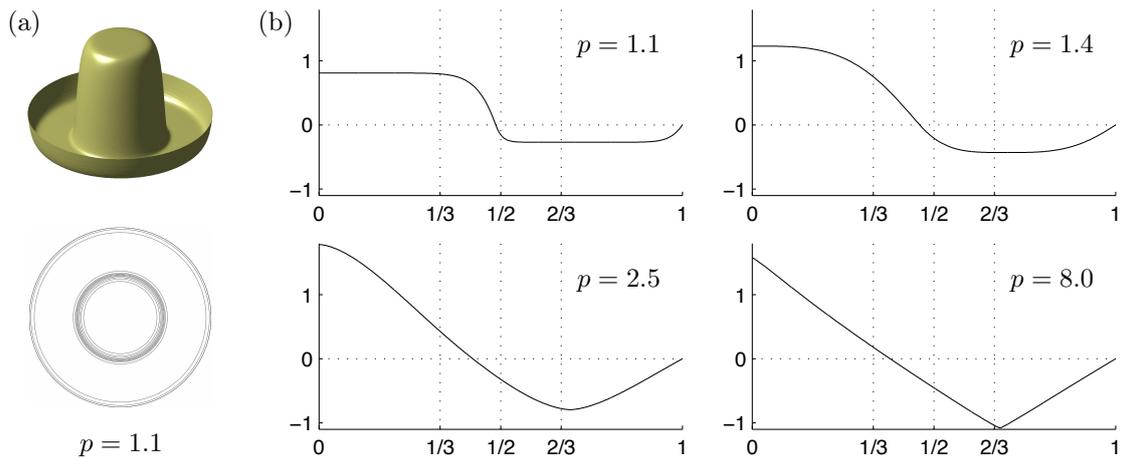

Figure~\ref{fig:disk_em}(b) shows an alternative shape of
$e_\mathrm{M}$. With such an initial path CMPA converged for
$p=1.1$ and $p=1.2$ to a radially symmetric function we call
$u_2^\mathrm{rad}$ (but for higher values of $p$ to the oddly
symmetric function $u_2$). Figure~\ref{fig:disk_u2rad}(a) shows
$u_2^\mathrm{rad}$ for $p=1.1$.

Figure~\ref{fig:disk_em}(c) shows yet another choice of $e_\mathrm{M}$
(radially symmetric). With this intermediate point of the initial path
and for $p=1.3$ and $p=1.4$ (but not larger) CMPA seems to
converge to a radially symmetric function first but after many
iterations the path slips down and the algorithm converges eventually
to the oddly symmetric $u_2$. The graph in Fig.~\ref{fig:disk_em}(d)
shows how the maximum value of the Dirichlet functional $I$ along the
path develops during the run of the algorithm (for $p=1.3$). The
horizontal axis shows the number of iterations. The flat part between
iterations 70 and 260 indicates that the path is staying close to a
critical point. When now the norm of the descent direction
$w_{z^\mathrm{max}}^{}$ given in \pref{eq:descent_direction} computed at
the ``highest'' point $z^\mathrm{max}$ of the path gets small enough,
we stop the algorithm and save this highest point. Since it displays a
radial symmetry, we call it $u_2^\mathrm{rad}$ again.

The eigenvalues $\lambda_2^\mathrm{rad}$ corresponding to the found
$u_2^\mathrm{rad}$ are also listed in Table~\ref{tab:disk_lambda}(a).

Figure~\ref{fig:disk_u2rad}(b) shows profiles of the eigenfunction
$u_2^\mathrm{rad}$ computed by the radial 2D method for several values
of $p$. The eigenvalues $\lambda_1$ and $\lambda_2^\mathrm{rad}$
computed by this method for these and other values of $p$ are listed
in Table~\ref{tab:disk_lambda}(b). The convergence of CMPA does
not seem to be sensitive to the choice of $e_\mathrm{M}$ in this case.

By comparing the values of $\lambda_1$ and $\lambda_2^\mathrm{rad}$ in
Table~\ref{tab:disk_lambda}(a) with those in
Table~\ref{tab:disk_lambda}(b) which were computed by the two
different numerical methods we observe that their first three digits
coincide in almost all the cases. Also, the profiles of $u_1$,
$u_2^\mathrm{rad}$ are very close for both methods, respectively
(cf.~Fig~\ref{fig:disk_u2rad}(a) and the top left graph in (b) for
$u_2^\mathrm{rad}$ and $p=1.1$). We conclude that these are numerical
approximations of the same eigenvalue-eigenfunction pairs.

The behavior of CMPA suggests that although $u_2^\mathrm{rad}$ is
a constrained mountain pass point of $I$ among radially symmetric
functions, it is not a constrained mountain pass point with no
assumption on the symmetry (cf.~Fig.~\ref{fig:disk_em}(d)). The case
of $p=1.1$ and $p=1.2$ when CMPA with $e_\mathrm{M}$ from
Fig.~\ref{fig:disk_em}(b) converged to a radially symmetric function
and the path did not slip off to asymmetric functions with lower
values of $I$ seems to contradict this. However, we assume that this
was caused by the ``flat'' shape of the landscape of $I$ close to
$u_2^\mathrm{rad}$ for $p$ close to 1 and by numerical inaccuracies.

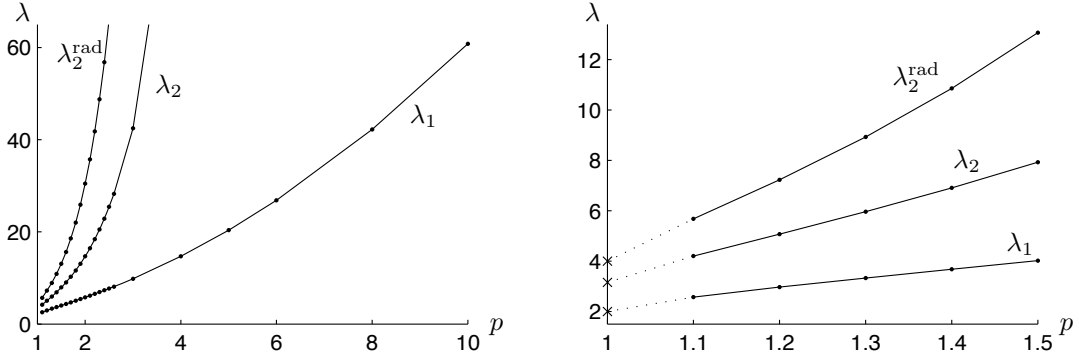
\begin{figure}[t]
  \centering
  \setlength{\unitlength}{1mm}
  \begin{picture}(150,50)(0,2)
    \put(0,0){\includegraphics[width=7.3cm]{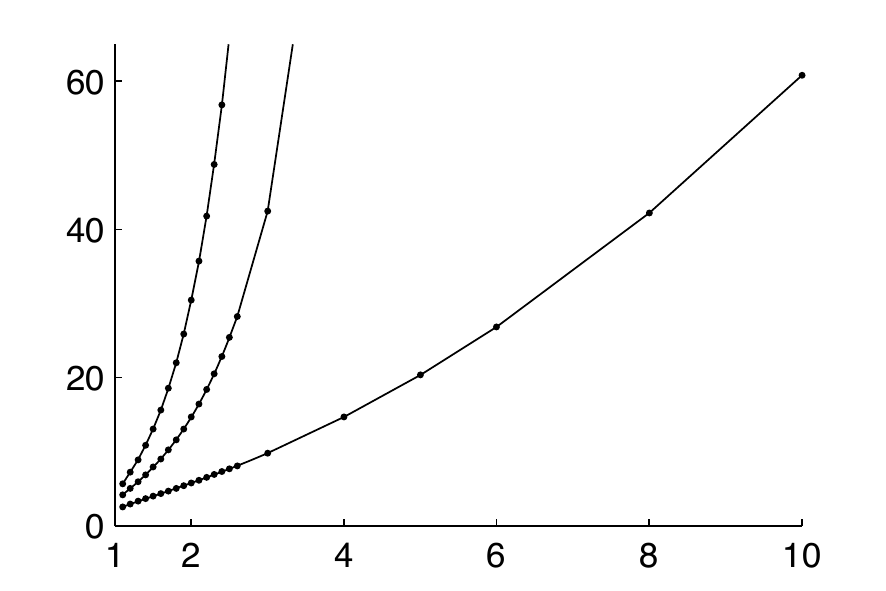}}
    \put(6.5,45.5){$\lambda$}
    \put(69,4.5){$p$}
    \put(58.5,32){$\lambda_1$}
    \put(25,36){$\lambda_2$}
    \put(12,40){$\lambda_2^\mathrm{rad}$}
    \put(75,0){\includegraphics[width=7.3cm]{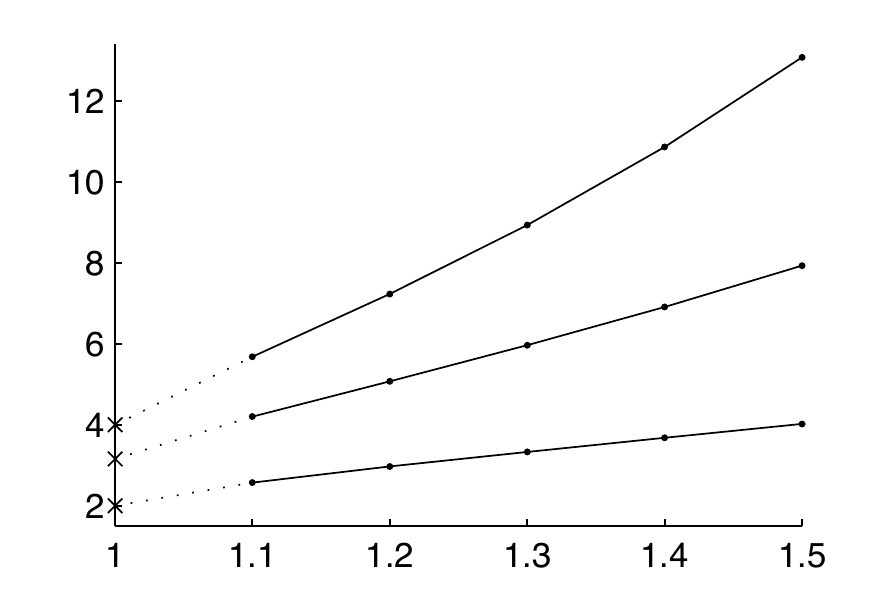}}
    \put(81.5,45.5){$\lambda$}
    \put(144,4.5){$p$}
    \put(137,15){$\lambda_1$}
    \put(130,26){$\lambda_2$}
    \put(122,37){$\lambda_2^\mathrm{rad}$}
  \end{picture}  
  \caption{Dependence of the numerically computed eigenvalues for the disk on $p$. The
    three cross symbols in the graph on the right mark the values
    $h_1(\Omega)=2$, $h_2(\Omega)\approx 3.1543$, and
    $h_2^\mathrm{rad}(\Omega)=4$.}
  \label{fig:disk_lambda_plot1}
\end{figure}

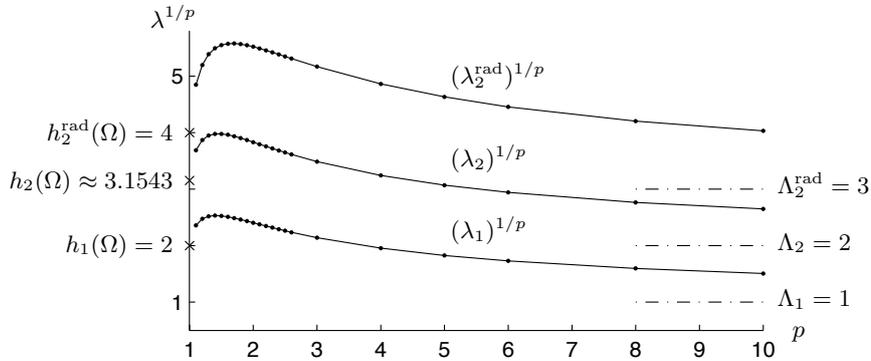
\begin{figure}[t]
  \centering
  \setlength{\unitlength}{1mm}
  \begin{picture}(130,50)(-15,2)\small
    \put(0,0){\includegraphics[width=9.73cm]{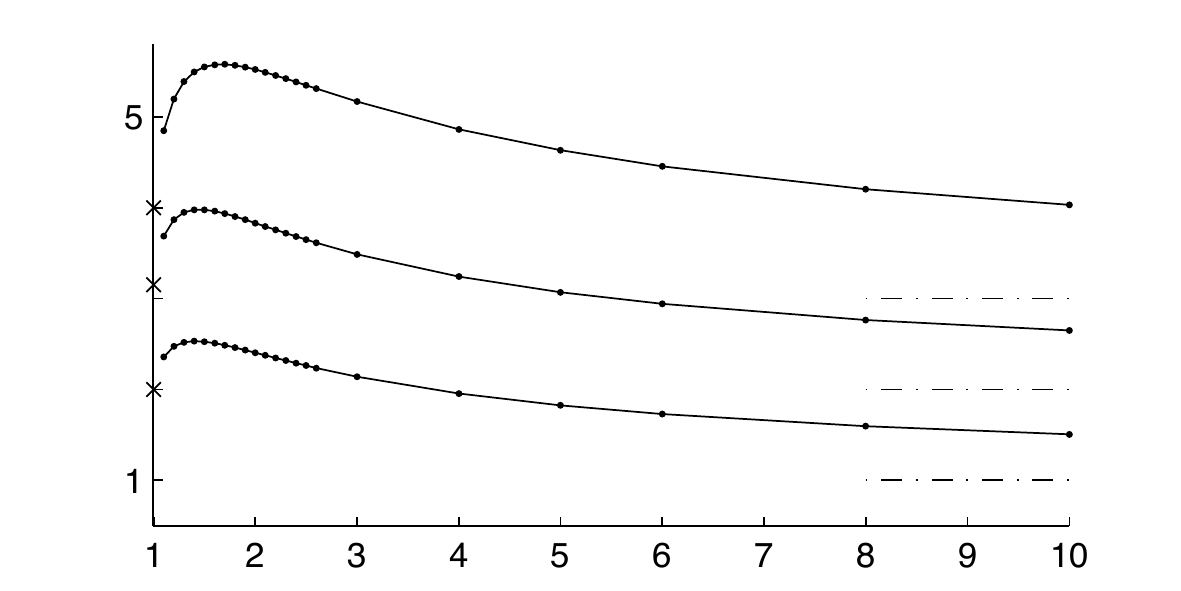}}
    \put(7.5,45.5){$\lambda^{1/p}$}
    \put(92,4.5){$p$}
    \put(47,17.5){$(\lambda_1)^{1/p}$}
    \put(47,27){$(\lambda_2)^{1/p}$}
    \put(47,38){$(\lambda_2^\mathrm{rad})^{1/p}$}
    \put(-14.5,14.5){\makebox(25,4)[r]{$h_1(\Omega)=2$}}
    \put(-14.5,23){\makebox(25,4)[r]{$h_2(\Omega)\approx 3.1543$}}
    \put(-14.5,29.5){\makebox(25,4)[r]{$h_2^\mathrm{rad}(\Omega)=4$}}
    \put(90,8.5){$\Lambda_1=1$}
    \put(90,16){$\Lambda_2=2$}
    \put(90,23.5){$\Lambda_2^\mathrm{rad}=3$}
  \end{picture}  
  \caption{Dependence of the numerically computed eigenvalues for the disk
    raised to $1/p$ on $p$.}
  \label{fig:disk_lambda_plot2}
\end{figure}

The dependence of $\lambda_1$, $\lambda_2$, and
$\lambda_2^\mathrm{rad}$ on $p$ is presented in
Figs.~\ref{fig:disk_lambda_plot1} and
\ref{fig:disk_lambda_plot2}. First of all we observe that for all the
values of $p$ considered the inequality
$\lambda_2<\lambda_2^\mathrm{rad}$ holds. Hence this is a numerical
evidence that the second eigenfunction for the disk is not radially
symmetric not only for small $p$ but for a large range of $p$.

Second, we can observe the following asymptotic behavior:
\begin{center}\renewcommand{\arraystretch}{1.3}\vspace{1ex}
  \begin{tabular}{|c|c|c|c|}
    \cline{2-4}
    \multicolumn{1}{c|}{} & $\lambda_1$ & $\lambda_2$ & $\lambda_2^\mathrm{rad}$ \\\hline
    $\lim_{p\to1+}\lambda$ & 2 & 3.1543 & 4 \\\hline
    $\lim_{p\to\infty}\lambda^{1/p}$ & 1 & 2 & 3 \\\hline
  \end{tabular}\renewcommand{\arraystretch}{1}\vspace{1ex}
\end{center}
Theoretical results for $\lambda_1$ and $\lambda_2$ were summarize in
Sec.~\ref{sec:background}. The values $h_1(\Omega)$,
$\Lambda_1(\Omega)$, and $\Lambda_2(\Omega)$ for the disk are easy to
compute. In \cite{Parini} it was proved that $h_2(\Omega)$ for the
disk equals the first Cheeger constant for the half-disk which is
approximately 3.1543. We can observe (Fig.~\ref{fig:disk_u1}) that
$u_1$ converges to 1 for $p\to 1$ as explained in \cite{KawohlFridman}
and to the distance function to the boundary for $p\to\infty$ as
explained in \cite{JuutinenLindqvistManfredi}. In
Fig.~\ref{fig:disk_u2} we observe that for $p\to1$ the function $u_2$
is getting close to the indicator function of the Cheeger set for the
half-disk on each nodal domain.

In \cite{BenDrabGirg} a numerical evidence is given leading to the
conjecture for the asymptotic behavior of $\lambda_2^\mathrm{rad}$
given in the above table. Our numerical results (at least for $p\to1$)
support this conjecture. To motivate these values and the profiles of
$u_2^\mathrm{rad}$ in Fig.~\ref{fig:disk_u2rad} we make the following
two remarks:
\begin{plapnum_remark}
  1.~For $r\in(0,1)$ let $D(r)$ be the disk of radius $r$ centered at
  the origin, $A(r)=\Omega\setminus\overline{D(r)}$ an annulus. It is
  easy to show that for $r=1/2$ both $D$ and $A$ have the same Cheeger
  constant $h_2^\mathrm{rad}(\Omega):=h_1(D(1/2))=h_1(A(1/2))=4$ (see,
  e.g., \cite{KrejcirikPratelli} for a result about the Cheeger constant of an
  annulus). The function $u_2^\mathrm{rad}$ with its profile shown in
  Fig.~\ref{fig:disk_u2rad} seems to get close to the indicator
  function of $D(1/2)$ and $A(1/2)$ on each nodal domain for $p\to1$.

  2.~Under the assumption of radial symmetry two largest disjoint
  disks of the same radius inscribed in $\Omega$ have radius
  1/3. Hence we define $\Lambda_2^\mathrm{rad}=\frac{1}{1/3}=3$. The
  function $u_2^\mathrm{rad}$ with its profile shown in
  Fig.~\ref{fig:disk_u2rad} seems to get close on each nodal domain to
  a multiple of the function giving the distance to the boundary on
  $D(1/3)$ and $A(1/3)$ for large $p$.
\end{plapnum_remark}

\subsection{Square}
\begin{figure}[t]
  \centering
  \setlength{\unitlength}{1mm}
  \begin{picture}(150,56)
    \put(0,0){\squonegraphs{1.1}{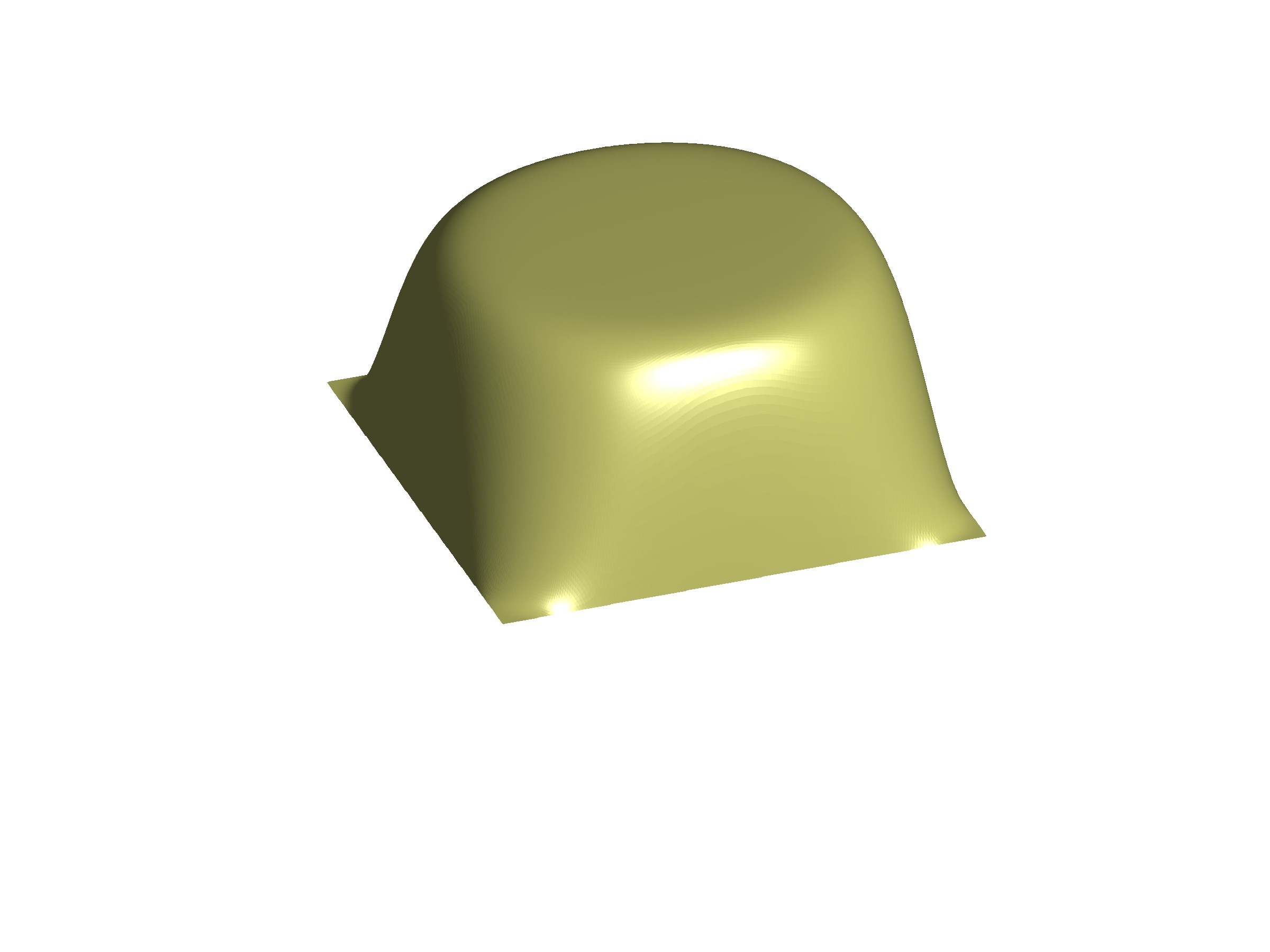}}
    \put(30,0){\squonegraphs{1.4}{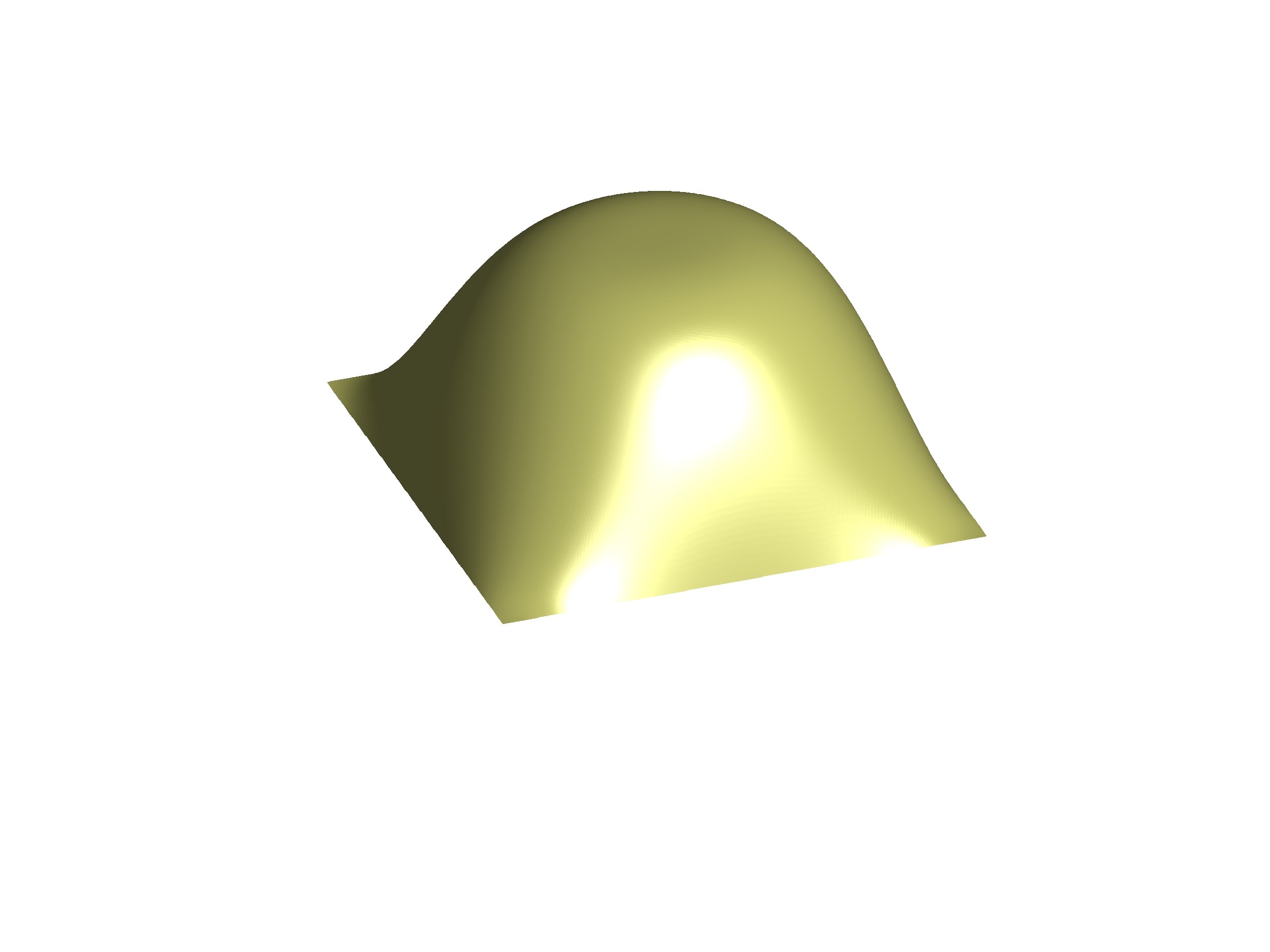}}
    \put(60,0){\squonegraphs{2.0}{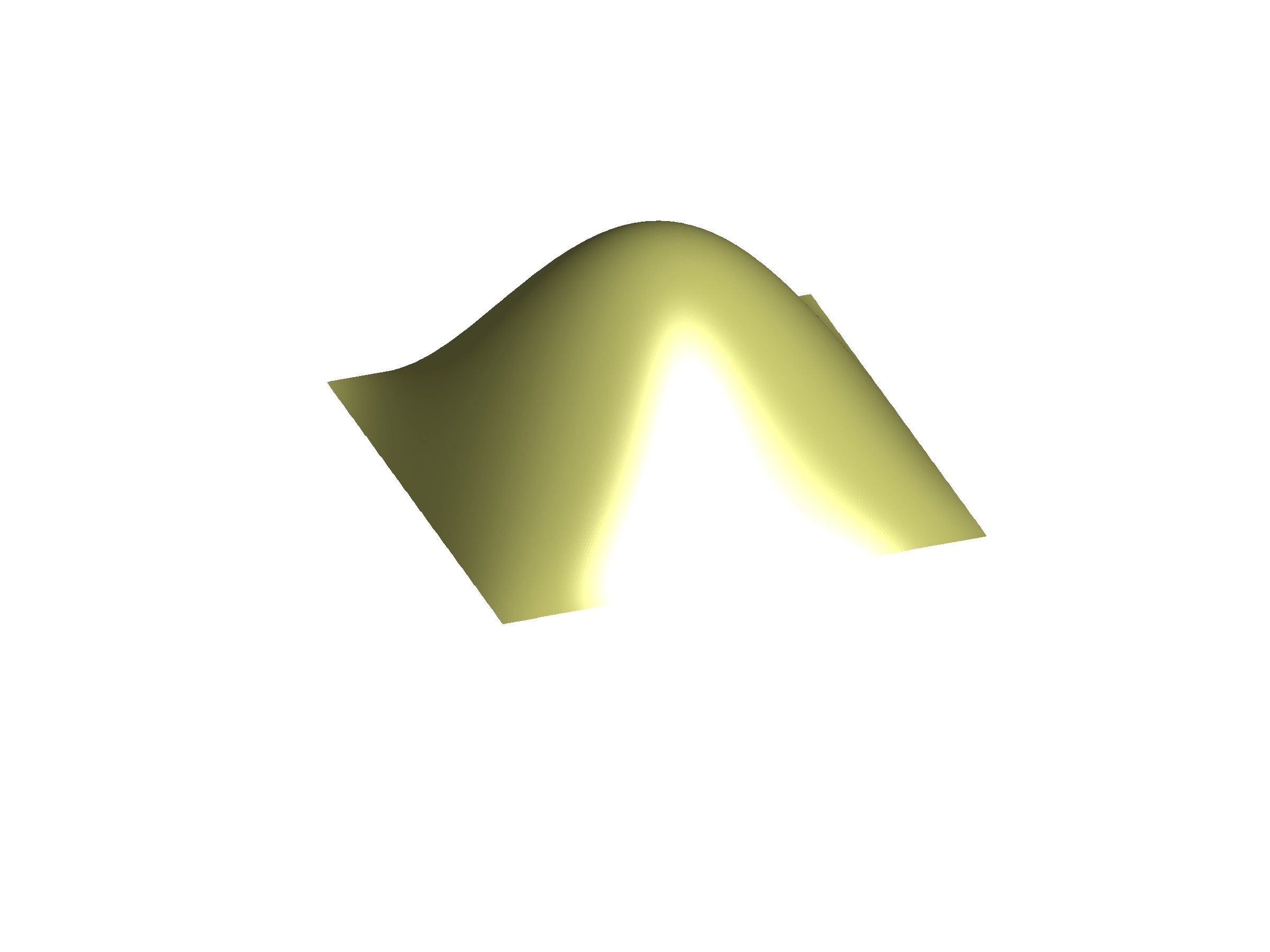}}
    \put(90,0){\squonegraphs{2.5}{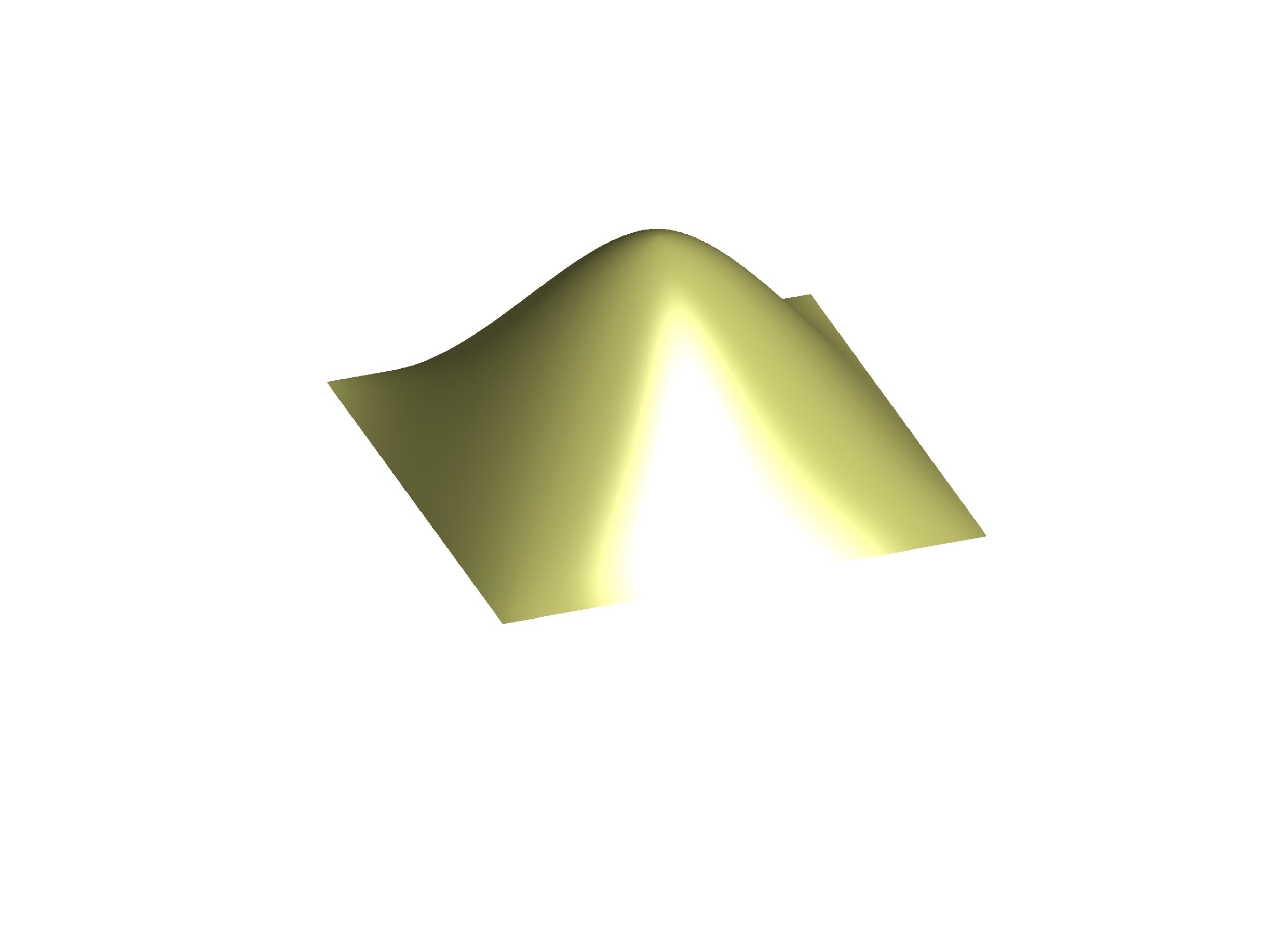}}
    \put(120,0){\squonegraphs{8.0}{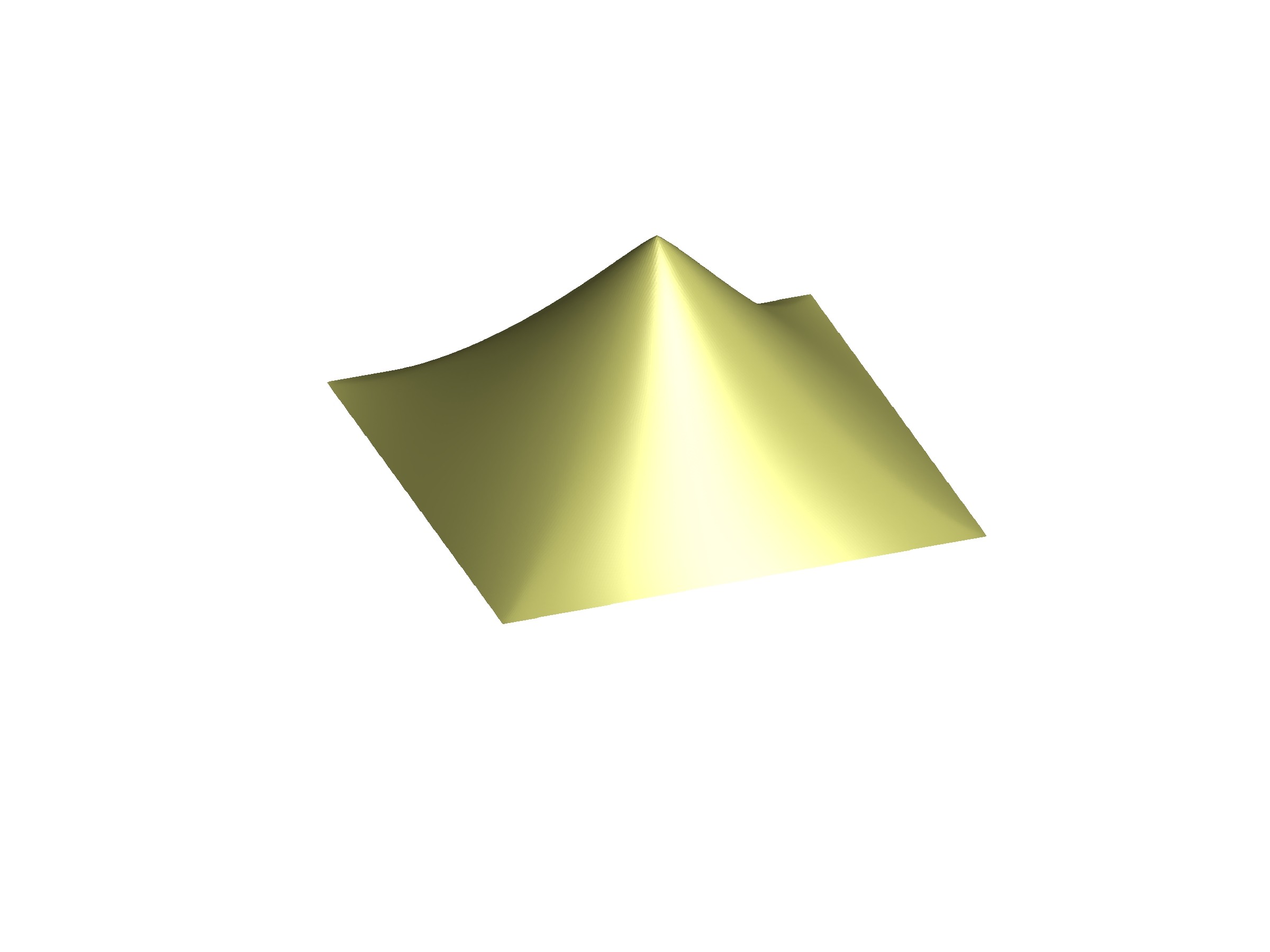}}
  \end{picture}
  \caption{The numerically computed first eigenfunction $u_1$ for the square.}
  \label{fig:square_u1}
\end{figure}

\begin{figure}[h]
  \centering
  \setlength{\unitlength}{1mm}
  \begin{picture}(150,60)
    \put(0,0){\squtwographs{1.1}{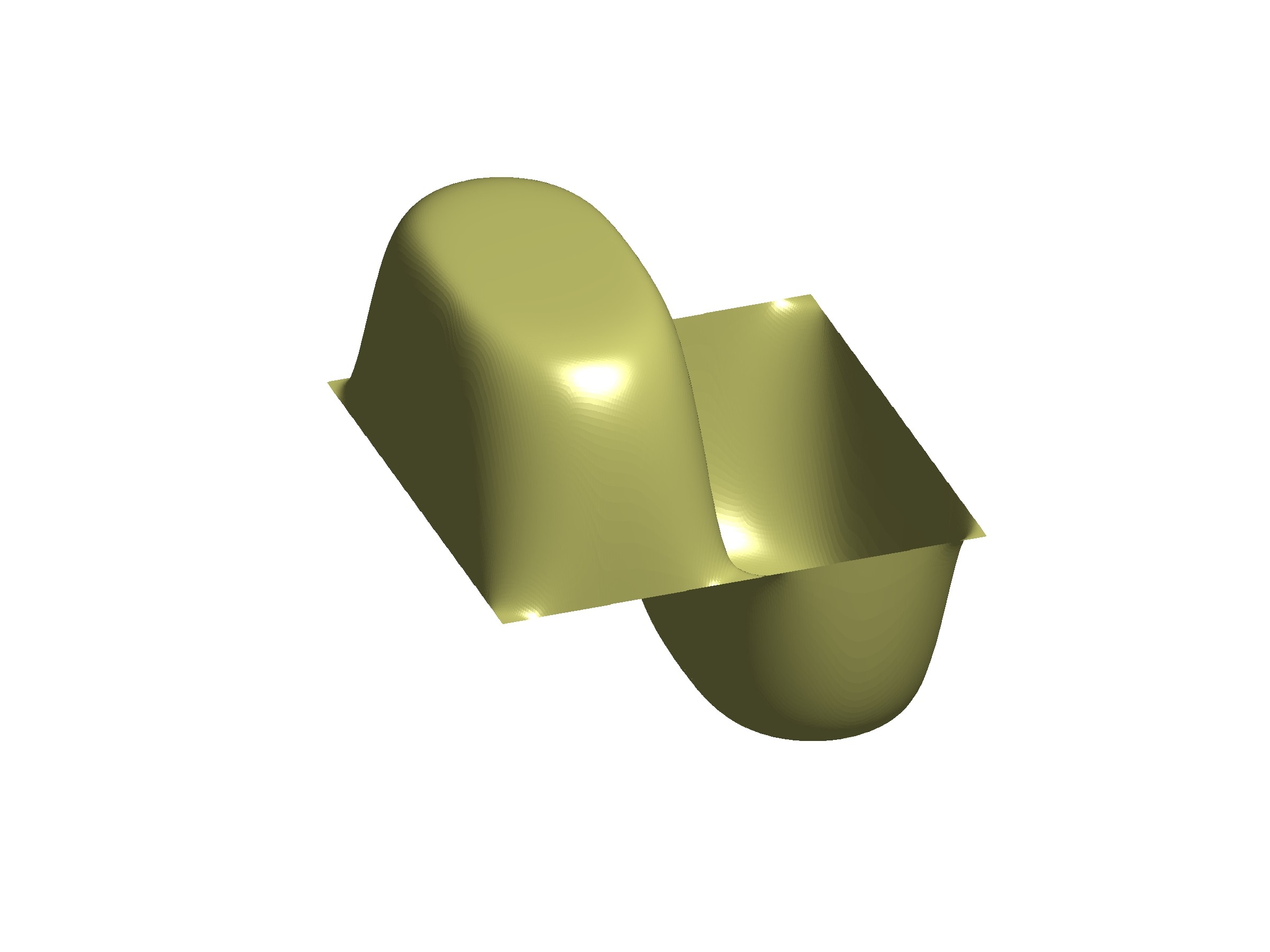}}
    \put(30,0){\squtwographs{1.4}{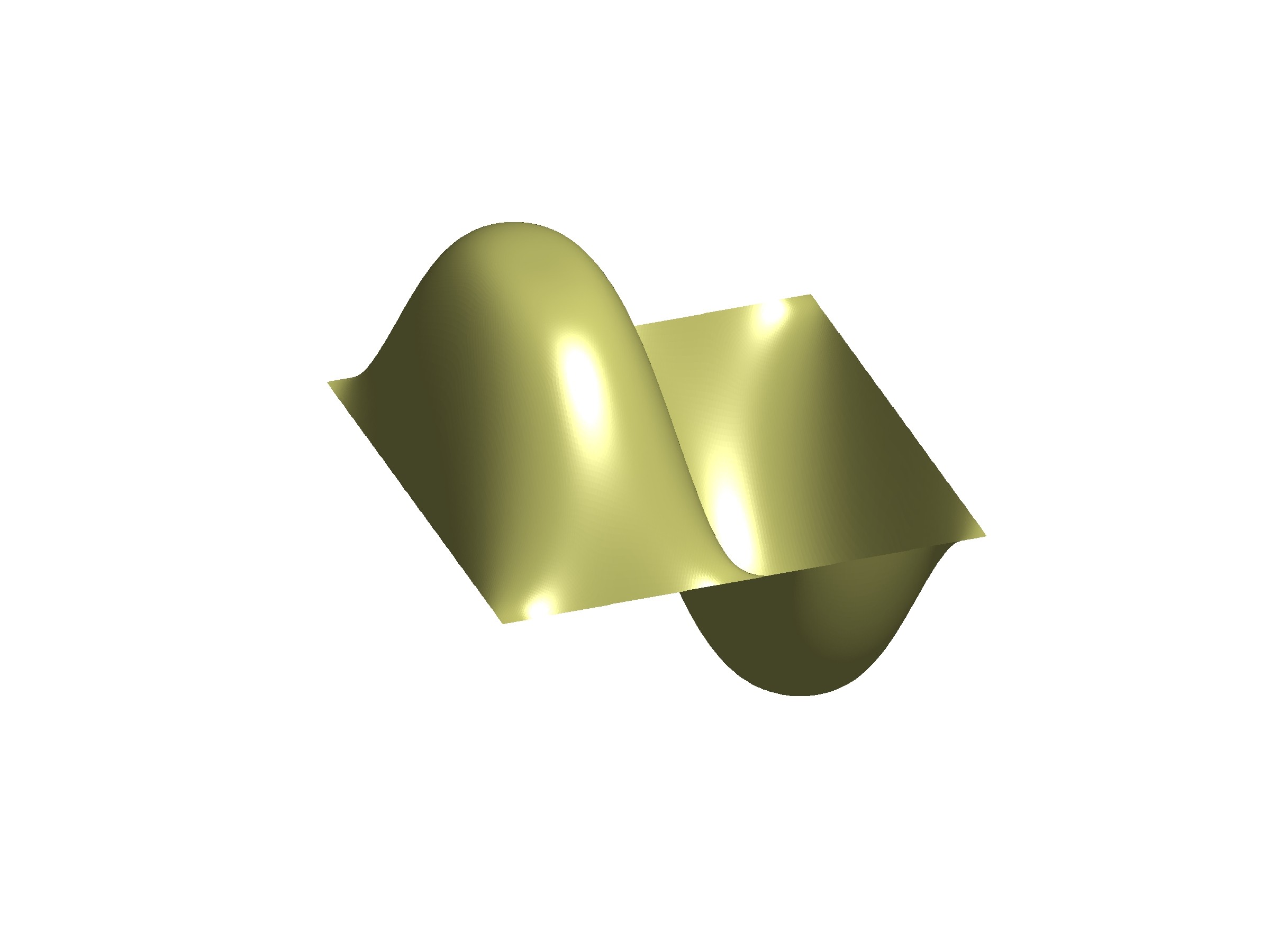}}
    \put(60,0){\squtwographs{2.0}{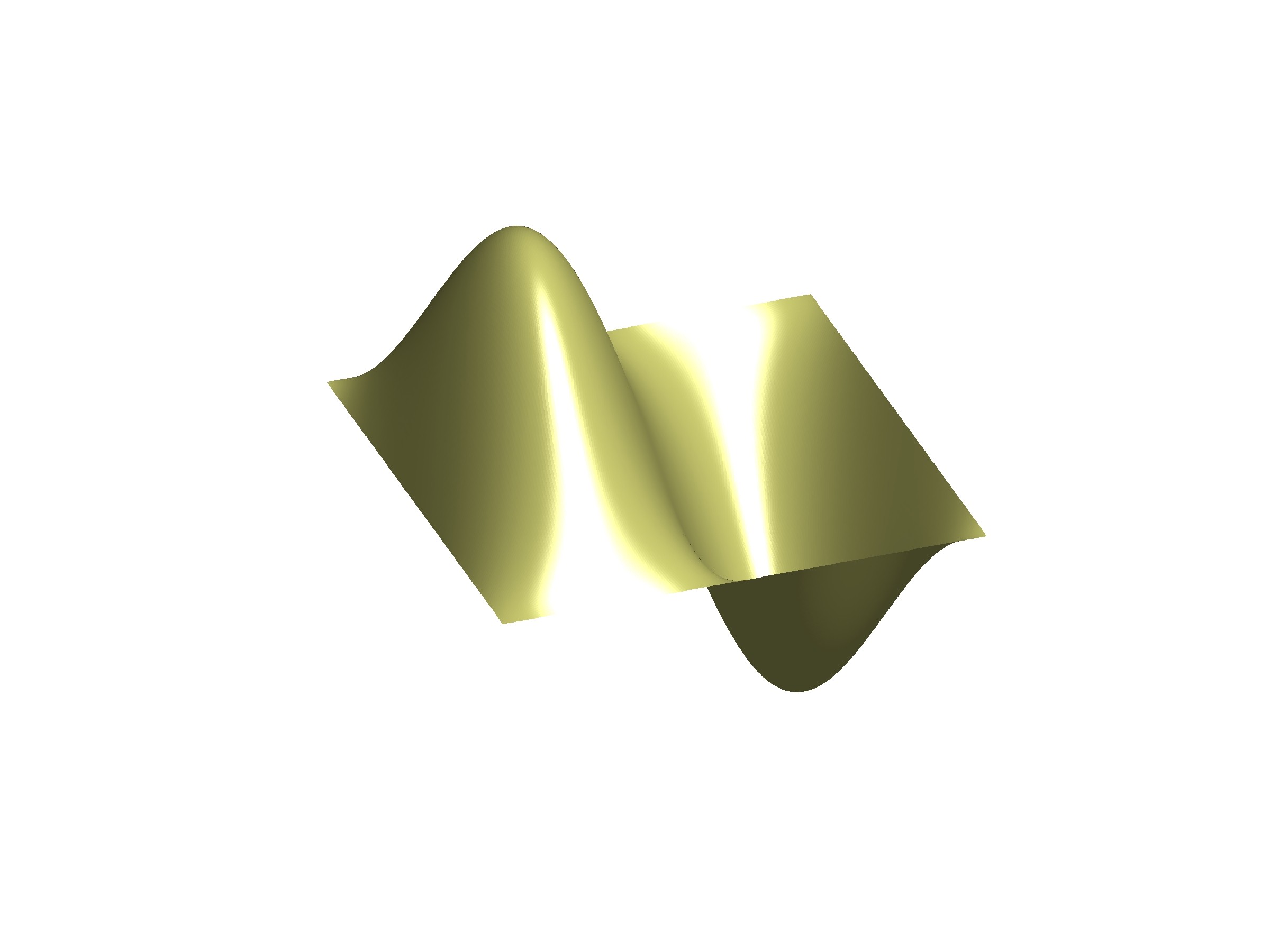}}
    \put(90,0){\squtwographs{2.5}{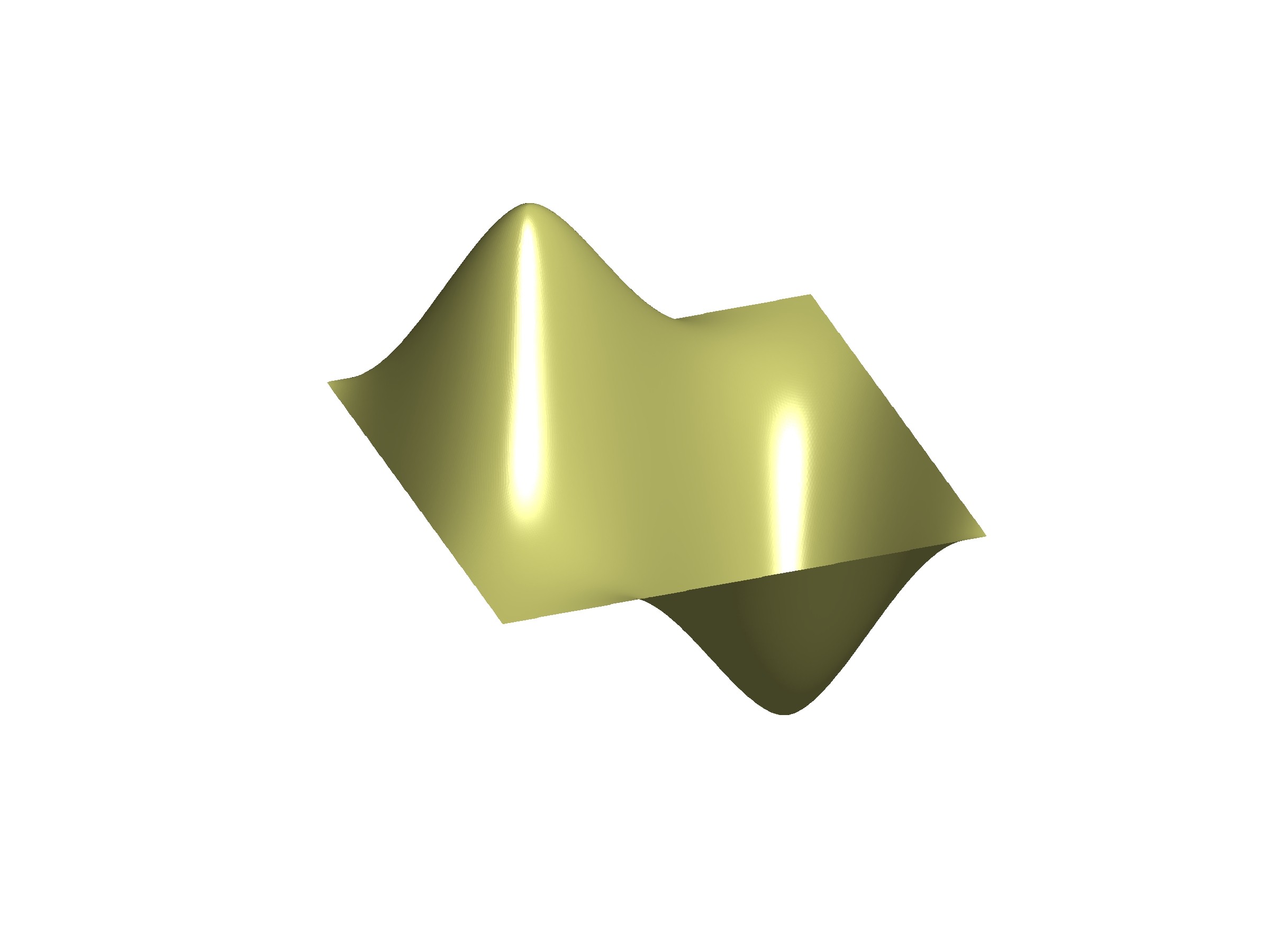}}
    \put(120,0){\squtwographs{8.0}{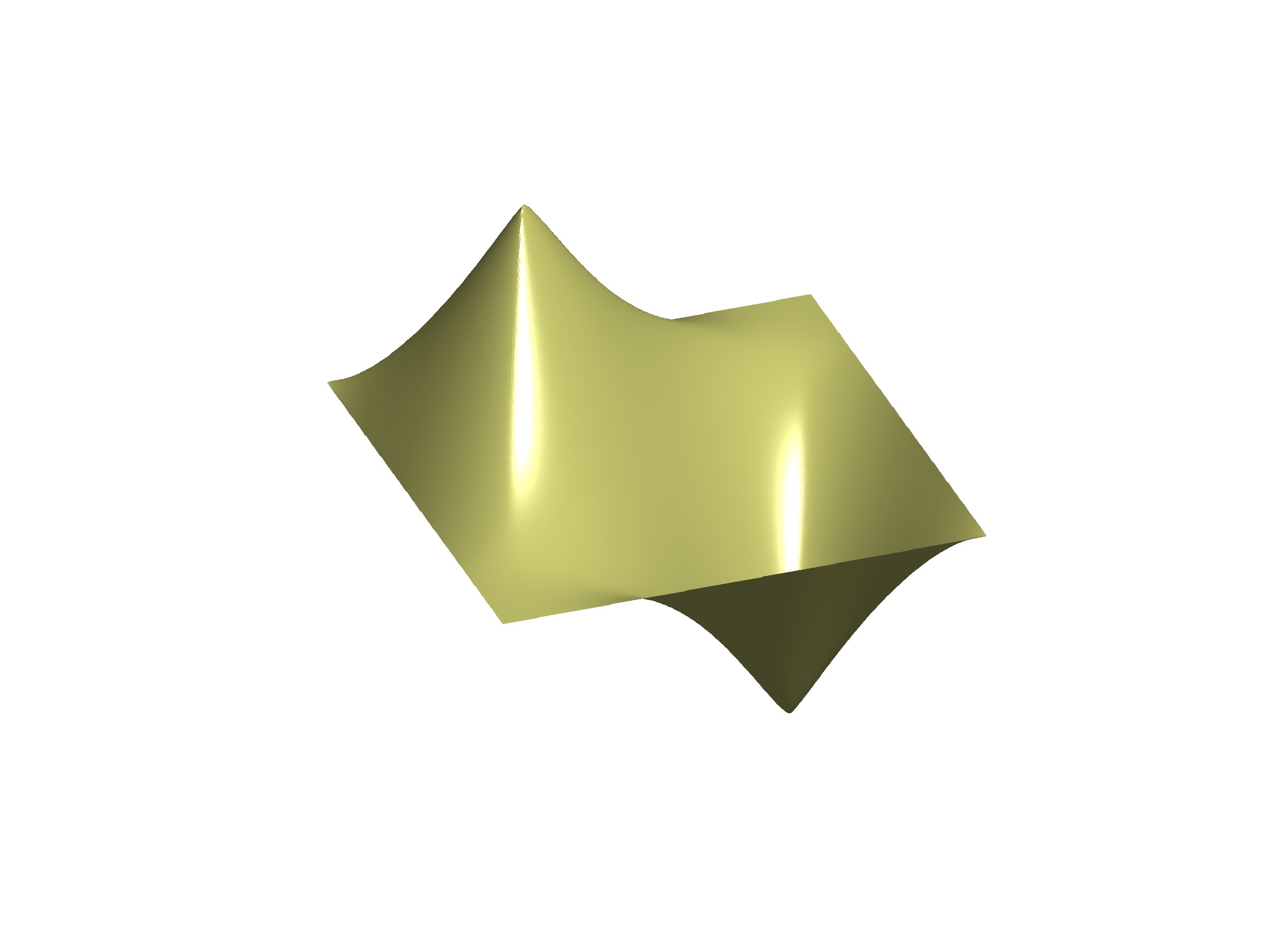}}
  \end{picture}
  \caption{The numerically computed second eigenfunction $u_2$ for the square.}
  \label{fig:square_u2}
\end{figure}

Let
\begin{displaymath}
  \Omega=\left\{(x_1,x_2)\in\R^2\,\big|\, x_1,x_2\in(0,2)\right\}.
\end{displaymath}
This domain was discretized using 83,968 triangles. Figures
\ref{fig:square_u1} and \ref{fig:square_u2} show the eigenfunctions
$u_1$ and $u_2$ computed for several values of $p$,
respectively. Table \ref{tab:square_lambda} lists the corresponding
values of $\lambda_1$ and $\lambda_2$.

\begin{table}[b]
  \centering\vspace{1ex}
  \begin{tabular}[t]{|c|c|c|c|c|c|c|c|c|}
    \cline{1-4}\cline{6-9}
    $p$ & $\lambda_1$ & $\lambda_2$ & $\lambda_{\mathcal{S}_2}$ & \qquad\qquad & $p$ & $\lambda_1$ & $\lambda_2$ & $\lambda_{\mathcal{S}_1}$ \\\cline{1-4}\cline{6-9}
    \multicolumn{9}{c}{\vspace{\tableheadsep}}\\\cline{1-4}\cline{6-9}
    1.1 & 2.3649 & 3.7586 & 3.8702 & & 2.0 & 4.9349 & 12.338 & 12.338 \\\cline{1-4}\cline{6-9}
    1.2 & 2.6934 & 4.5012 & 4.6179 & & 2.1 & 5.2139 & 13.684 & 13.744 \\\cline{1-4}\cline{6-9}
    1.3 & 2.9986 & 5.2500 & 5.3715 & & 2.2 & 5.4952 & 15.144 & 15.282 \\\cline{1-4}\cline{6-9}
    1.4 & 3.2834 & 6.0385 & 6.1621 & & 2.3 & 5.7791 & 16.725 & 16.961 \\\cline{1-4}\cline{6-9}
    1.5 & 3.5611 & 6.8835 & 7.0053 & & 2.4 & 6.0658 & 18.438 & 18.797 \\\cline{1-4}\cline{6-9}
    1.6 & 3.8356 & 7.7971 & 7.9118 & & 2.5 & 6.3552 & 20.293 & 20.802 \\\cline{1-4}\cline{6-9}
    1.7 & 4.1092 & 8.7897 & 8.8903 & & 3.0 & 7.8452 & 32.107 & 33.956 \\\cline{1-4}\cline{6-9}
    1.8 & 4.3830 & 9.8708 & 9.9490 & & 4.0 & 11.038 & 74.757 & 85.447 \\\cline{1-4}\cline{6-9}
    1.9 & 4.6581 & 11.050 & 11.095 & & 5.0 & 14.497 & 163.59 & 205.08 \\\cline{1-4}\cline{6-9}
    2.0 & 4.9349 & 12.338 & 12.338 & & 6.0 & 18.194 & 343.77 & 477.60 \\\cline{1-4}\cline{6-9}
    \multicolumn{4}{c}{} & & 8.0 & 26.221 & 1,402.1 & 2,443.4 \\\cline{6-9}
    \multicolumn{4}{c}{} & & 10.0 & 34.990 & 5,339.0 & 11,888 \\\cline{6-9}
  \end{tabular}\vspace{2ex}
  \caption{Eigenvalues for the square.}
  \label{tab:square_lambda}
\end{table}

Various choices of the intermediate path point $e_\mathrm{M}$ were
used to compute $u_2$. Only in case $p=2$ different choices of
$e_\mathrm{M}$ caused CMPA to converged to different functions
$u_2$. Since for the square the eigenspace corresponding to the second
eigenvalue of the Laplace operator is two-dimensional, CMPA converges
to some member of this eigenspace depending on the shape of the
initial path. Figure \ref{fig:square_u2} shows one such
eigenfunction. However, even for $p=2$ this has no influence on the
computed value of $\lambda_2$.

We say that a function has symmetry $\mathcal{S}_1$ (odd symmetry
about $x_1=1$ and even symmetry about $x_2=1$) if it belongs to
\begin{displaymath}
  \mathcal{S}_1 := \{u:\Omega\to\R\,|\,u(x_1,x_2)=-u(2-x_1,x_2), u(x_1,x_2)=u(x_1,2-x_2)\},
\end{displaymath}
and symmetry $\mathcal{S}_2$ (odd symmetry about $x_1=x_2$ and even
symmetry about $x_1=2-x_2$) if it belongs to
\begin{displaymath}
  \mathcal{S}_2 := \{u:\Omega\to\R\,|\,u(x_1,x_2)=-u(x_2,x_1),u(x_1,x_2)=u(2-x_2,2-x_1)\}.
\end{displaymath}
As it was observed in \cite{YaoZhou1}, $u_2$ changes its symmetry at
$p=2$ from $\mathcal{S}_1$ for $p<2$ to $\mathcal{S}_2$ for $p>2$. Let
$\lambda_{\mathcal{S}_i}$ denote the smallest eigenvalue with an
eigenfunction belonging to $\mathcal{S}_i$ where $i\in\{1,2\}$. The
values of $\lambda_{\mathcal{S}_i}$ can be computed using CDM on
$\Omega$ with additional boundary conditions $u(1,x_2)=0$ for
$x_2\in(0,2)$ or $u(x,x)=0$ for $x\in(0,2)$, respectively, or as the
first eigenvalue on the half-domain $\Omega^\mathrm{half}_i$, where
\begin{align}
  \Omega^\mathrm{half}_1 &:= \left\{(x_1,x_2)\in\R^2\,\big|\,x_1\in(0,1), x_2\in(0,2)\right\},\notag \\
  \Omega^\mathrm{half}_2 &:= \left\{(x_1,x_2)\in\R^2\,\big|\,x_1\in(0,2), x_2\in(0,x_1)\right\}.\notag
\end{align}

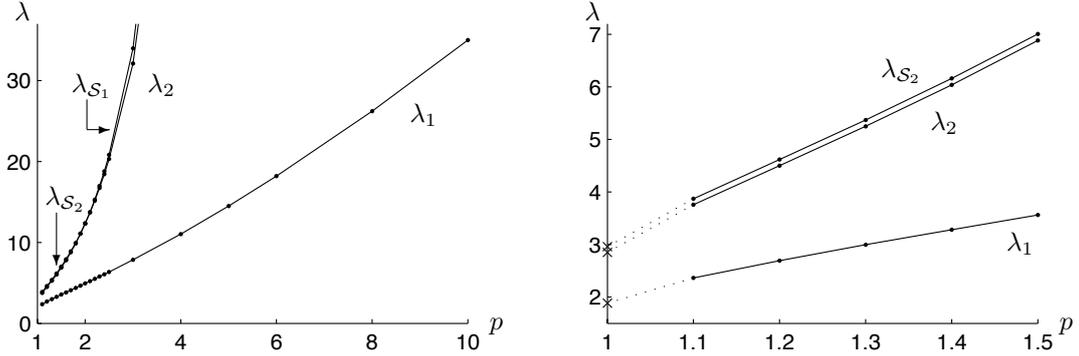
\begin{figure}[h]
  \centering
  \setlength{\unitlength}{1mm}
  \begin{picture}(150,50)(0,2)
    \put(0,0){\includegraphics[width=7.3cm]{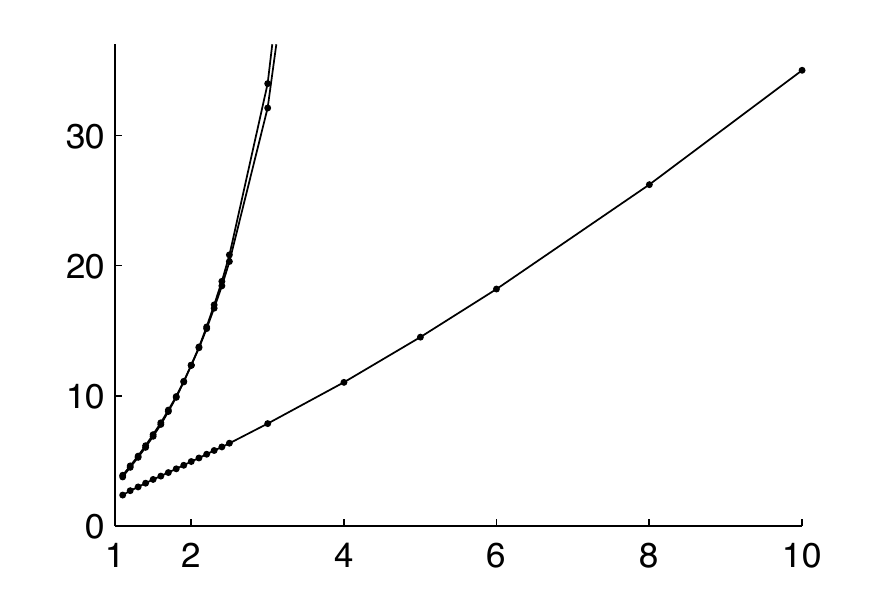}}
    \put(6.5,45.5){$\lambda$}
    \put(69,4.5){$p$}
    \put(58.5,32){$\lambda_1$}
    \put(24,36){$\lambda_2$}
    \put(10.5,21){$\lambda_{\mathcal{S}_2}$}
    \put(12,20){\vector(0,-1){7}}
    \put(14,36){$\lambda_{\mathcal{S}_1}$}
    \put(16,31){\line(0,1){4}}
    \put(16,31){\vector(1,0){3}}
    \put(75,0){\includegraphics[width=7.3cm]{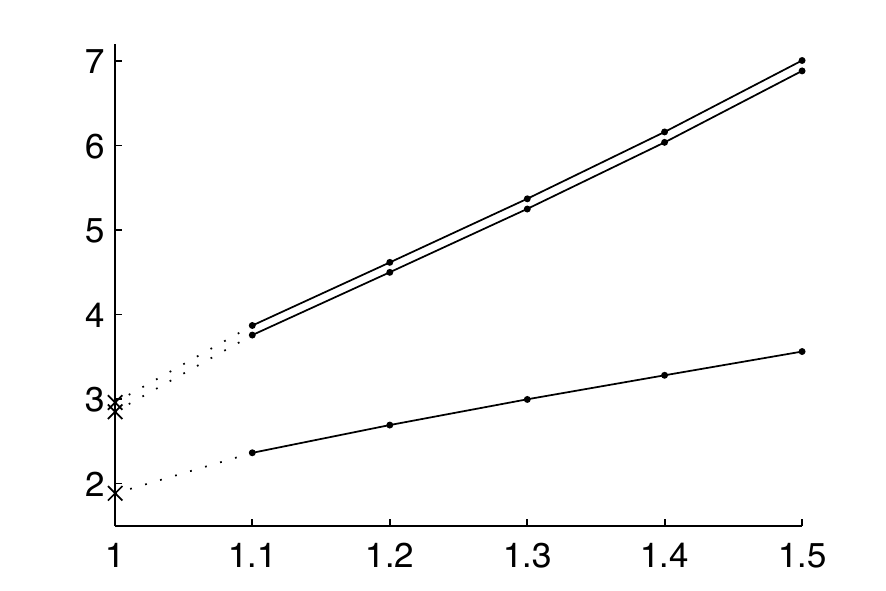}}
    \put(81.5,45.5){$\lambda$}
    \put(144,4.5){$p$}
    \put(137,15){$\lambda_1$}
    \put(127,31){$\lambda_2$}
    \put(120.5,38){$\lambda_{\mathcal{S}_2}$}
  \end{picture}  
  \caption{Dependence of the numerically computed eigenvalues for the
    square on $p$. The three cross symbols in the graph on the right
    mark the values of $h_1$ for $\Omega$, $\Omega^\mathrm{half}_1$,
    and $\Omega^\mathrm{half}_2$.}
  \label{fig:square_lambda_plot1}
\end{figure}

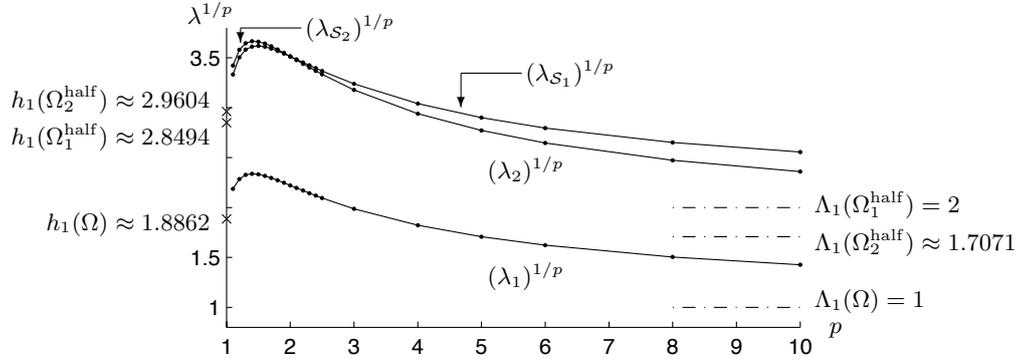
\begin{figure}[h]
  \centering
  \setlength{\unitlength}{1mm}
  \begin{picture}(130,50)(-15,2)\small
    \put(0,0){\includegraphics[width=9.73cm]{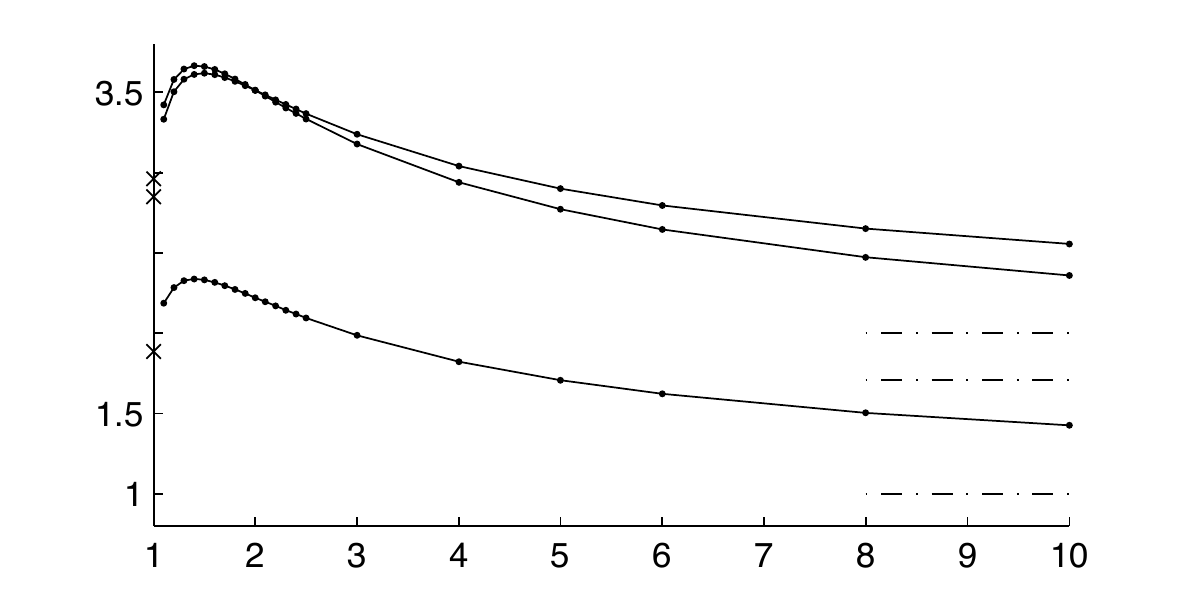}}
    \put(7.5,45.5){$\lambda^{1/p}$}
    \put(92,4.5){$p$}
    \put(47,11){$(\lambda_1)^{1/p}$}
    \put(47,25){$(\lambda_2)^{1/p}$}
    \put(52,38){$(\lambda_{\mathcal{S}_1})^{1/p}$}
    \put(43.5,39){\vector(0,-1){4.5}}
    \put(43.5,39){\line(1,0){8}}
    \put(23,44){$(\lambda_{\mathcal{S}_2})^{1/p}$}
    \put(14.5,45){\vector(0,-1){2}}
    \put(14.5,45){\line(1,0){8}}
    \put(-14.5,17){\makebox(25,4)[r]{$h_1(\Omega)\approx 1.8862$}}
    \put(-14.5,28.5){\makebox(25,4)[r]{$h_1(\Omega^\mathrm{half}_1)\approx 2.8494$}}
    \put(-14.5,33.5){\makebox(25,4)[r]{$h_1(\Omega^\mathrm{half}_2)\approx 2.9604$}}
    \put(90,8){$\Lambda_1(\Omega)=1$}
    \put(90,15.5){$\Lambda_1(\Omega^\mathrm{half}_2)\approx 1.7071$}
    \put(90,20.5){$\Lambda_1(\Omega^\mathrm{half}_1)=2$}
  \end{picture}  
  \caption{Dependence of the numerically computed eigenvalues for the square
    raised to $1/p$ on $p$.}
  \label{fig:square_lambda_plot2}
\end{figure}

Our numerical observations regarding these eigenvalues and $\lambda_2$
are summarized in Table \ref{tab:square_observations}(a) and the
computed values are listed in Table \ref{tab:square_lambda}. We stress
that $\lambda_2$ was computed with no a priori assumptions on
symmetry. The dependence of the eigenvalues $\lambda_1$, $\lambda_2$,
$\lambda_{\mathcal{S}_1}$ and $\lambda_{\mathcal{S}_2}$ on $p$ is
further plotted in Figures \ref{fig:square_lambda_plot1} and
\ref{fig:square_lambda_plot2}.

\begin{table}[h]
  \begin{center}\renewcommand{\arraystretch}{1.3}
    \begin{tabular}{r|c|c|c|cr|c|c|c|c|}
      \cline{3-4}\cline{8-10}
      \multicolumn{1}{c}{(a)} & \multicolumn{1}{c|}{} & $\lambda_{\mathcal{S}_1}$ & $\lambda_{\mathcal{S}_2}$ &
      \rule{2mm}{0pt} &
      \multicolumn{1}{c}{(b)} & \multicolumn{1}{c|}{} & $\lambda_1$ &
      $\lambda_{\mathcal{S}_1}$ & $\lambda_{\mathcal{S}_2}$\\\cline{2-4}\cline{7-10}
      & $p<2$ & $=\lambda_2$ & $>\lambda_2$ &&& $\lim_{p\to1+}\lambda$ &
      $1+\frac{1}{2}\sqrt{\pi}$ & $\frac{4-\pi}{3-\sqrt{1+2\pi}}$ & $1+\frac{1}{2}(\sqrt{2}+\sqrt{2\pi})$
      \\\cline{2-4}\cline{7-10}
      & $2<p$ & $>\lambda_2$ & $=\lambda_2$ &&& $\lim_{p\to\infty}\lambda^{1/p}$ &
      1 & 2 & $1+\frac{1}{2}\sqrt{2}$ \\\cline{2-4}\cline{7-10}
    \end{tabular}\renewcommand{\arraystretch}{1}\vspace{2ex}
  \end{center}
  \caption{The smallest eigenvalues $\lambda_{\mathcal{S}_1}$ and $\lambda_{\mathcal{S}_2}$
    under symmetry assumptions for the square.
    (a) Numerical comparison with $\lambda_2$.
    (b) Asymptotic behavior: the first row shows values of $h_1$, the second row values of $\Lambda_1$
    for $\Omega$, $\Omega^\mathrm{half}_1$, and $\Omega^\mathrm{half}_2$, respectively.}
  \label{tab:square_observations}
\end{table}

These figures and Table \ref{tab:square_observations}(b) also explain
the asymptotic behavior as $p\to 1$ and $p\to\infty$. While the table
shows the limit values as given by the theory, the graphs indicate
convergence to these values (at least for $p\to 1$; for $p\to\infty$
it seems a larger range of $p$ would be needed). The Cheeger constants
$h_1$ shown in the first row of Table \ref{tab:square_observations}(b)
have been computed according to \cite{KawohlFridman},
\cite{LachandRobertKawohl} by
\begin{equation}
  \label{eq:square_cheeger}
  h_1\big((0,a)\times(0,b)\big)=\frac{4-\pi}{a+b-\sqrt{(a-b)^2+\pi ab}} \qquad\text{for }a,b>0.
\end{equation}
The evaluation of $\Lambda_1$ in the second row is straightforward.

\subsection{Rectangle}\label{rectangle}
\begin{figure}[h]
  \centering
  \setlength{\unitlength}{1mm}
  \begin{picture}(150,33)
    \put(0,0){\reutwographs{3.6}{out_u2_p_3-6}}
    \put(37.5,0){\reutwographs{3.7}{out_u2_p_3-7}}
    \put(75,0){\reutwographs{3.8}{out_u2_p_3-8}}
    \put(112.5,0){\reutwographs{8.0}{out_u2_p_8-0}}
  \end{picture}  
  \caption{The numerically computed second eigenfunction $u_2$ for the rectangle.}
  \label{fig:rect_u2}
\end{figure}
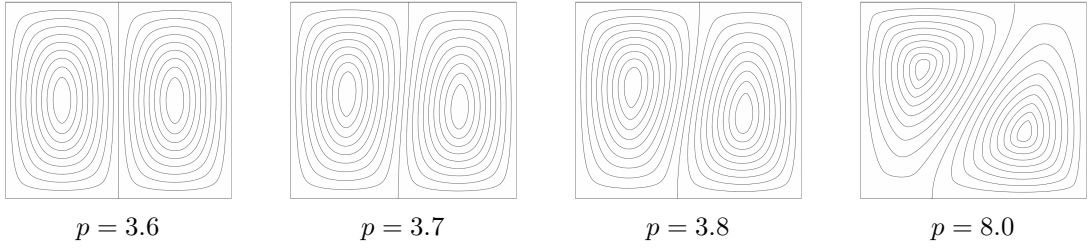

Let
\begin{displaymath}
  \Omega=\left\{(x_1,x_2)\in\R^2\,\big|\, x_1\in(0,2), x_2\in(0,1.75)\right\}.
\end{displaymath}
This domain was discretized using 77,312 triangles. The shape of the
first eigenfunction $u_1$ and the graph of the first eigenvalue
$\lambda_1(\Omega;p)$ are similar to those for the square. However,
the symmetry properties of the second eigenfunction $u_2$ are
different: According to our numerical observations, for $p\leq 3.6$
the eigenfunction $u_2$ preserves an odd symmetry about $x_1=1$ and an
even symmetry about $x_2=0.875$ (which we call $\mathcal{S}_1$ as in
the case of the square). For $p\geq 3.7$ this symmetry is lost and
$u_2$ maintains an odd symmetry with respect to $(1,0.875)$, the
center of $\Omega$. The contour lines of $u_2$ for several values of
$p$ are shown in Fig.~\ref{fig:rect_u2}.

\begin{figure}[h]
  \centering
  \setlength{\unitlength}{1mm}
  \begin{picture}(140,42)\small
    \put(0,0){\includegraphics[width=7.5cm]{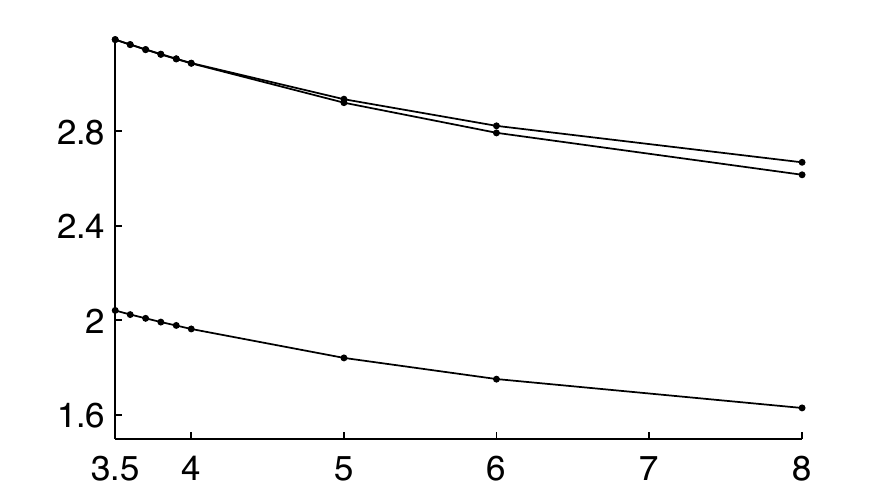}}
    \put(4.5,39){$\lambda^{1/p}$}
    \put(52,10){$(\lambda_1)^{1/p}$}
    \put(52,23.5){$(\lambda_2)^{1/p}$}
    \put(52,31.5){$(\lambda_{\mathcal{S}_1})^{1/p}$}
    \put(71,3.5){$p$}
    \put(80,2){\normalsize
      \begin{tabular}[b]{|c|c|c|c|}
        \multicolumn{4}{c}{$\Omega=(0,2)\times(0,1.75)$} \\
        \multicolumn{4}{c}{\vspace{\tableheadsep}}\\\hline
        $p$ & $\lambda_1$ & $\lambda_2$ & $\lambda_{\mathcal{S}_1}$ \\ \hline
        \multicolumn{4}{c}{\vspace{\tableheadsep}}\\\hline
        3.5 & 12.166 & 57.874 & 57.874 \\\hline
        3.6 & 12.680 & 63.436 & 63.436 \\\hline
        3.7 & 13.206 & 69.491 & 69.492 \\\hline
        3.8 & 13.743 & 76.059 & 76.083 \\\hline
        3.9 & 14.292 & 83.175 & 83.256 \\\hline
        4.0 & 14.853 & 90.881 & 91.059 \\\hline
        8.0 & 49.531 & 2,192.9 & 2,574.6 \\\hline
      \end{tabular}
    }
  \end{picture}
  \caption{Comparison of the numerically computed second eigenvalue
    $\lambda_2$ and the smallest eigenvalue $\lambda_{\mathcal{S}_1}$
    under the symmetry $\mathcal{S}_1$ for the rectangle
    $\Omega$ ($\lambda_1$ is shown for reference).}
  \label{fig:rect_lambda2}
\end{figure}
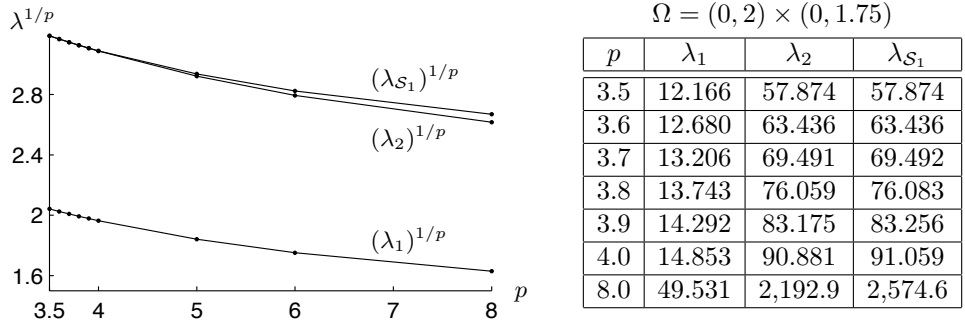

\begin{table}[h]
  \centering
  \begin{tabular}{c|c|c|c|c|cc|c|c|c|c|}
    \multicolumn{1}{c}{(a)} & \multicolumn{4}{c}{Rectangle $(0,2)\times(0,1.9)$} &
    \rule{8mm}{0pt} & \multicolumn{1}{c}{(b)} &
    \multicolumn{4}{c}{Rectangle $(0,2)\times(0,1.6)$}\\
    \multicolumn{11}{c}{\vspace{\tableheadsep}}\\\cline{2-5}\cline{8-11}
    & $p$ & $\lambda_1$ & $\lambda_2$ & $\lambda_{\mathcal{S}_1}$ &&
    & $p$ & $\lambda_1$ & $\lambda_2$ & $\lambda_{\mathcal{S}_1}$ \\\cline{2-5}\cline{8-11}
    \multicolumn{11}{c}{\vspace{\tableheadsep}}\\\cline{2-5}\cline{8-11}
    & 2.44 & 6.5926 & 20.0177 & 20.0177 &&& 5.6 & 36.077 & 383.4648 & 383.4648
    \\\cline{2-5}\cline{8-11}
    & 2.46 & 6.6579 & 20.4281 & 20.4281 &&& 5.8 & 38.898 & 453.2332 & 453.2333
    \\\cline{2-5}\cline{8-11}
    & 2.48 & 6.7234 & 20.8451 & 20.8457 &&& 6.0 & 41.892 & 535.2007 & 535.2009
    \\\cline{2-5}\cline{8-11}
    & 2.50 & 6.7891 & 21.2683 & 21.2708 &&& 6.2 & 45.068 & 631.4438 & 631.4478
    \\\cline{2-5}\cline{8-11}
    & 2.60 & 7.1205 & 23.4816 & 23.5124 &&& 6.4 & 48.437 & 744.1846 & 744.4026
    \\\cline{2-5}\cline{8-11}
  \end{tabular}\vspace{2ex}
  \caption{Comparison of the numerically computed second eigenvalue
    $\lambda_2$ and the smallest eigenvalue $\lambda_{\mathcal{S}_1}$
    under the symmetry $\mathcal{S}_1$ for other rectangles.}
  \label{tab:rect_more_lambda2}
\end{table}

For $p\geq 3.7$ the smallest eigenvalue $\lambda_{\mathcal{S}_1}$
corresponding to an eigenfunction with symmetry $\mathcal{S}_1$ is
larger than $\lambda_2$ (cf.\ Fig.~\ref{fig:rect_lambda2}). This
eigenpair can be computed by CDM on $\Omega$ with an additional
boundary condition $u(1,x_2)=0$ for $x_2\in(0,1.75)$ or as the first
eigenpair on the half-rectangle
\begin{displaymath}
  \Omega^\mathrm{half}=\left\{(x_1,x_2)\in\R^2\,\big|\, x_1\in(0,1),
    x_2\in(0,1.75)\right\}.
\end{displaymath}

Our conjecture is that for a rectangle $R=(0,a)\times(0,b)$ with
$0<b<a$
there exists $p_0>2$ such that $u_2$ has two nodal domains which for
$p<p_0$ are rectangles with sides $a/2$ and $b$. For $p>p_0$ the nodal
domains are not rectangular and $u_2$ has only an odd symmetry with
respect to the center of $R$. According to our numerical observations,
$p_0$ gets larger the larger the ratio $a/b$: Besides $\Omega$ we ran
the computation for two other rectangles. For $R=(0,2)\times(0,1.9)$
the loss of symmetry $\mathcal{S}_1$ of $u_2$ is observable
approximately between $p=2.44$ and $p=2.48$ and for
$R=(0,2)\times(0,1.6)$ between $p=5.6$ and $p=6.0$ (cf.\
Table~\ref{tab:rect_more_lambda2}). As $p$ grows and crosses $p_0$,
the nodal line which is straight for $p<p_0$ gets distorted. This
distortion is faster for smaller ratios $a/b$ and slower for larger
ratios.

\subsection{Triangle with height 1}\label{triangle1}
\begin{figure}[p]
  \centering
  \setlength{\unitlength}{1mm}
  \begin{picture}(150,56)
    \put(0,0){\troneuonegraphs{1.1}{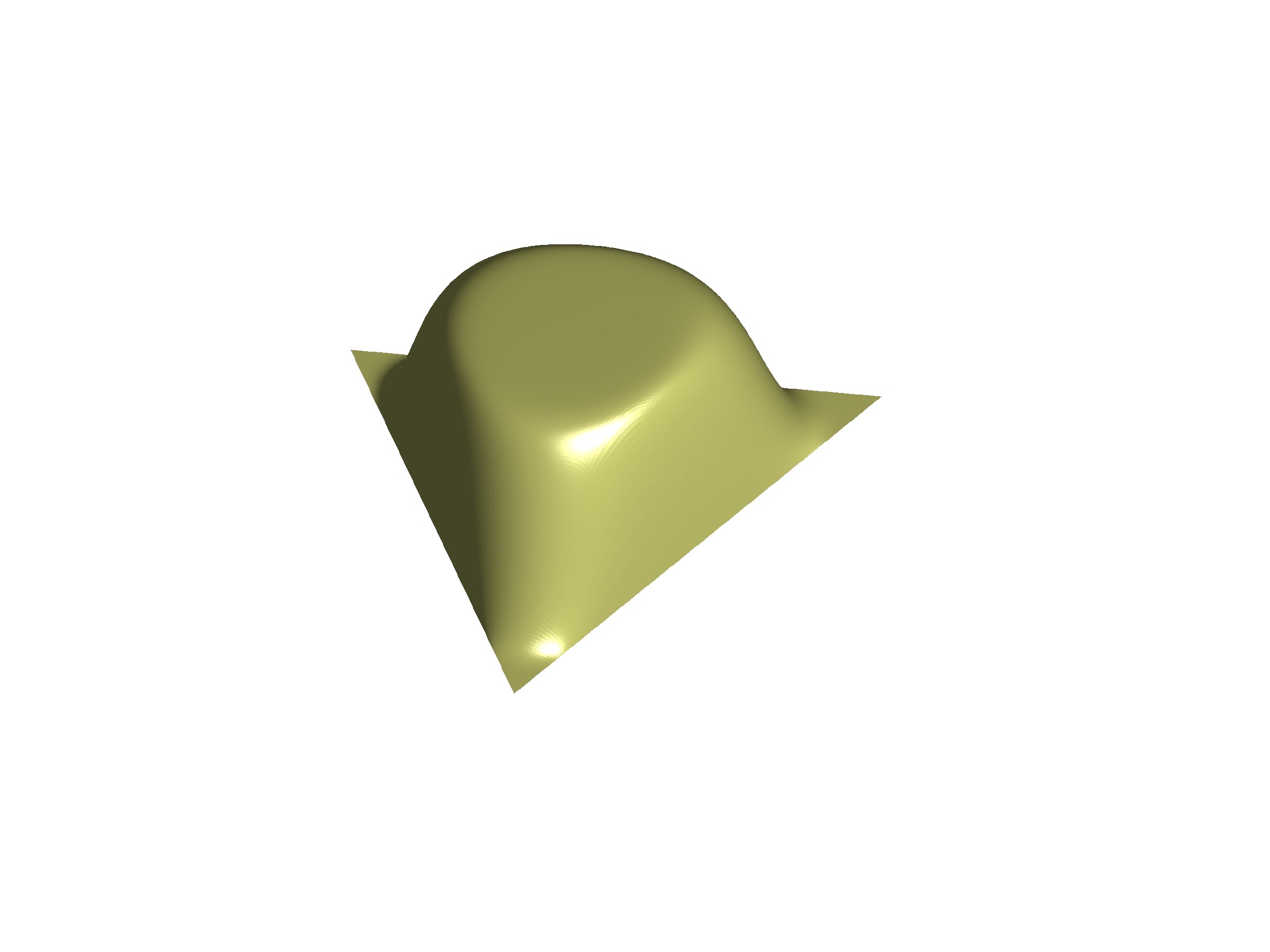}}
    \put(30,0){\troneuonegraphs{1.4}{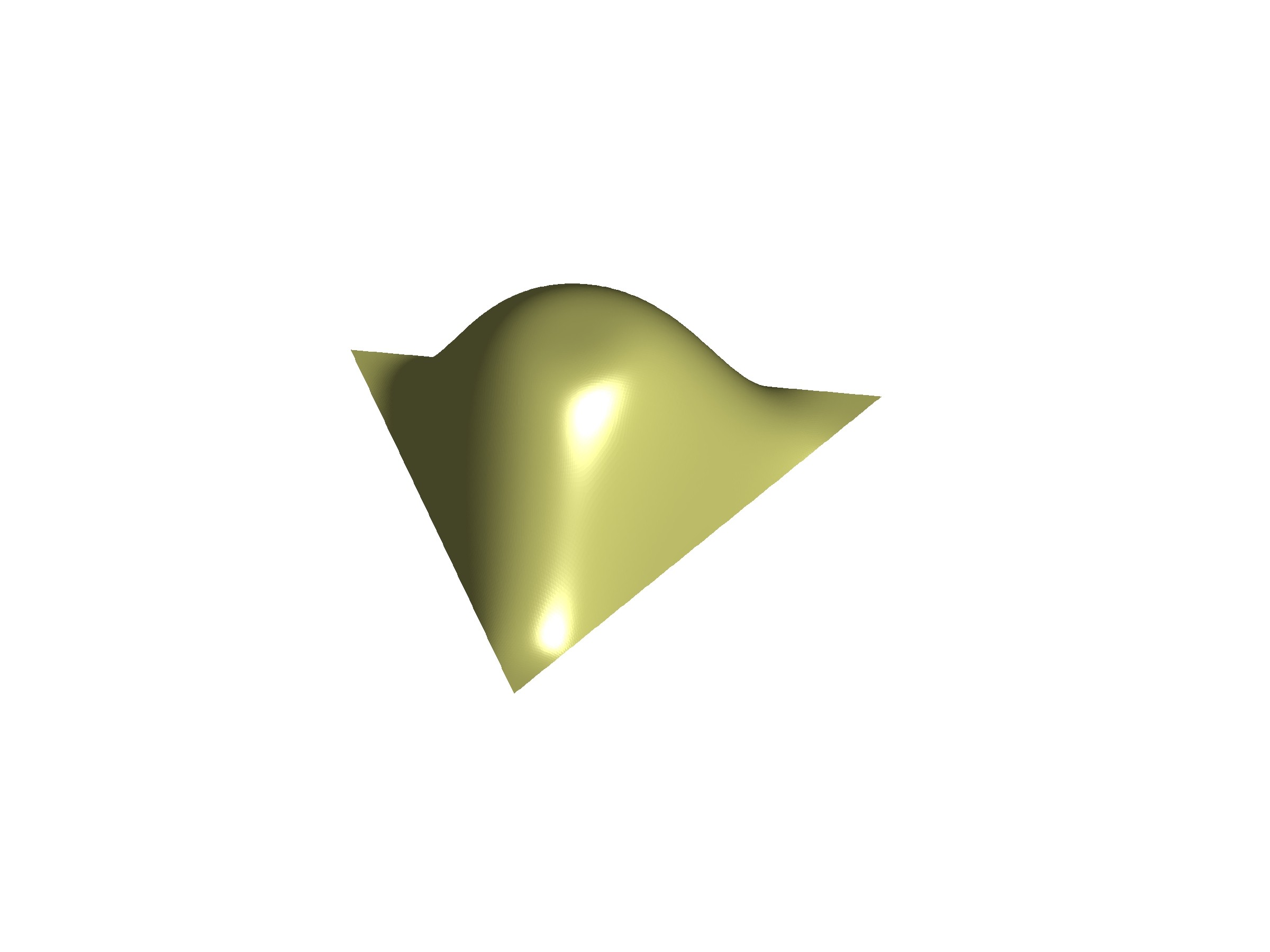}}
    \put(60,0){\troneuonegraphs{2.0}{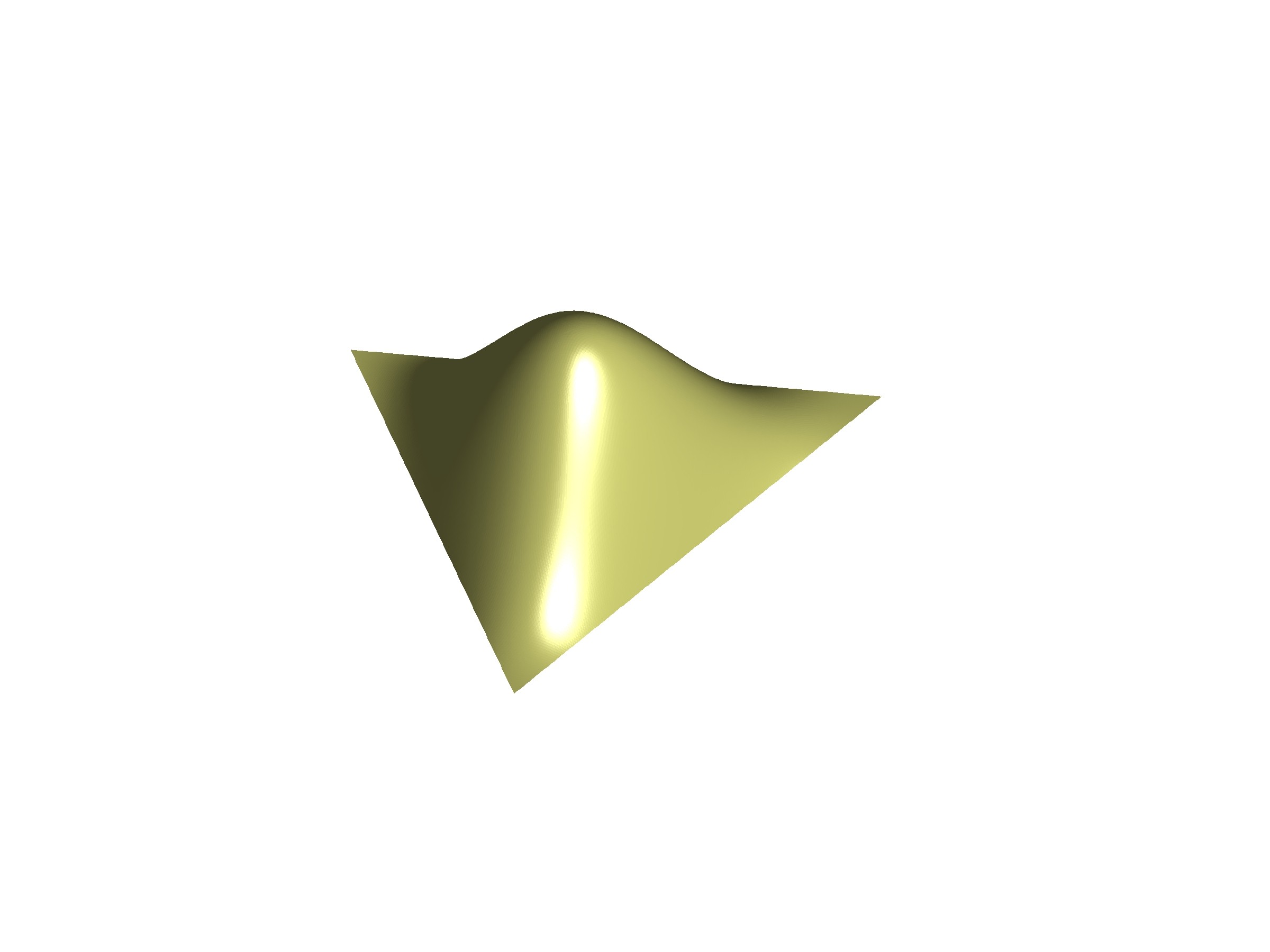}}
    \put(90,0){\troneuonegraphs{2.5}{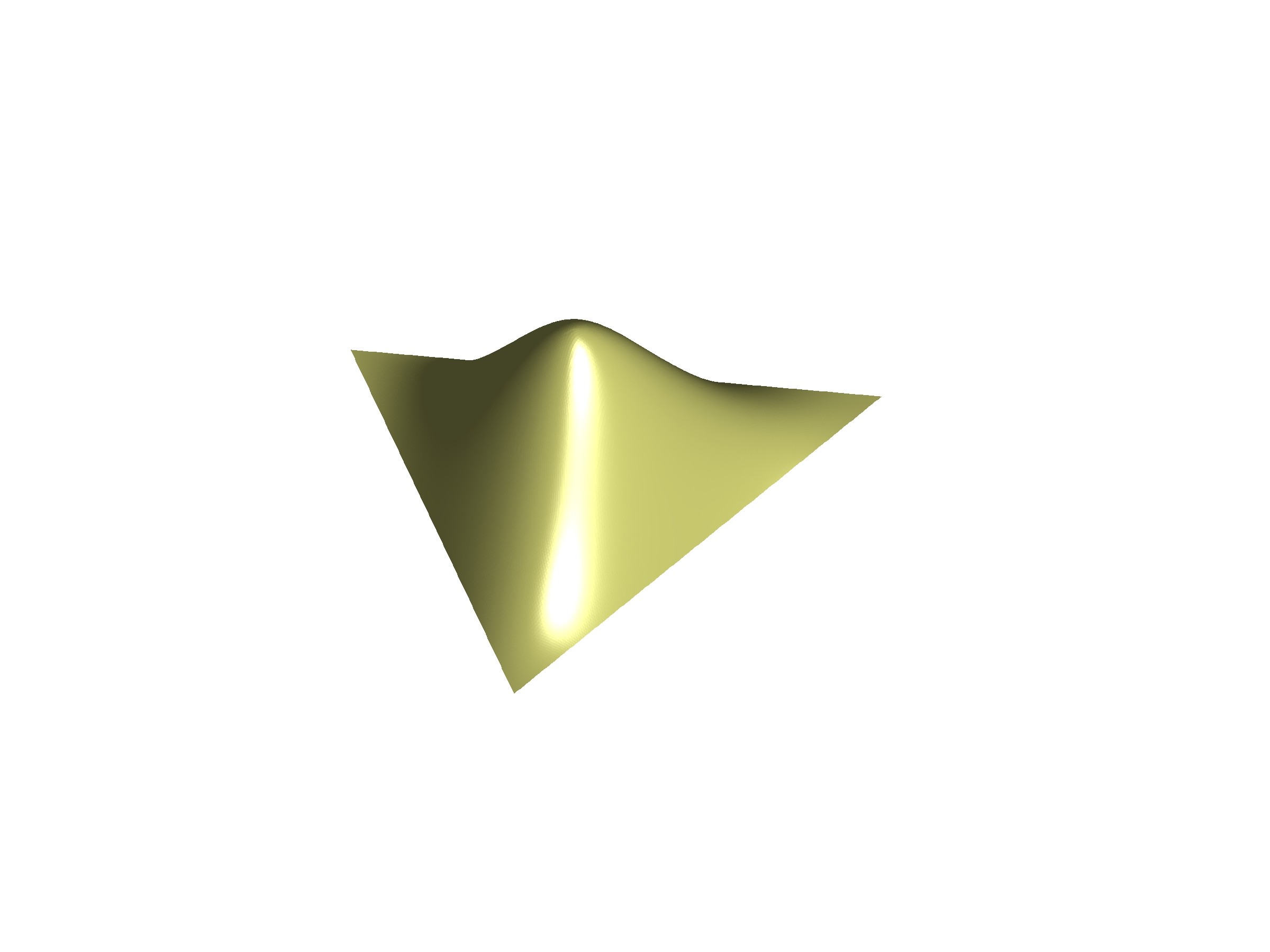}}
    \put(120,0){\troneuonegraphs{8.0}{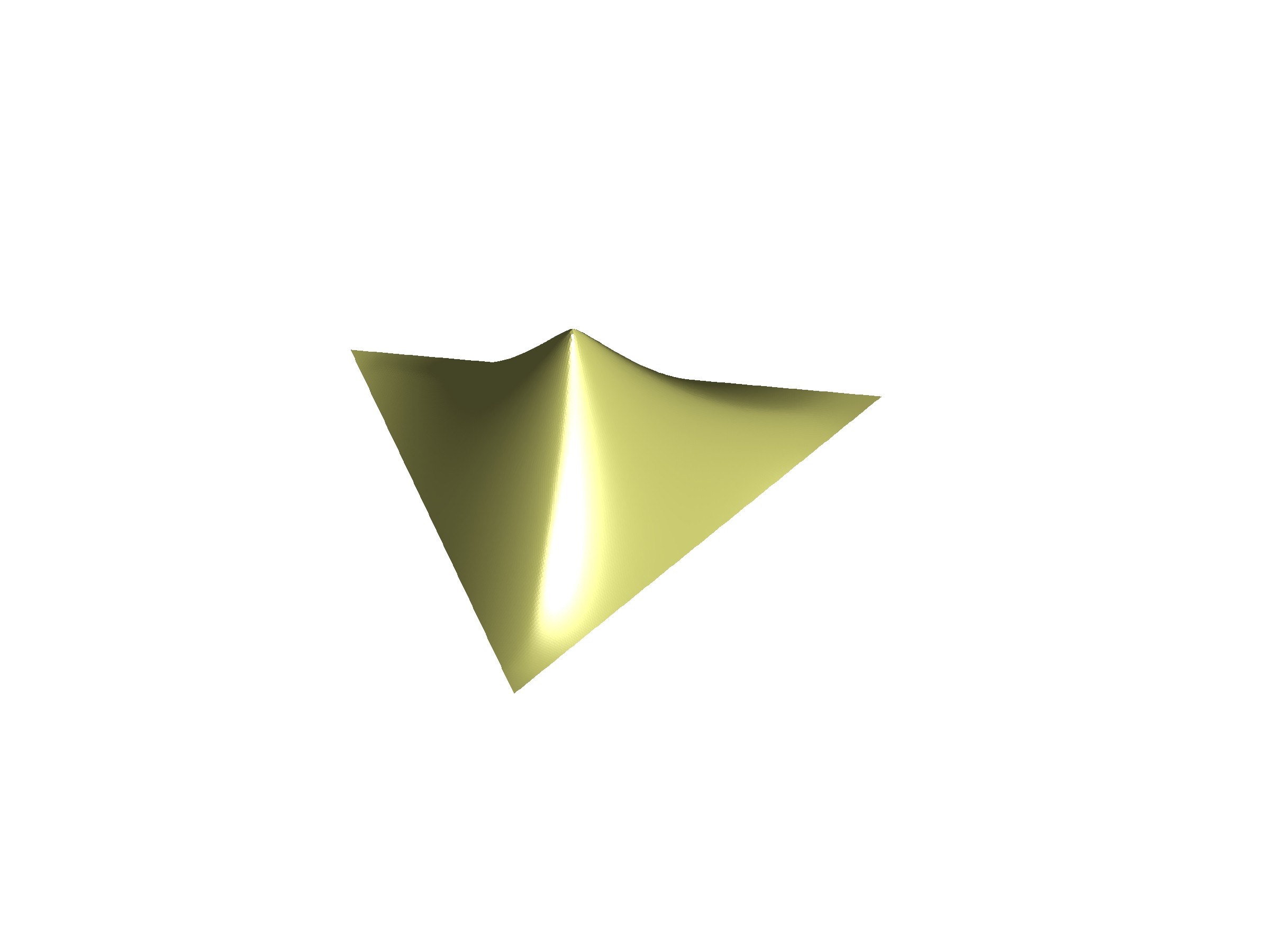}}
    \put(-.5,14){\rotatebox{90}{\scriptsize base}}
    \put(2.5,40){\rotatebox{-65}{\scriptsize base}}
  \end{picture}
  \caption{The numerically computed first eigenfunction $u_1$ for the triangle with height 1.}
  \label{fig:triangle1_u1}
\end{figure}

\begin{figure}[p]
  \centering
  \setlength{\unitlength}{1mm}
  \begin{picture}(150,54)
    \put(0,0){\troneutwographs{1.1}{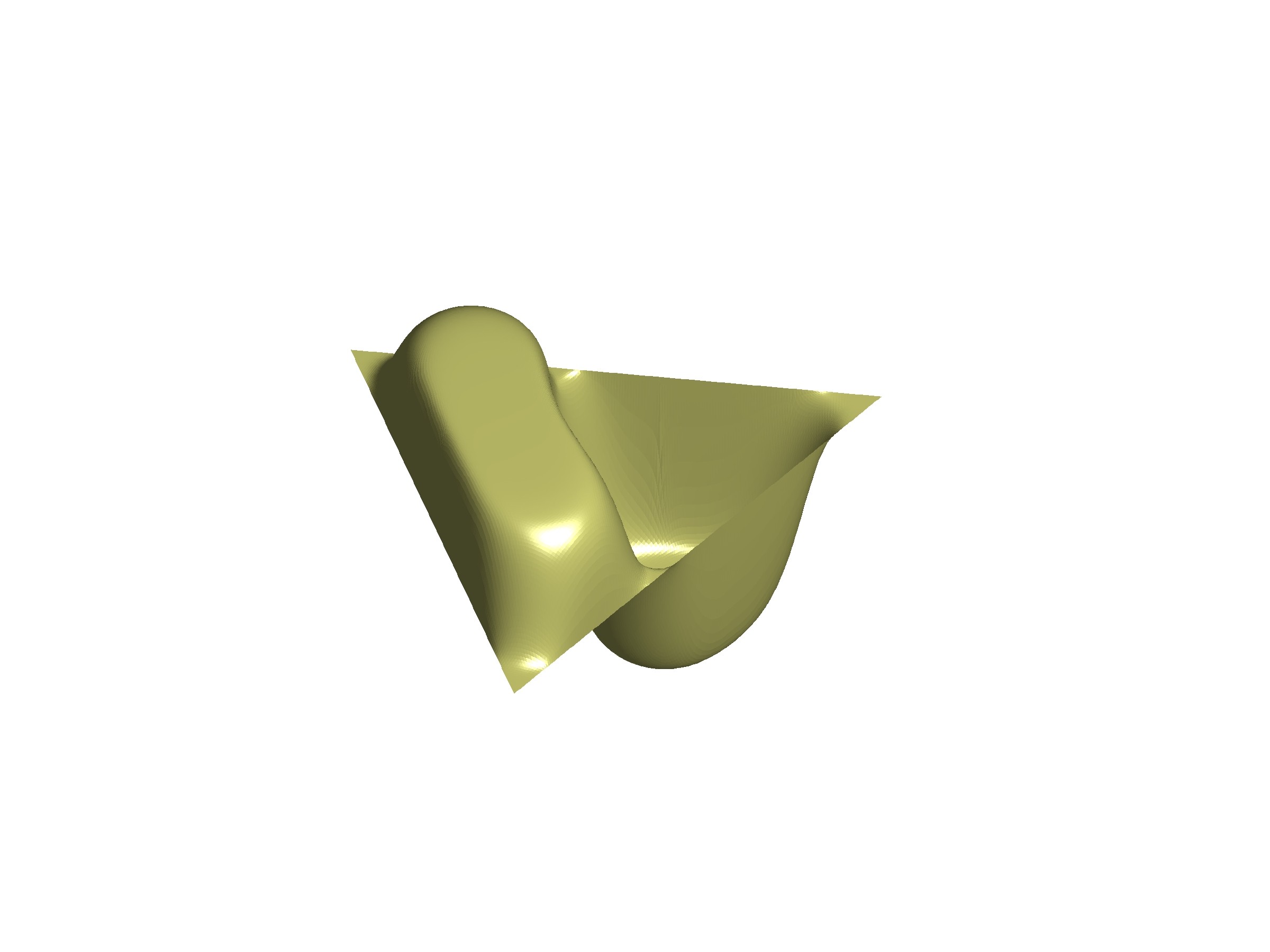}}
    \put(30,0){\troneutwographs{2.6}{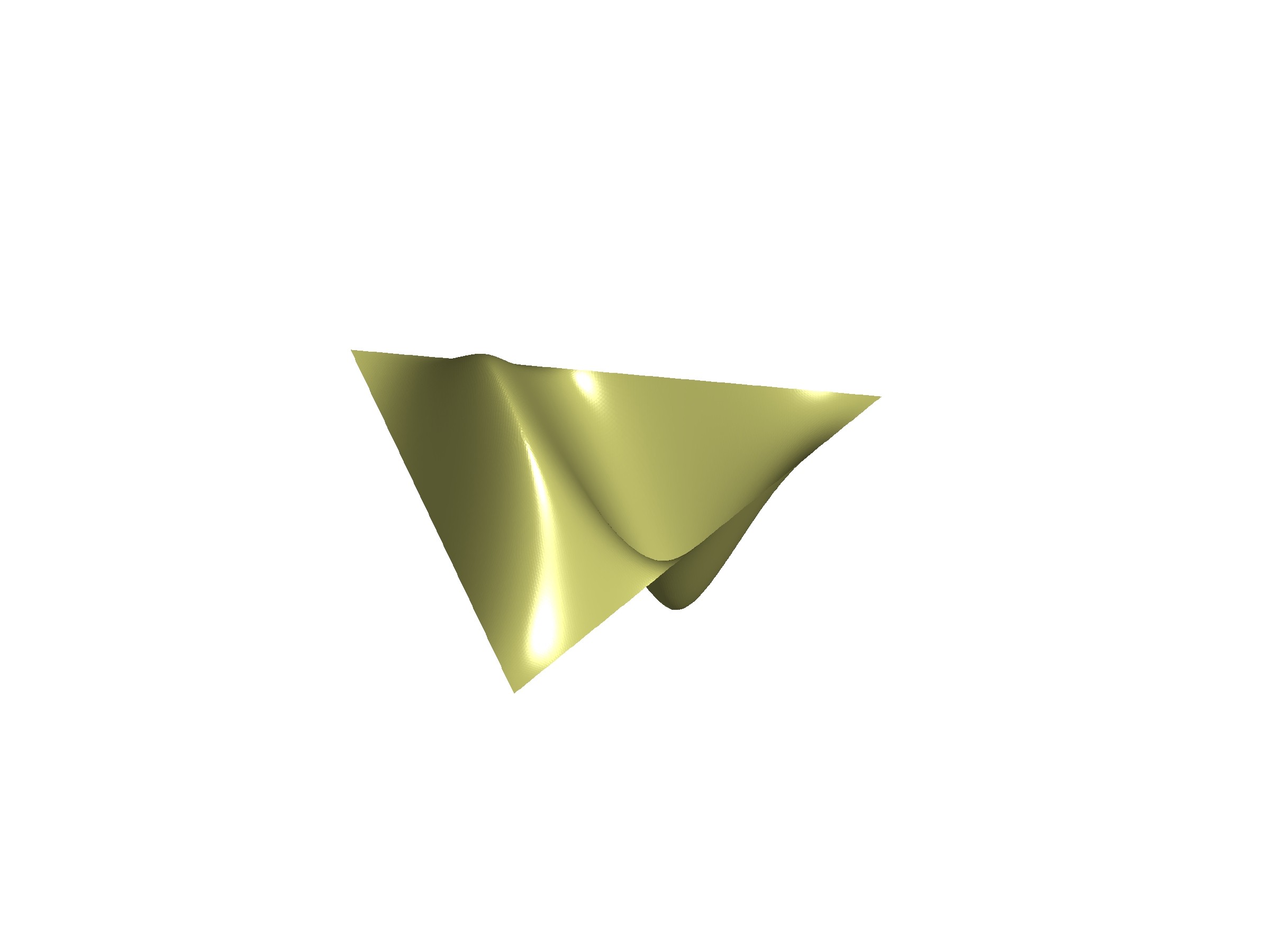}}
    \put(60,0){\troneutwographs{2.7}{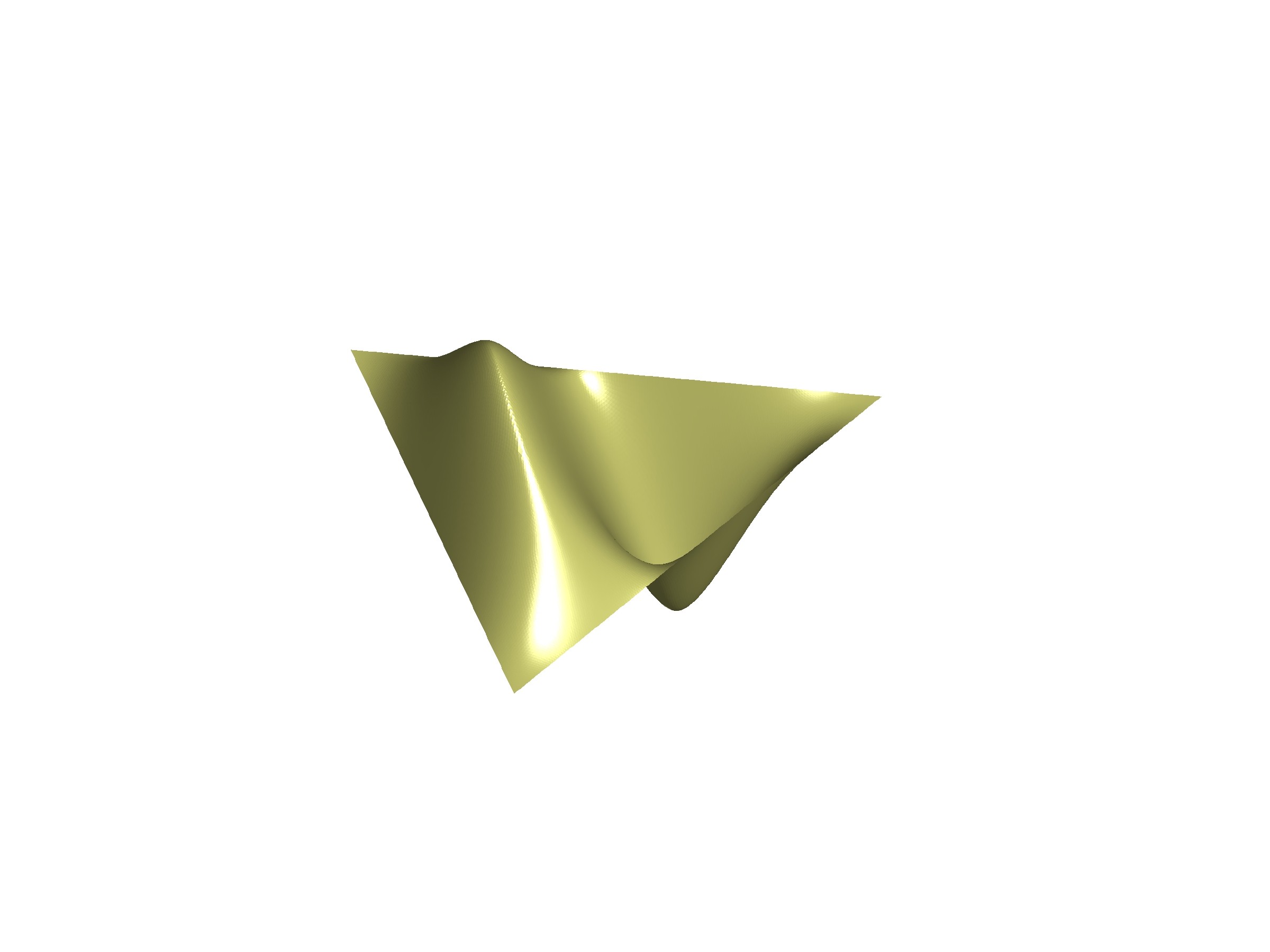}}
    \put(90,0){\troneutwographs{3.0}{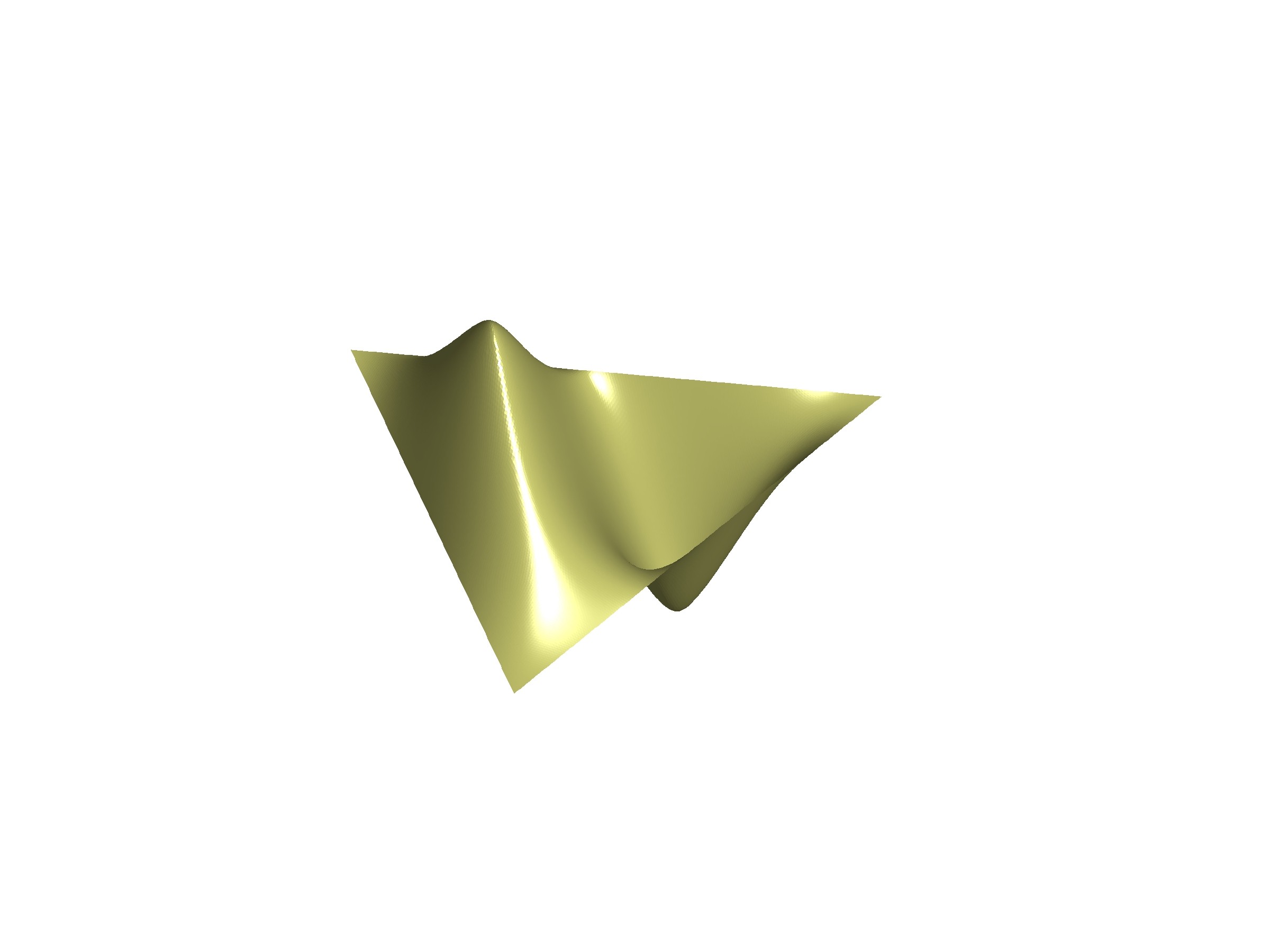}}
    \put(120,0){\troneutwographs{8.0}{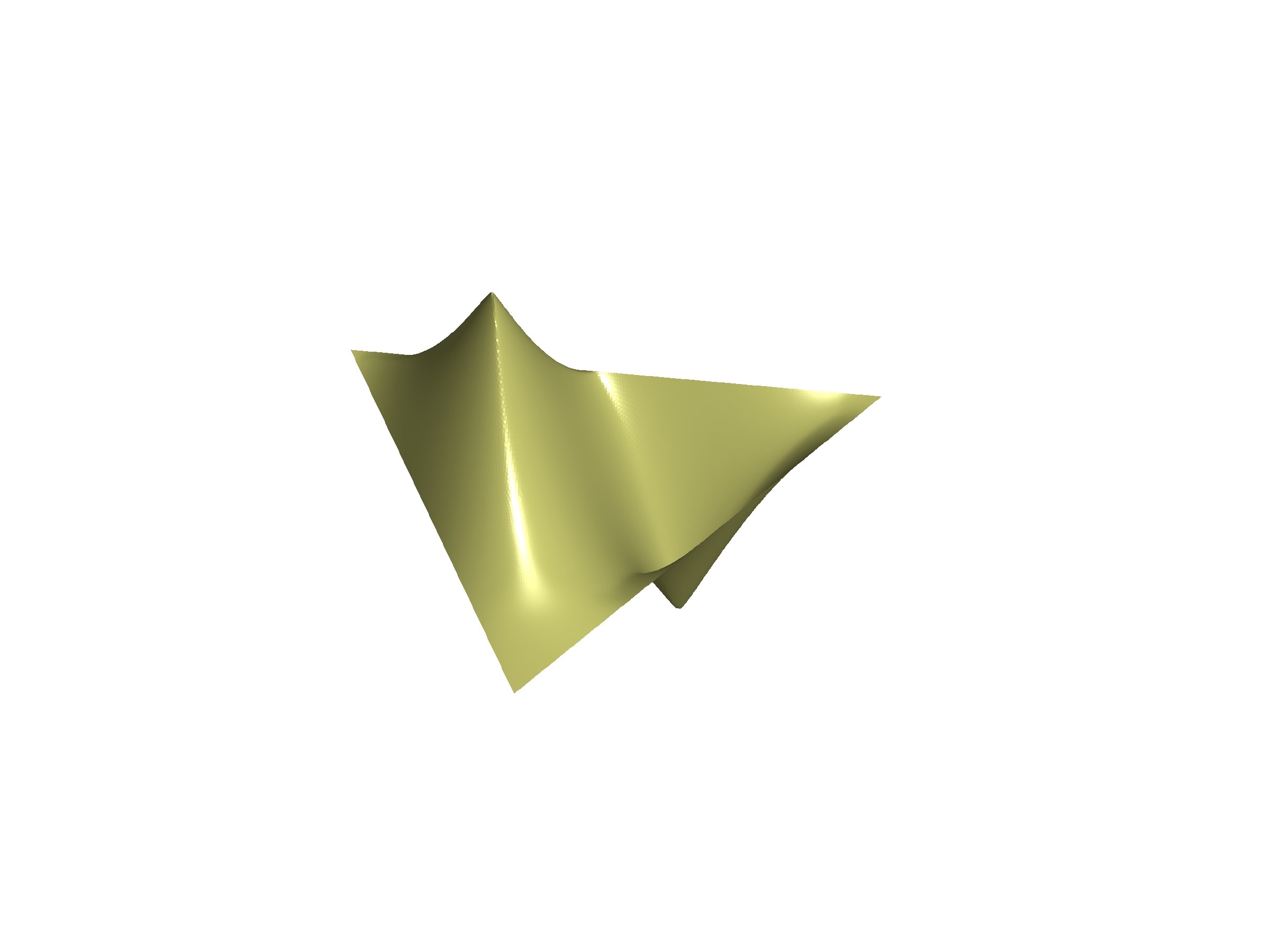}}
  \end{picture}
  \caption{The numerically computed second eigenfunction $u_2$ for the triangle with height 1.}
  \label{fig:triangle1_u2}
\end{figure}

Let
\begin{displaymath}\textstyle
  \Omega=\left\{(x_1,x_2)\in\R^2\,\big|\, x_1\in(0,1),|x_2|<\frac{1}{2}(1-x_1)\right\}
\end{displaymath}
be an isosceles triangle with base 1 and height 1. It was discretized
using 38,912 triangles. Figures \ref{fig:triangle1_u1} and
\ref{fig:triangle1_u2} show the eigenfunctions $u_1$ and $u_2$ for
several values of $p$, respectively. Table \ref{tab:triangle1_lambda}
lists the corresponding values of $\lambda_1$ and $\lambda_2$.

\begin{table}[b]
  \centering\vspace{1ex}
  \begin{tabular}[t]{|c|c|c|c|c|c|c|c|}
    \cline{1-3}\cline{5-8}
    $p$ & $\lambda_1$ & $\lambda_2$ & \qquad\qquad & $p$ & $\lambda_1$ & $\lambda_2$ & $\lambda_{2,\mathcal{S}_\mathrm{E}}$ \\\cline{1-3}\cline{5-8}
    \multicolumn{8}{c}{\vspace{\tableheadsep}}\\\cline{1-3}\cline{5-8}
    1.1 & 8.0143 & 12.188 & & 2.6 & 122.02 & 356.35 & 356.35 \\\cline{1-3}\cline{5-8}
    1.2 & 10.208 & 16.211 & & 2.7 & 142.81 & 435.98 & 435.99 \\\cline{1-3}\cline{5-8}
    1.3 & 12.673 & 21.009 & & 2.8 & 166.94 & 532.61 & 532.78 \\\cline{1-3}\cline{5-8}
    1.4 & 15.515 & 26.847 & & 2.9 & 194.90 & 649.76 & 650.31 \\\cline{1-3}\cline{5-8}
    1.5 & 18.822 & 33.998 & & 3.0 & 227.29 & 791.69 & 792.92 \\\cline{1-3}\cline{5-8}
    1.6 & 22.683 & 42.774 & & 3.5 & 483.05 & 2,093.5 & 2,107.6 \\\cline{1-3}\cline{5-8}
    1.7 & 27.196 & 53.546 & & 4.0 & 1,006.3 & 5,425.7 & 5,498.4 \\\cline{1-3}\cline{5-8}
    1.8 & 32.471 & 66.762 & & 5.0 & 4,183.4 & 34,911 & 35,924 \\\cline{1-3}\cline{5-8}
    1.9 & 38.634 & 82.963 & & 6.0 & 16,688 & $2.1571\cdot 10^5$ & $2.2561\cdot 10^5$
    \\\cline{1-3}\cline{5-8}
    2.0 & 45.831 & 102.80 & & 8.0 & $2.4510\cdot 10^5$ & $7.6097\cdot 10^6$ & $8.2094\cdot 10^6$
    \\\cline{1-3}\cline{5-8}
    2.1 & 54.228 & 127.06 & & 10.0 & $3.3583\cdot 10^6$ & $2.5069\cdot 10^8$ & $2.7692\cdot 10^8$
    \\\cline{1-3}\cline{5-8}
    2.2 & 64.016 & 156.72 & \multicolumn{5}{c}{} \\\cline{1-3}
    2.3 & 75.415 & 192.92 & \multicolumn{5}{c}{} \\\cline{1-3}
    2.4 & 88.681 & 237.08 & \multicolumn{5}{c}{} \\\cline{1-3}
    2.5 & 104.10 & 290.87 & \multicolumn{5}{c}{} \\\cline{1-3}
  \end{tabular}\vspace{2ex}
  \caption{Eigenvalues for the triangle with height 1.}
  \label{tab:triangle1_lambda}
\end{table}

Various intermediate path points $e_\mathrm{M}$ were used to compute
$u_2$. However, the function that CMPA converged to did not depend on
this choice. The symmetry properties of the computed $u_2$ depend only
on the value $p$. For $p\leq 2.6$ it is \emph{even} in $x_2$, i.e., it
belongs to
\begin{equation}
  \label{eq:traingleSE}
  \mathcal{S}_\mathrm{E}:=\{u:\Omega\to\R\,|\,u(x_1,x_2)=u(x_1,-x_2) \}.
\end{equation}
For $p\geq 2.7$ this symmetry is lost by $u_2$ as the graphs in
Fig.~\ref{fig:triangle1_u2} show.

For $p=2.6$ the computation was repeated with intermediate path points
$e_\mathrm{M}$ 
without symmetry
$\mathcal{S}_\mathrm{E}$ but CMPA always converged to the function
shown in Fig.~\ref{fig:triangle1_u2} which displays symmetry
$\mathcal{S}_\mathrm{E}$.

\begin{figure}[p]
  \centering
  \setlength{\unitlength}{1mm}
  \begin{picture}(150,57)
    \put(28,53){(a)}
    \put(2,2){\includegraphics[width=60mm]{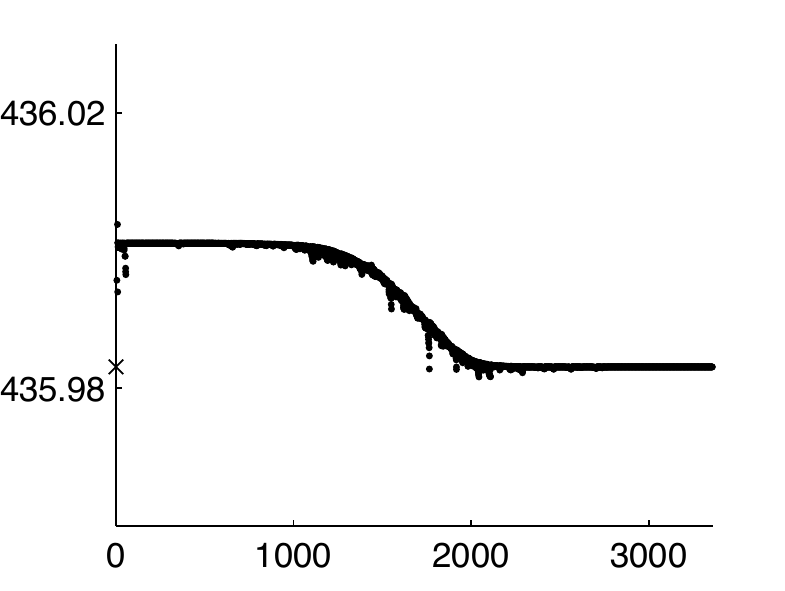}}
    \put(43,0){\small iterations}
    \put(7,46){$I(z^\mathrm{max})$}
    \put(37,30){$p=2.7$}
    \put(16,30.5){$\approx u_{2,\mathcal{S}_\mathrm{E}}$}
    \put(49,21){$u_2$}
    \put(7,19.5){\small $\lambda_2$}
    \put(60,0){\troneutwoextragraphs{3.5}{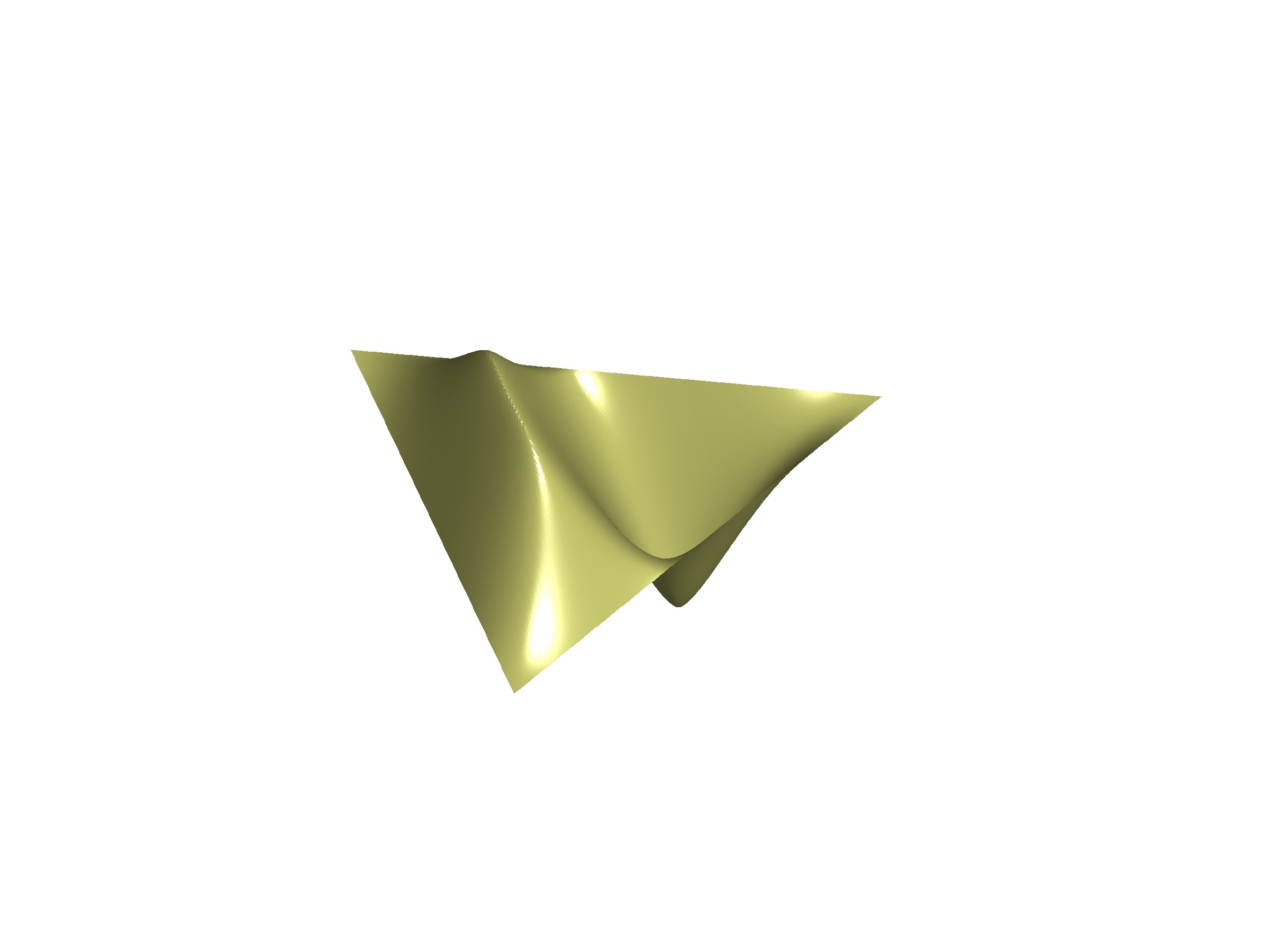}}
    \put(64,53){(b)}
    \put(90,0){\troneutwoextragraphs{6.0}{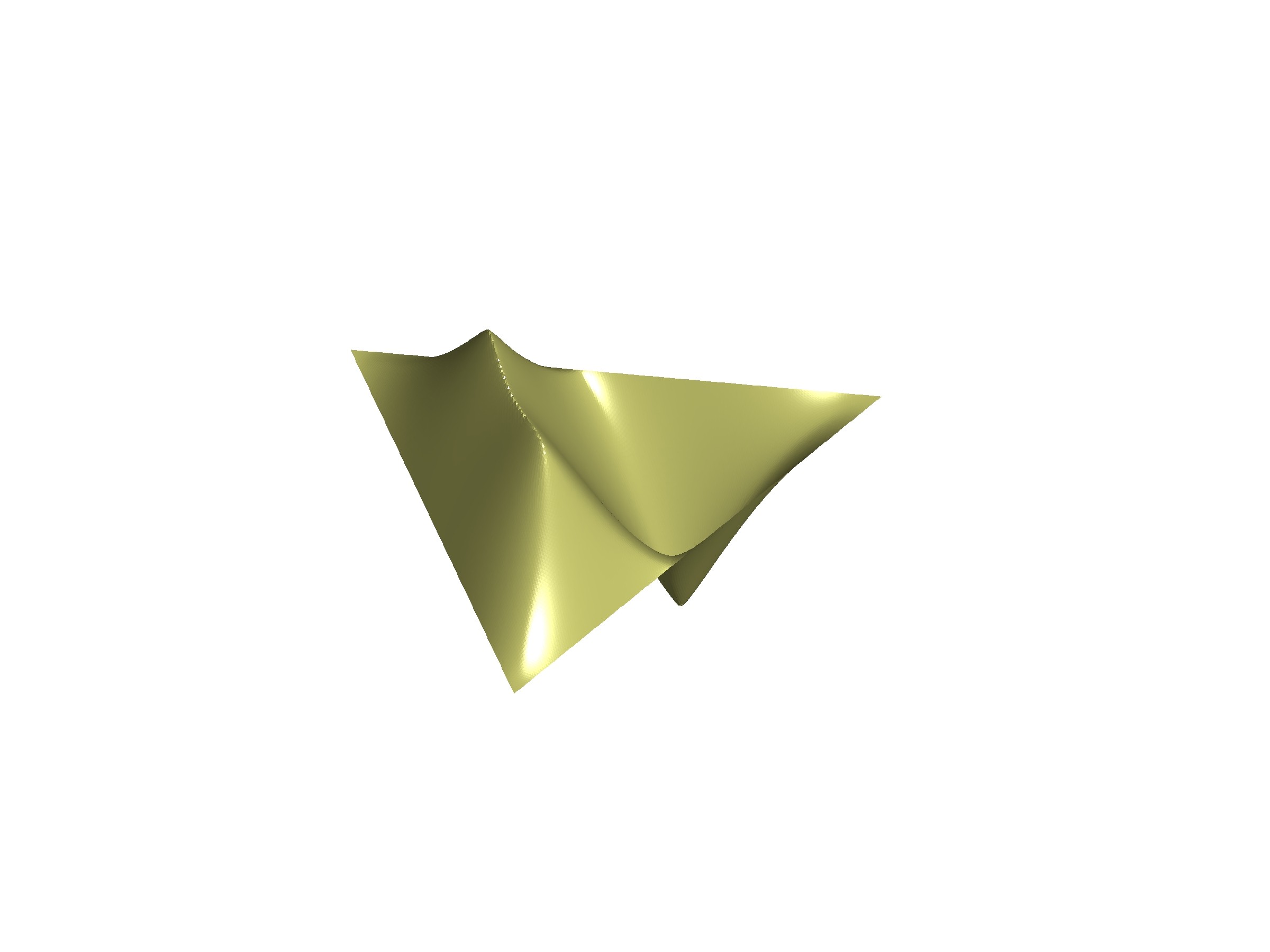}}
    \put(120,0){\troneutwoextragraphs{10.0}{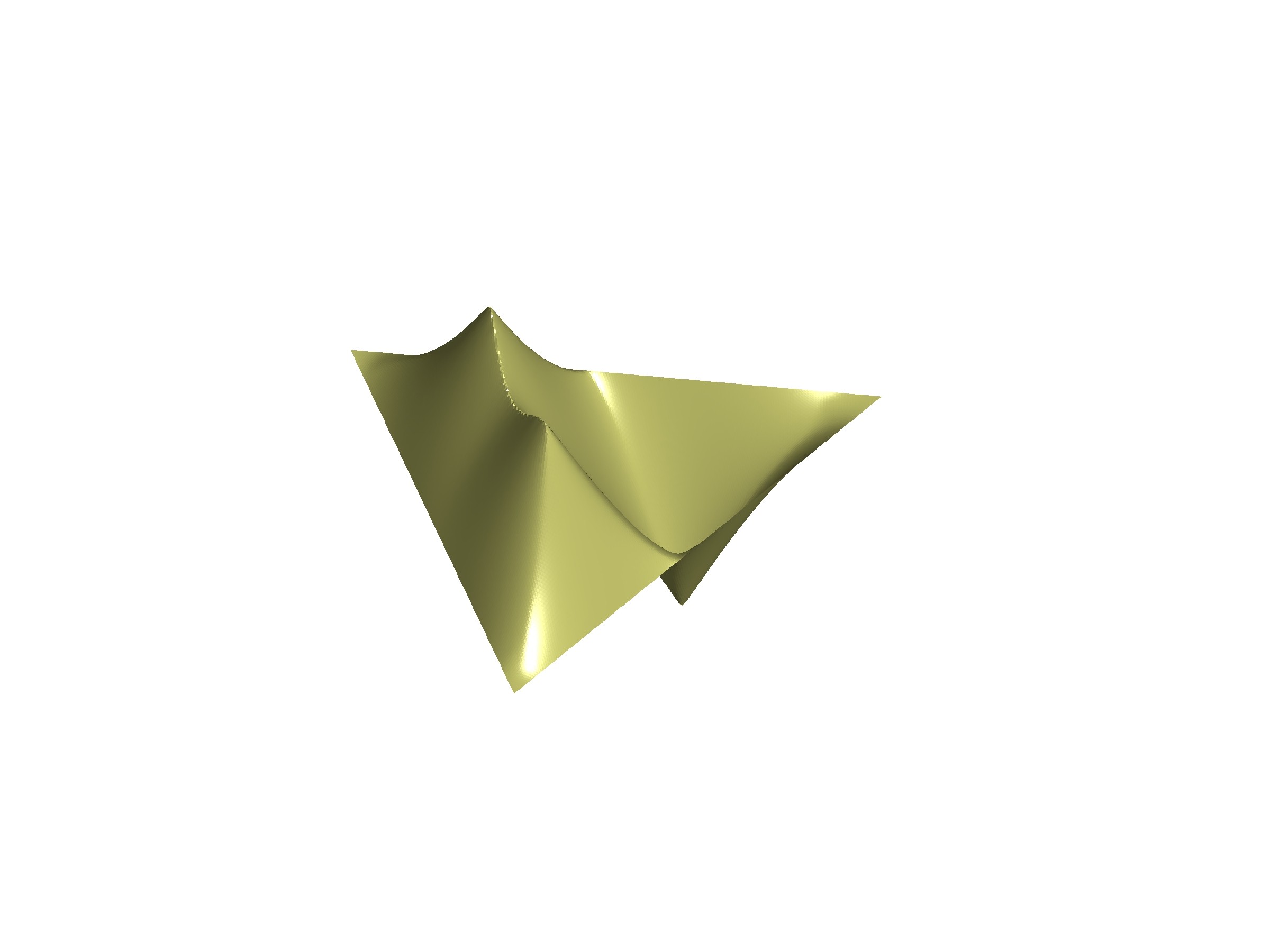}}
  \end{picture}
  \caption{Triangle with height 1: (a) Maximum value of the Dirichlet
    functional $I$ along the path during the run of CMPA for $p=2.7$
    and $e_\mathrm{M}\in\mathcal{S}_\mathrm{E}$. (b) The computed
    eigenfunction $u_{2,\mathcal{S}_\mathrm{E}}$ for $p=3.5, 6.0$, and
    10.0.}
  \label{fig:triangle1_u2extra}
\end{figure}
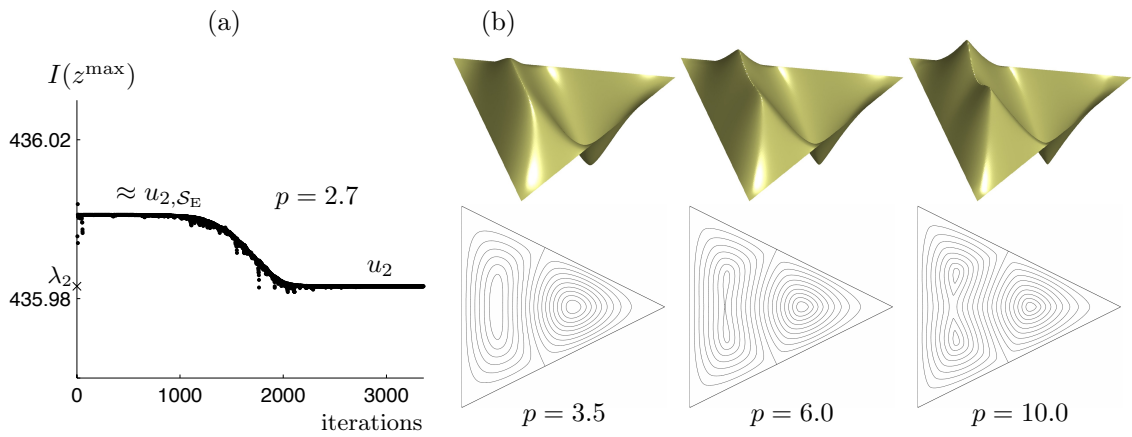

For $p=2.7$ a symmetric $e_\mathrm{M}\in\mathcal{S}_\mathrm{E}$ was
chosen. The graph in Fig.~\ref{fig:triangle1_u2extra}(a) shows how
the maximum of the Dirichlet functional $I$ along the path evolved
during this run of CMPA. The path connecting $u_1$ with $-u_1$ which
gets deformed at every step of CMPA seems to stay close to some
critical point having symmetry $\mathcal{S}_\mathrm{E}$ during the
first 1000 steps but then it slips down to lower values of $I$ and
stays close to another critical point. This is the asymmetric $u_2$
which the algorithm eventually converges to.

Even beyond $p=2.6$ there exist eigenfunctions with symmetry
$\mathcal{S}_\mathrm{E}$. Let $u_{2,\mathcal{S}_\mathrm{E}}$ denote a
sign-changing eigenfunction of the $p$-Laplace operator on $\Omega$
which lies in $\mathcal{S}_\mathrm{E}$ and has the smallest eigenvalue
(which we denote $\lambda_{2,\mathcal{S}_\mathrm{E}}$). As mentioned
above, for $p\leq 2.6$ we observed that
$u_2=u_{2,\mathcal{S}_\mathrm{E}}$ (up to scaling). To compute
$u_{2,\mathcal{S}_\mathrm{E}}$ for $p\geq 2.7$ consider the following
eigenvalue problem:
\begin{equation}
  \label{eq:evproblem_triangle1_modified}
  \begin{aligned}
    -\plap u &= \lambda |u|^{p-2} u &&\text{in }\Omega^\mathrm{half}, \\
    \textstyle\frac{\partial u}{\partial n} &= 0 &&\text{on }\Gamma_1, \\
    u &= 0 &&\text{on }\Gamma_2,
  \end{aligned}
\end{equation}
where
\begin{equation}
  \label{eq:halftriangle_def}
  \begin{aligned}
    \Omega^\mathrm{half}&=\left\{(x_1,x_2)\in\Omega\,|\,x_2>0\right\}, \\
    \Gamma_1&=\big\{(x_1,x_2)\in\partial\Omega^\mathrm{half}\,\big|\,x_2=0\big\}, \\
    \Gamma_2&=\partial\Omega^\mathrm{half}\setminus\Gamma_1.
  \end{aligned}
\end{equation}
Any eigenfunction solving this problem can be extended to an
eigenfunction on the whole $\Omega$ by even symmetry about
$x_2=0$. Since the first eigenfunction of the original problem
\pref{eq:evproblem} for the triangle $\Omega$ belongs to
$\mathcal{S}_\mathrm{E}$, its restriction to $\Omega^\mathrm{half}$ is
the first eigenfunction for
\pref{eq:evproblem_triangle1_modified}. Hence to compute
$u_{2,\mathcal{S}_\mathrm{E}}$ we just need to apply CMPA to problem
\pref{eq:evproblem_triangle1_modified} with paths which again connect
$u_1$ and $-u_1$. The modification of the finite element method to
take into account the natural boundary condition on $\Gamma_1$ is
straightforward.
The computed values of $\lambda_{2,\mathcal{S}_\mathrm{E}}$ are listed
in Table~\ref{tab:triangle1_lambda}. Figure
\ref{fig:triangle1_u2extra}(b) shows the corresponding eigenfunction
$u_{2,\mathcal{S}_\mathrm{E}}$ for selected values of $p$.

\begin{figure}[t]
  \centering
  \setlength{\unitlength}{1mm}
  \begin{picture}(150,50)(0,2)
    \put(0,0){\includegraphics[width=7.3cm]{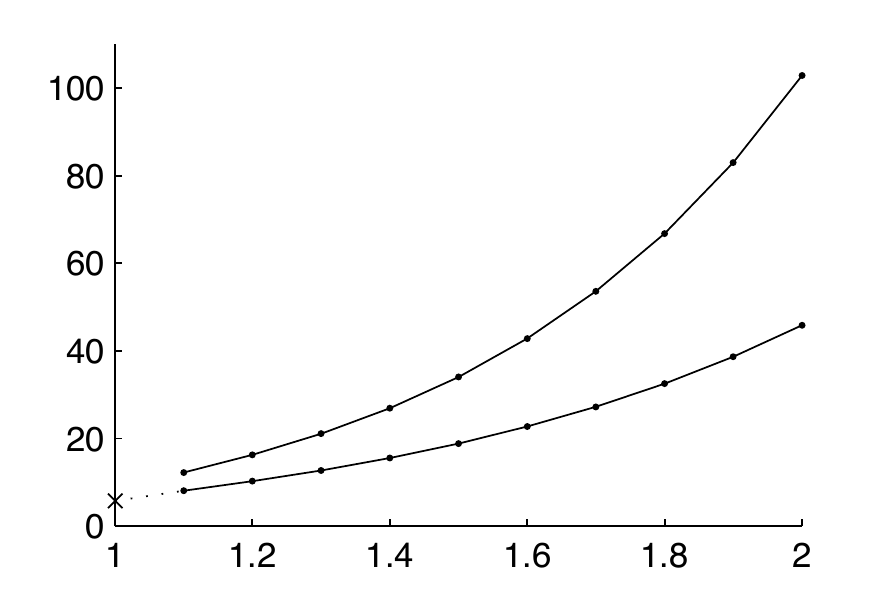}}
    \put(6.5,45.5){$\lambda$}
    \put(69,4.5){$p$}
    \put(56,21){$\lambda_1$}
    \put(56,37){$\lambda_2$}
    \put(0.2,7.5){\small $h_1(\Omega)$}
    \put(74,0){\includegraphics[width=7.3cm]{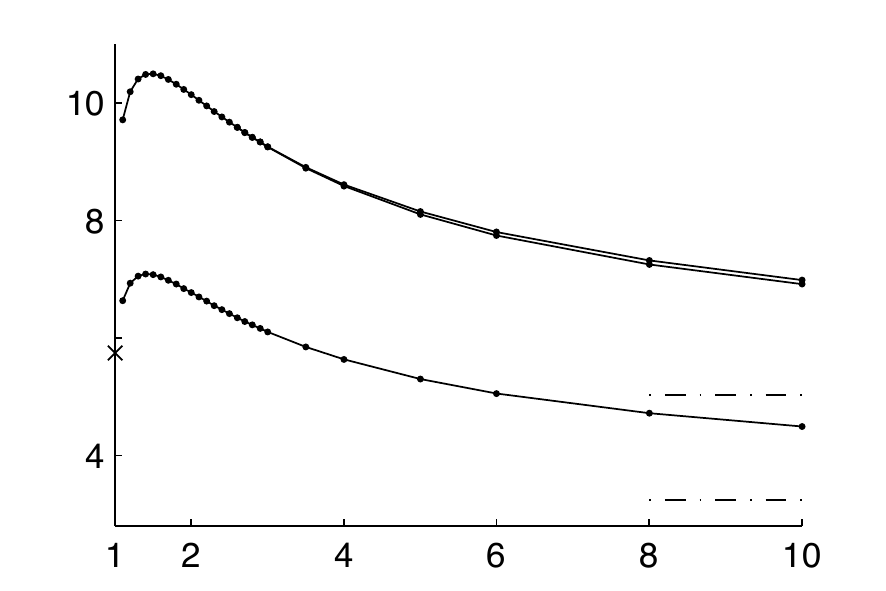}}
    \put(80.5,45.5){$\lambda^{1/p}$}
    \put(143,4){$p$}
    \put(111,11){$(\lambda_1)^{1/p}$}
    \put(111,24){$(\lambda_2)^{1/p}$}
    \put(121,24){\line(1,0){4}}
    \put(125,24){\vector(0,1){2.7}}
    \put(111,32.5){$(\lambda_{2,\mathcal{S}_\mathrm{E}})^{1/p}$}
    \put(124,32){\line(1,0){1}}
    \put(125,32){\vector(0,-1){3.8}}
    \put(74,19){\small $h_1(\Omega)$}
    \put(141.5,7){\small $\Lambda_1(\Omega)$}
    \put(141.5,15.5){\small $\Lambda_2(\Omega)$}
  \end{picture}  
  \caption{Dependence of the numerically computed eigenvalues for the
    triangle with height 1 on $p$.}
  \label{fig:triangle1_lambda_plot}
\end{figure}
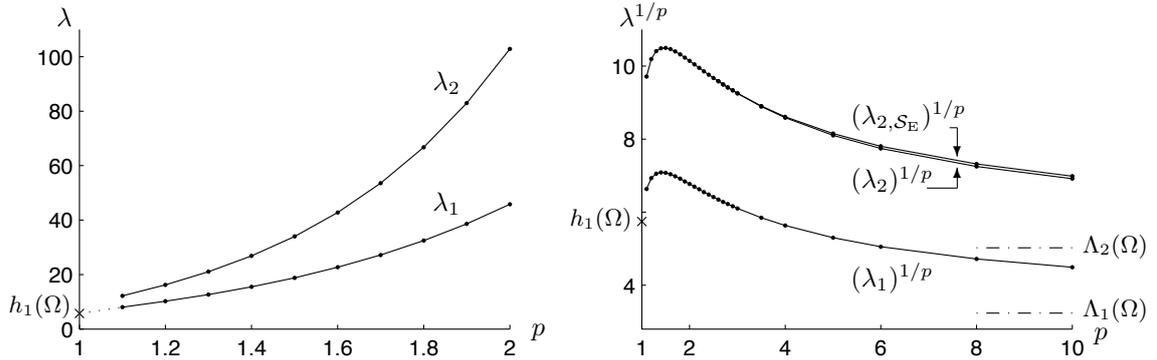

The dependence of the eigenvalues $\lambda_1$, $\lambda_2$, and
$\lambda_{2,\mathcal{S}_\mathrm{E}}$ on $p$ is further plotted in
Fig.~\ref{fig:triangle1_lambda_plot}. The figure also shows the limits
of $\lambda_1$ for $p\to1$ and $\infty$ and of $\lambda_2$ for
$p\to\infty$ which can be computed explicitly. As mentioned, for
example, in \cite{LachandRobertKawohl}, the Cheeger constant of a
triangle is given by
$h_1(\Omega)=(\per(\Omega)+\sqrt{4\pi|\Omega|})/(2|\Omega|)$ and
hence in our case $h_1(\Omega)=1+\sqrt{5}+\sqrt{2\pi}\approx
5.7427$. Simple computations yield
$\Lambda_1(\Omega)=1+\sqrt{5}\approx 3.2361$ and
$\Lambda_2(\Omega)=1+9/\sqrt{5}\approx 5.0249$.

\subsection{Triangle with height 3/4}\label{triangle3}
\begin{figure}[p]
  \centering
  \setlength{\unitlength}{1mm}
  \begin{picture}(150,56)
    \put(0,0){\trthreeuonegraphs{1.1}{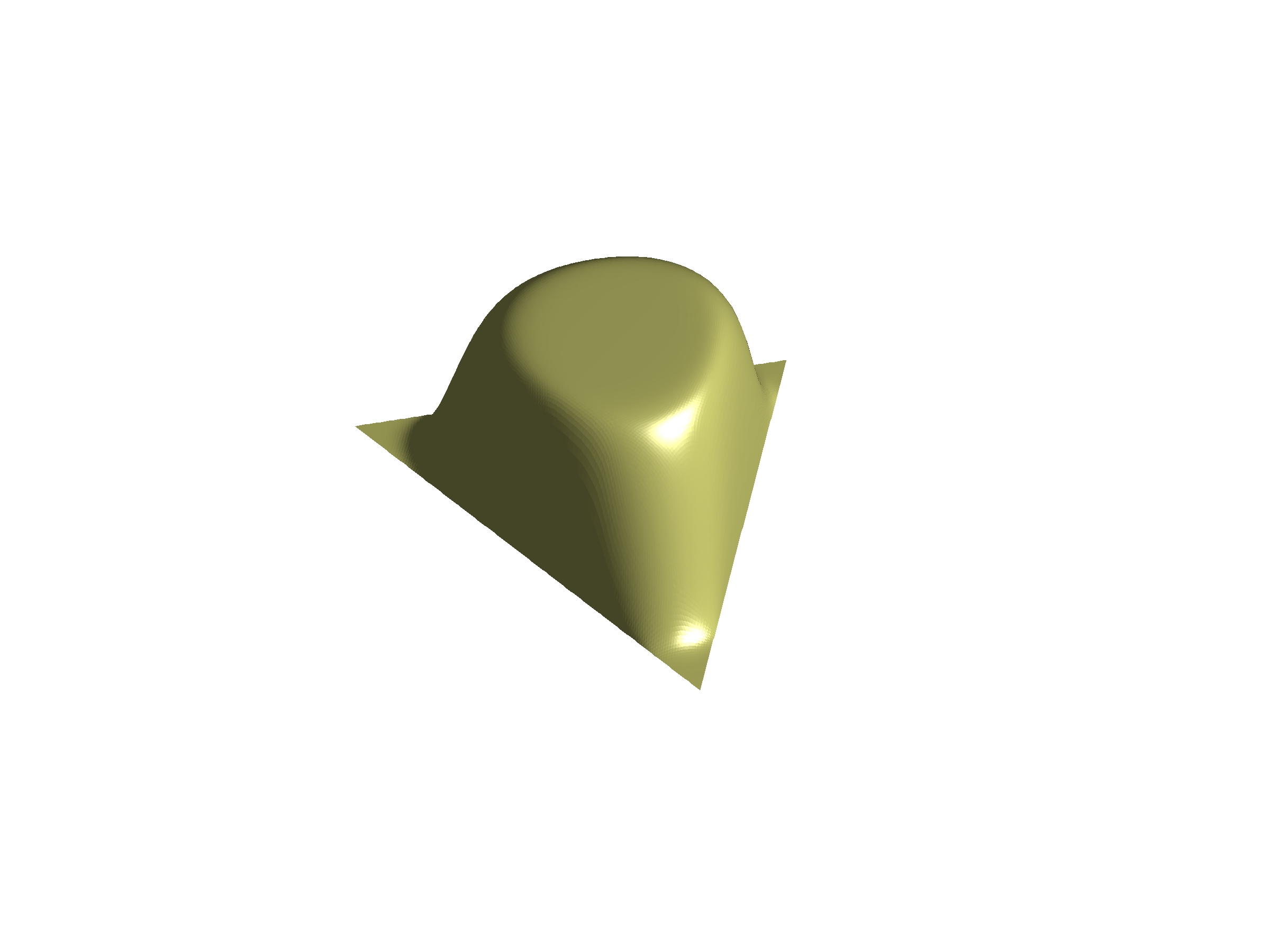}}
    \put(30,0){\trthreeuonegraphs{1.4}{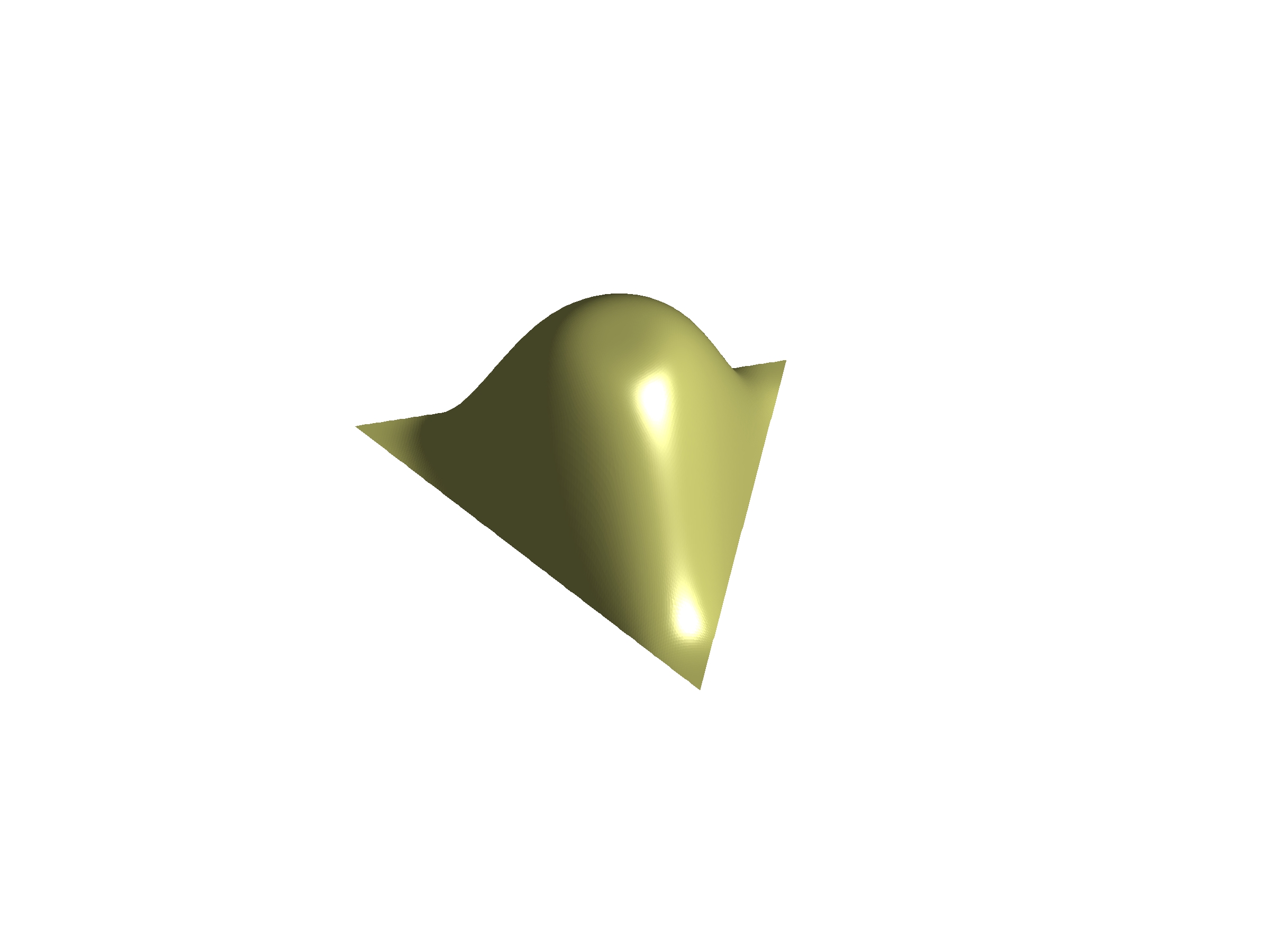}}
    \put(60,0){\trthreeuonegraphs{2.0}{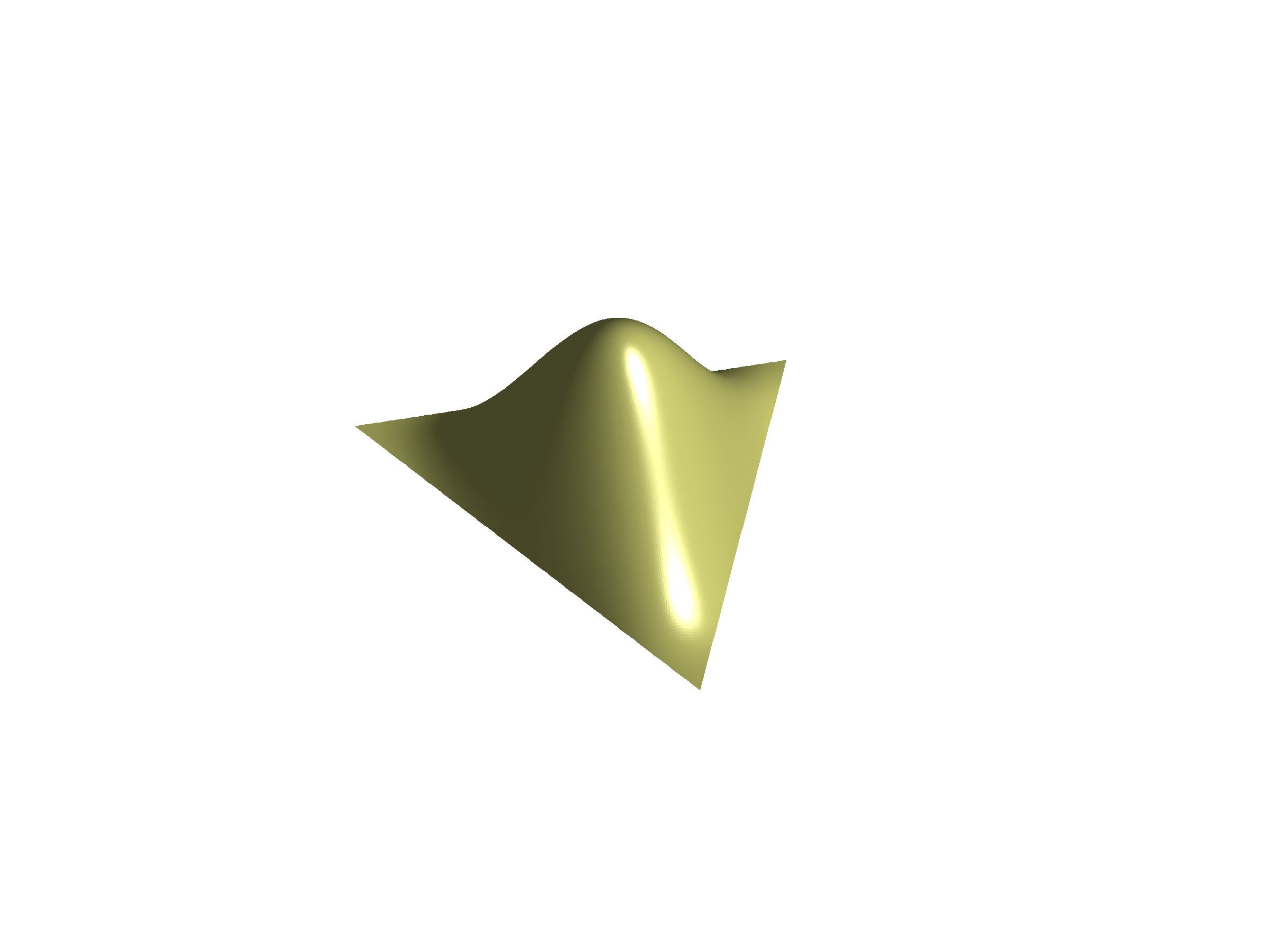}}
    \put(90,0){\trthreeuonegraphs{2.5}{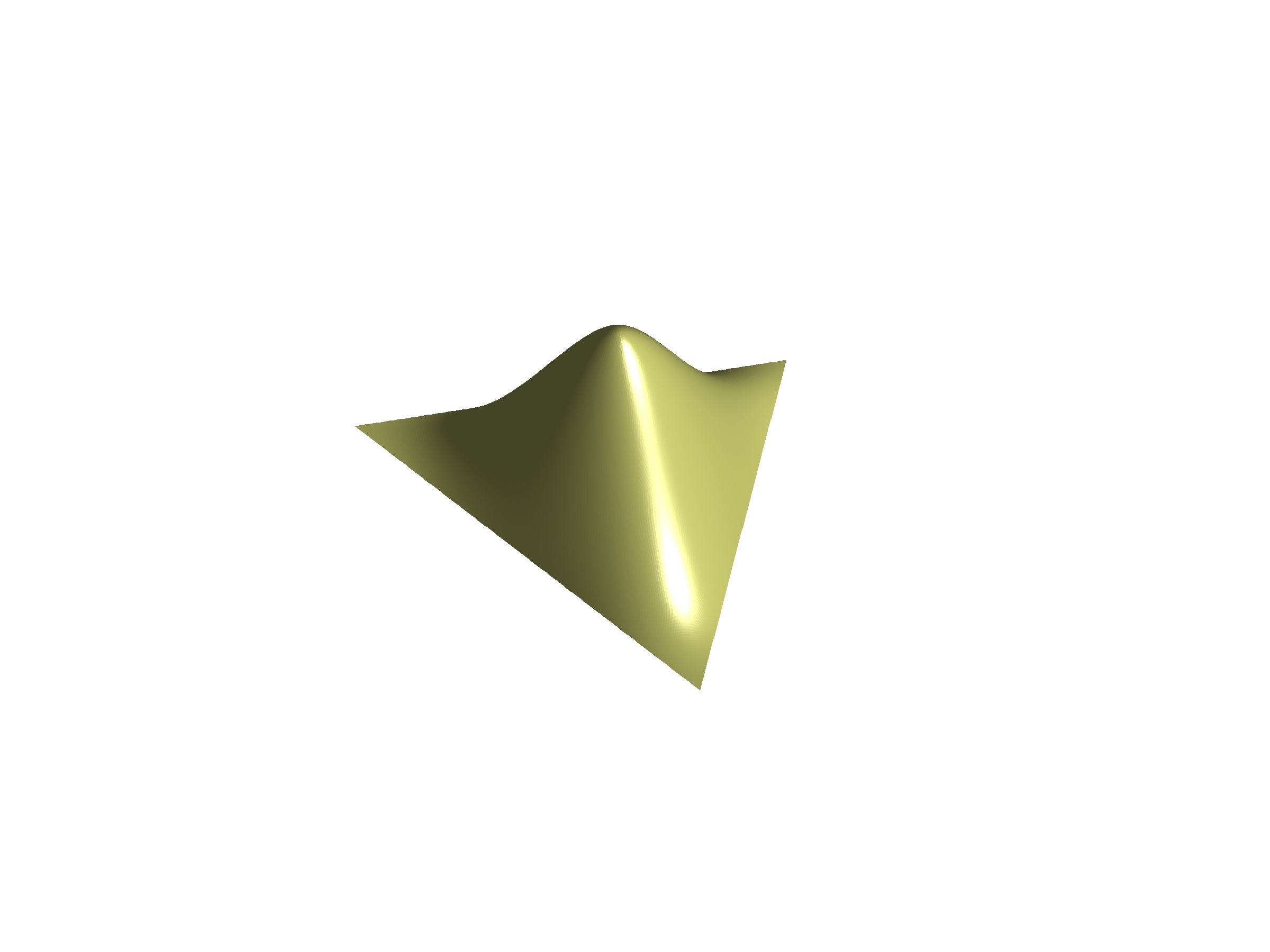}}
    \put(120,0){\trthreeuonegraphs{8.0}{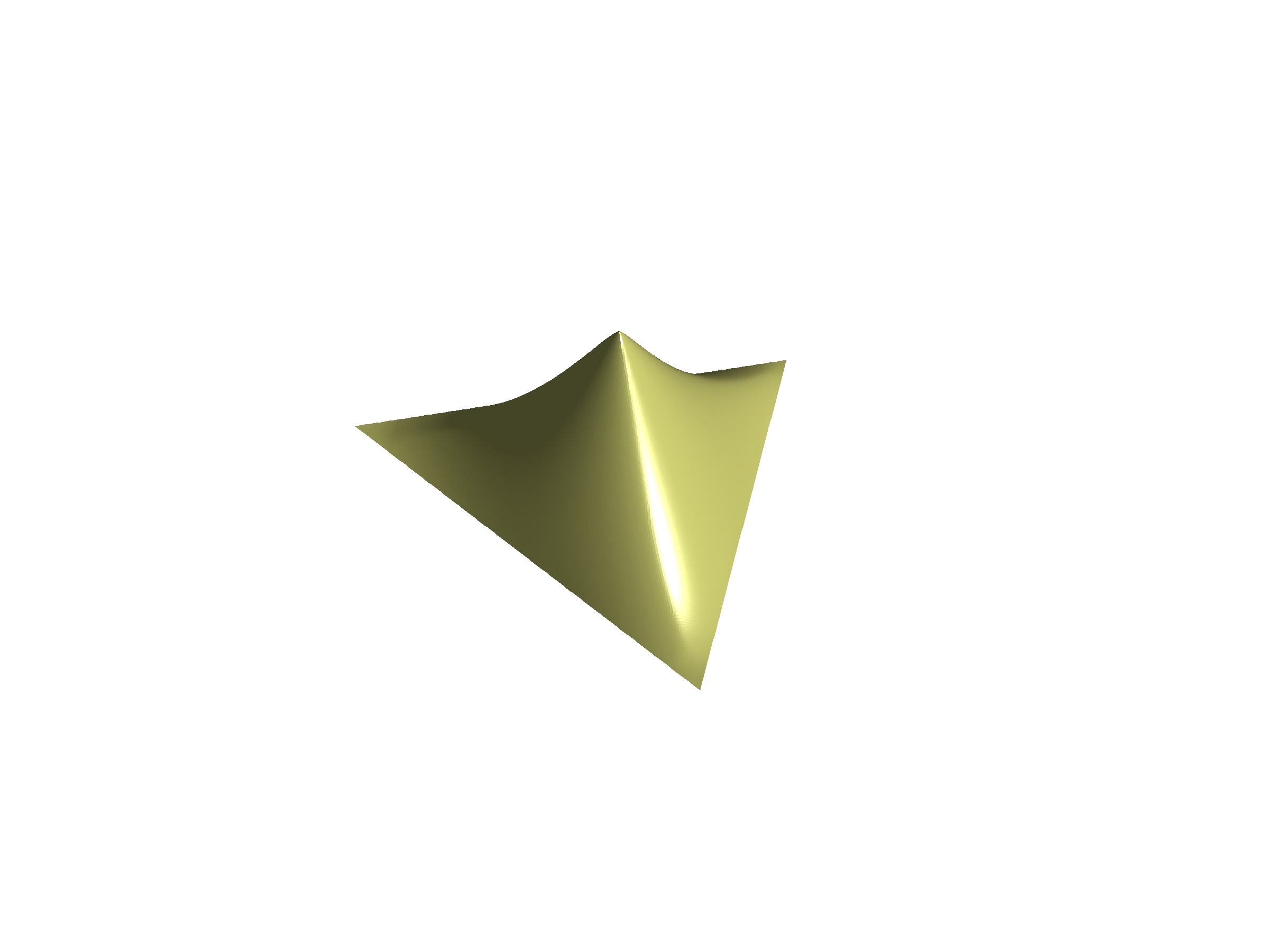}}
    \put(.5,14){\rotatebox{90}{\scriptsize base}}
    \put(9.5,35.5){\rotatebox{-37}{\scriptsize base}}
  \end{picture}
  \caption{The numerically computed first eigenfunction $u_1$ for the triangle with height 3/4.}
  \label{fig:triangle3_u1}
\end{figure}

\begin{figure}[p]
  \centering
  \setlength{\unitlength}{1mm}
  \begin{picture}(150,56)
    \put(0,0){\trthreeutwographs{1.1}{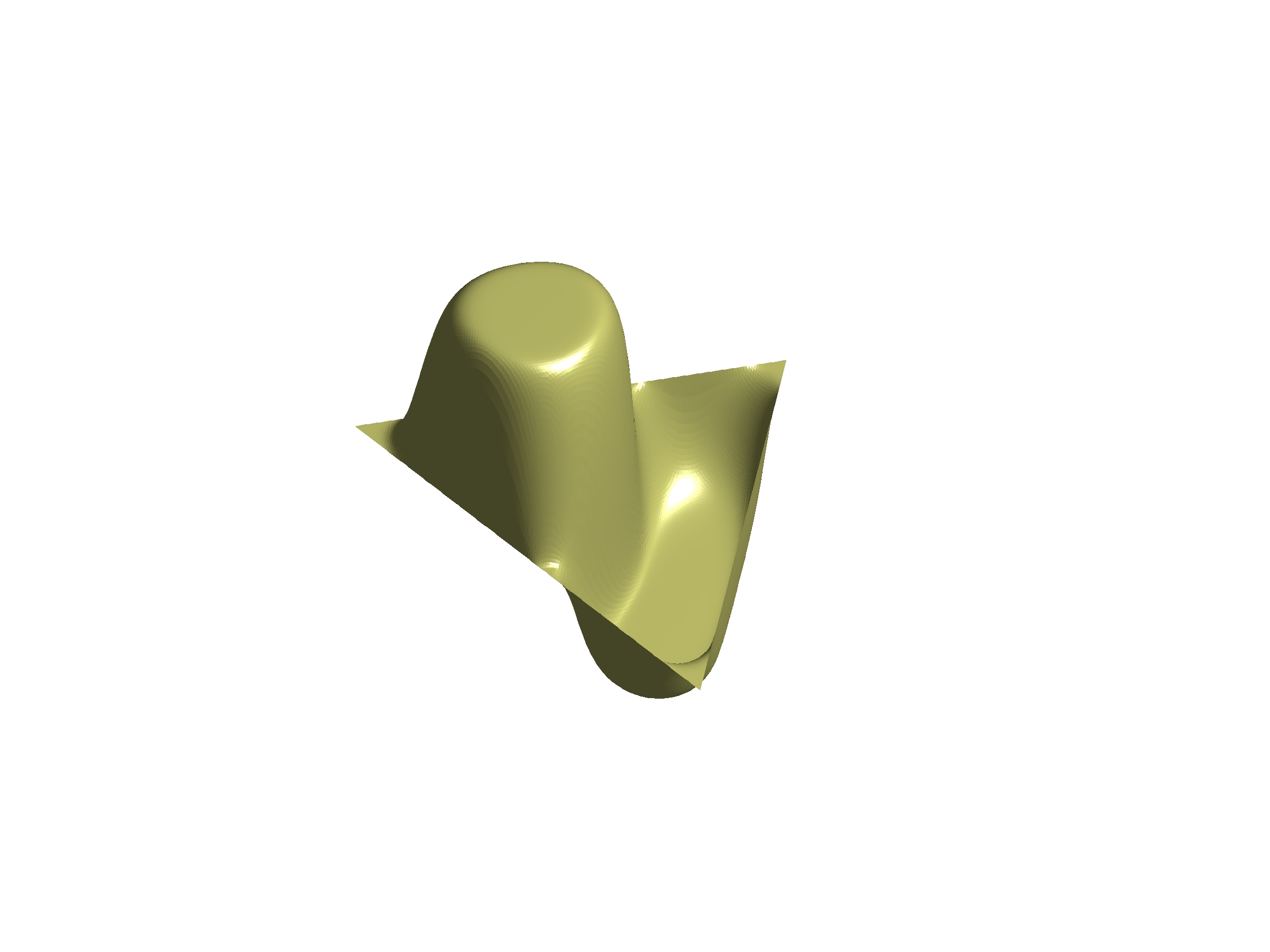}}
    \put(30,0){\trthreeutwographs{1.4}{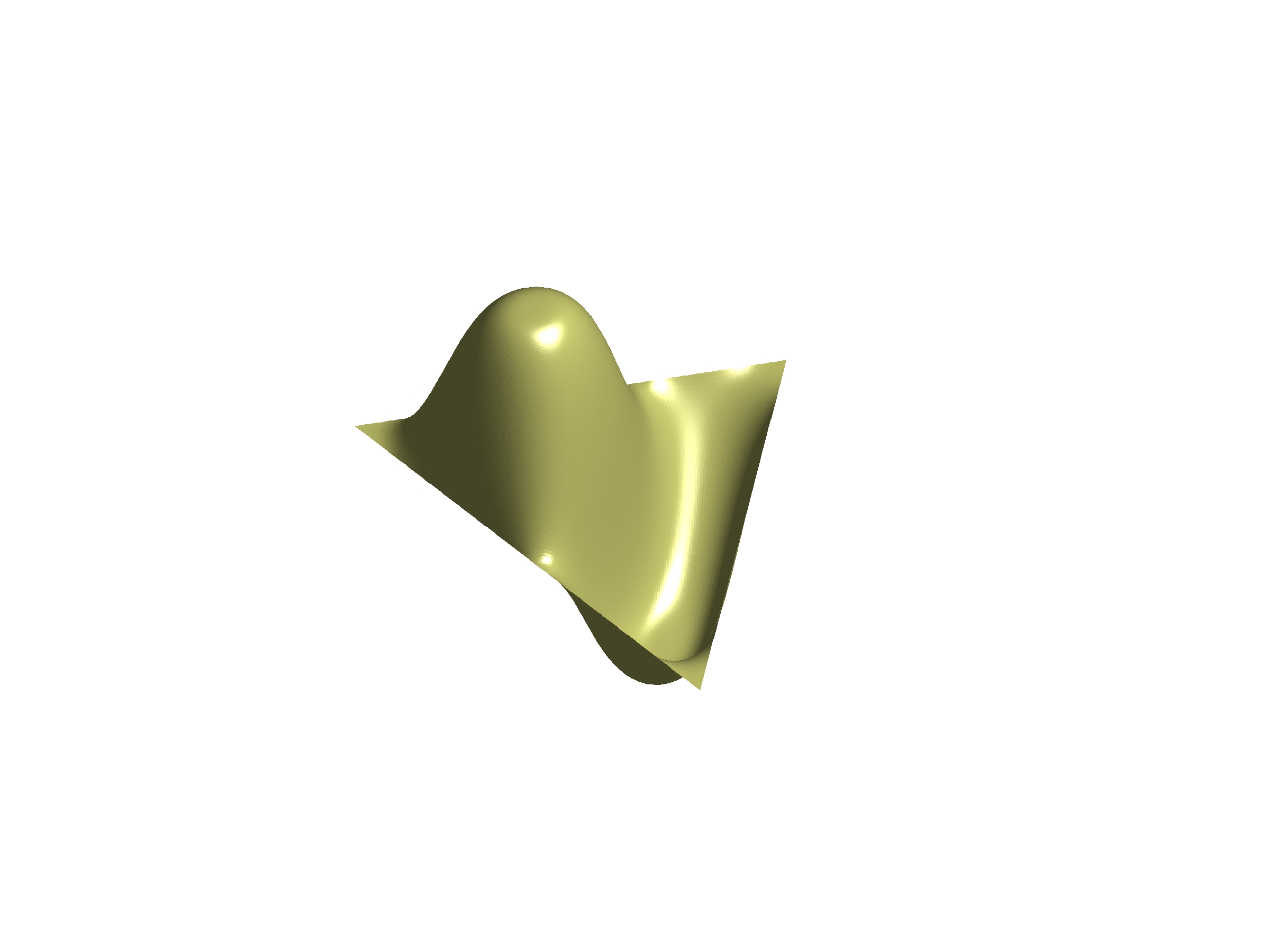}}
    \put(60,0){\trthreeutwographs{1.6}{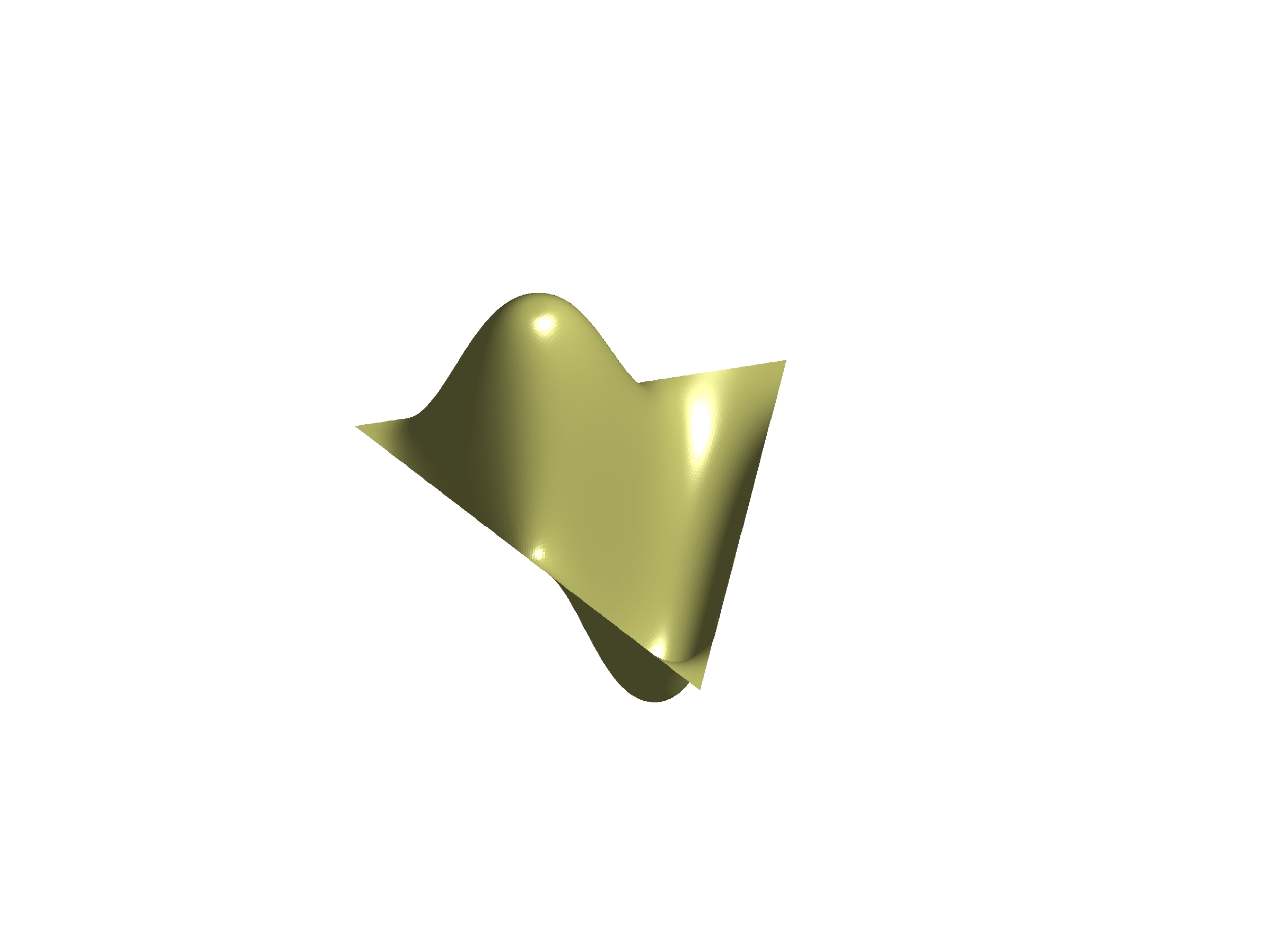}}
    \put(90,0){\trthreeutwographs{1.7}{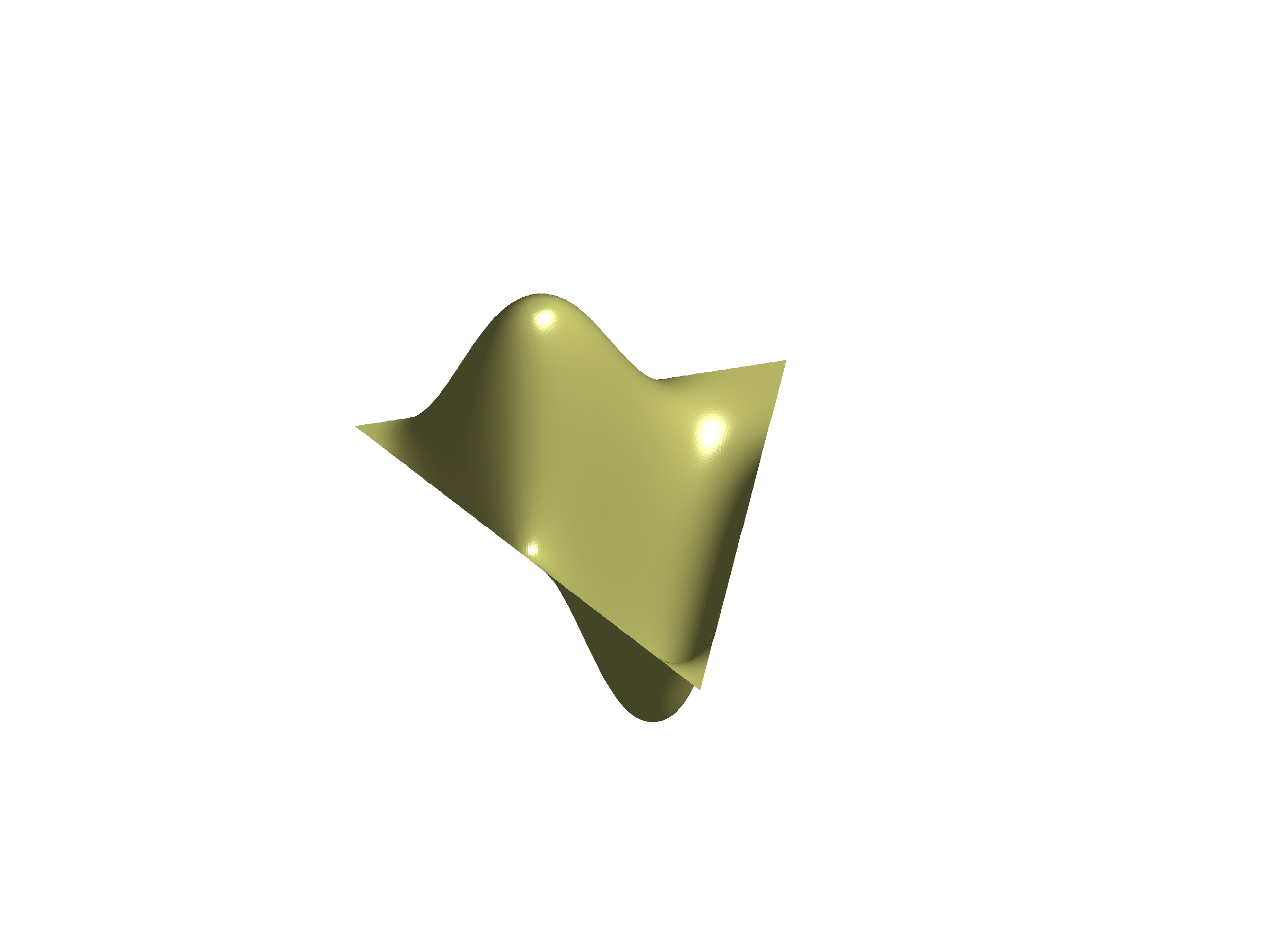}}
    \put(120,0){\trthreeutwographs{8.0}{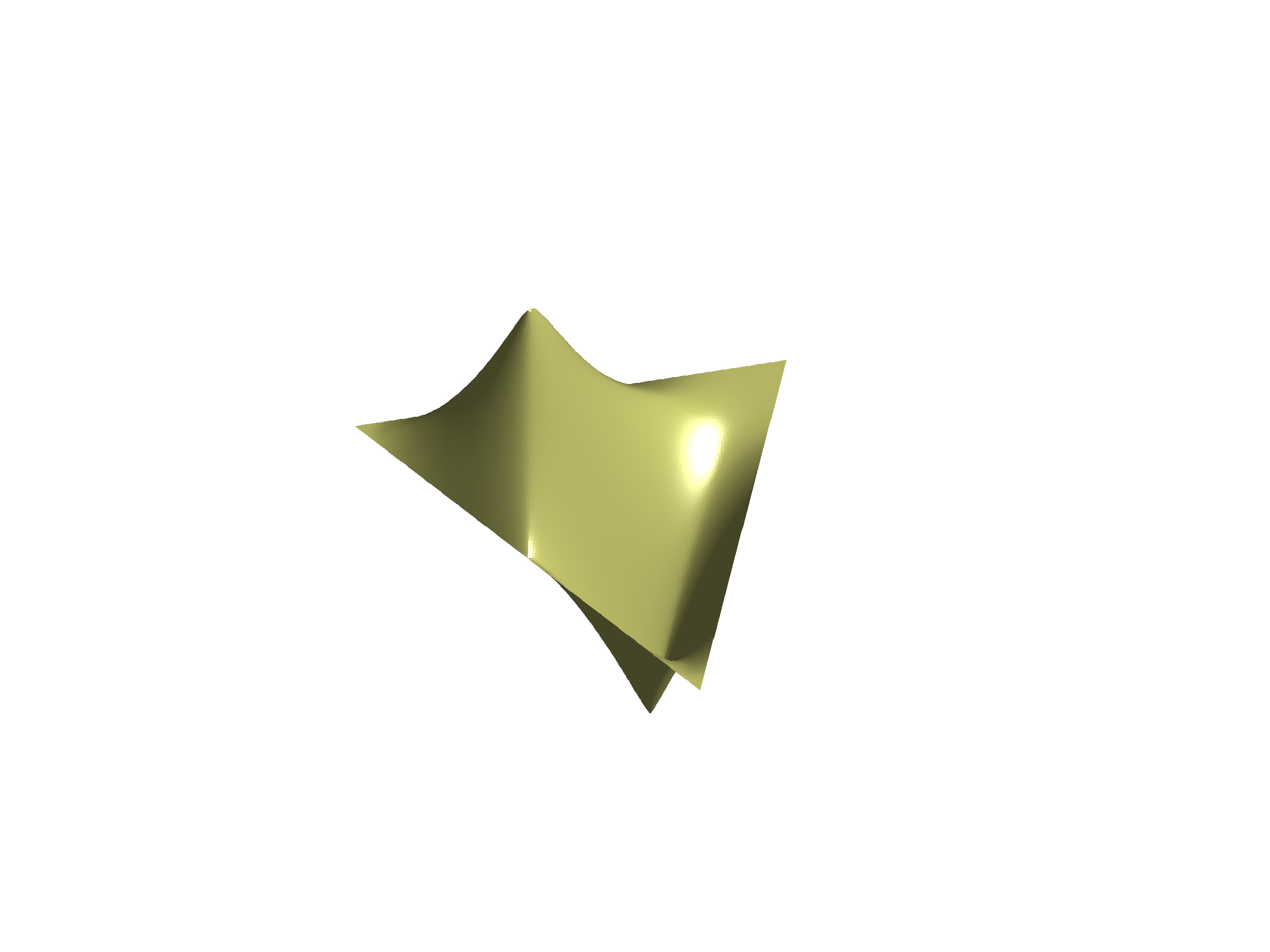}}
  \end{picture}
  \caption{The numerically computed second eigenfunction $u_2$ for the triangle with height 3/4.}
  \label{fig:triangle3_u2}
\end{figure}

Let
\begin{displaymath}\textstyle
  \Omega=\left\{(x_1,x_2)\in\R^2\,\big|\, x_1\in\left(0,\frac{3}{4}\right),
    |x_2|<\frac{2}{3}\left(\frac{3}{4}-x_1\right)\right\}
\end{displaymath}
be an isosceles triangle with base 1 and height 3/4. It was
discretized using 28,672 triangles. Figures \ref{fig:triangle3_u1} and
\ref{fig:triangle3_u2} show the eigenfunctions $u_1$ and $u_2$ for
several values of $p$, respectively. Table \ref{tab:triangle3_lambda}
lists the corresponding values of $\lambda_1$ and $\lambda_2$.

\begin{table}[b]
  \centering\vspace{1ex}
  \begin{tabular}[t]{|c|c|c|c|c|c|c|c|c|c|c|c|}
    \cline{1-4}\cline{6-8}\cline{10-12}
    $p$ & $\lambda_1$ & $\lambda_2$ & $\lambda_{\mathcal{S}_\mathrm{O}}$ & \qquad\qquad &
    $p$ & $\lambda_1$ & $\lambda_2$ & \qquad\qquad & $p$ & $\lambda_1$ & $\lambda_2$
    \\\cline{1-4}\cline{6-8}\cline{10-12}
    \multicolumn{12}{c}{\vspace{\tableheadsep}}\\\cline{1-4}\cline{6-8}\cline{10-12}
    1.1 & 9.389 & 14.38 & 14.50 && 1.8 & 42.07 & 86.84 && 2.5 & 149.1 & 413.0 \\\cline{1-4}\cline{6-8}\cline{10-12}
    1.2 & 12.13 & 19.41 & 19.52 && 1.9 & 50.78 & 109.2 && 3.0 & 350.0 & 1,196 \\\cline{1-4}\cline{6-8}\cline{10-12}
    1.3 & 15.28 & 25.53 & 25.62 && 2.0 & 61.11 & 137.1 && 4.0 & 1,789 & 9,351 \\\cline{1-4}\cline{6-8}\cline{10-12}
    1.4 & 18.97 & 33.10 & 33.17 && 2.1 & 73.36 & 171.6 && 5.0 & 8,591 & $6.871\cdot 10^4$ \\\cline{1-4}\cline{6-8}\cline{10-12}
    1.5 & 23.35 & 42.52 & 42.55 && 2.2 & 87.85 & 214.4 && 6.0 & $3.958\cdot 10^4$ & $4.849\cdot 10^5$ \\\cline{1-4}\cline{6-8}\cline{10-12}
    1.6 & 28.55 & 54.22 & 54.22 && 2.3 & 105.0 & 267.2 && 8.0 & $7.752\cdot 10^5$ & $2.232\cdot 10^7$ \\\cline{1-4}\cline{6-8}\cline{10-12}
    1.7 & 34.73 & 68.76 & 68.76 && 2.4 & 125.2 & 332.5 && 10.0 & $1.416\cdot 10^7$ & $9.602\cdot 10^8$ \\\cline{1-4}\cline{6-8}\cline{10-12}
  \end{tabular}\vspace{2ex}
  \caption{Eigenvalues for the triangle with height 3/4.}
  \label{tab:triangle3_lambda}
\end{table}

\begin{figure}[p]
  \centering
  \setlength{\unitlength}{1mm}
  \begin{picture}(150,59)
    \put(15,0){\trthreeutwoextragraphs{$u_2$}{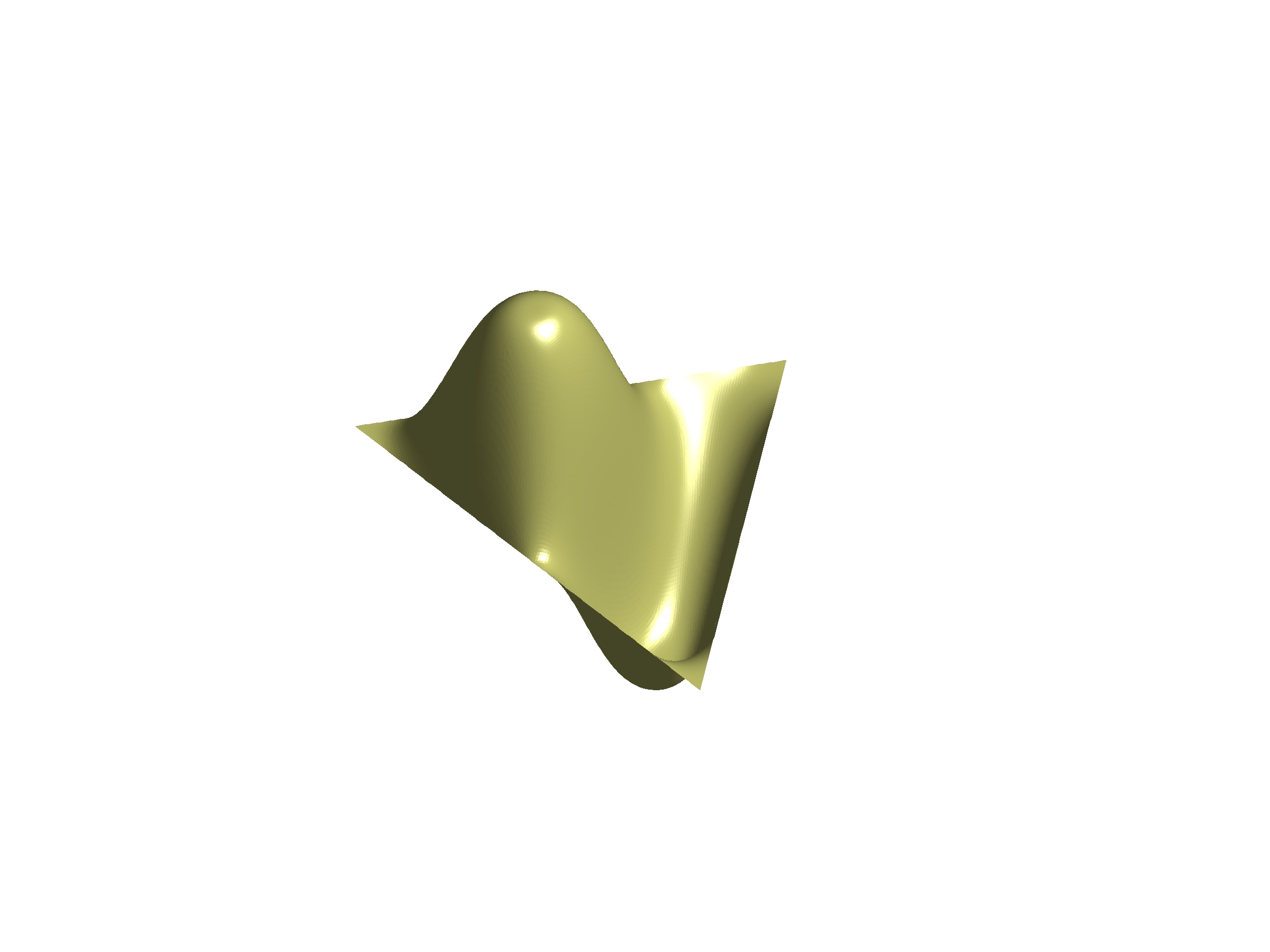}}
    \put(45,0){\trthreeutwoextragraphs{$u_{\mathcal{S}_\mathrm{O}}$}{out_u1_p_1-5}}
    \put(75,0){\trthreeutwoextragraphs{$u_{(2,\mathcal{S}_\mathrm{E})}$}{out_u2_p_1-5b}}
    \put(107,30){\parbox{30mm}{
        \begin{tabular}{r@{$\ =\ $}l}
          $p$ & 1.5\\
          \rule{0pt}{10mm}$\lambda_2$ & 42.52\\
          $\lambda_{\mathcal{S}_\mathrm{O}}$ & 42.55\\
          $\lambda_{(2,\mathcal{S}_\mathrm{E})}$ & 44.42
        \end{tabular}}}
  \end{picture}
  \caption{Higher eigenfunctions for the triangle with height 3/4 for
    $p=1.5$: $u_2$ and $u_{(2,\mathcal{S}_\mathrm{E})}$ computed as
    constrained local mountain pass points by CMPA with no a priori
    assumptions on symmetry, $u_{\mathcal{S}_\mathrm{O}}$ computed by
    CDM enforcing symmetry $\mathcal{S}_\mathrm{O}$.}
  \label{fig:triangle3_u2extra}
\end{figure}
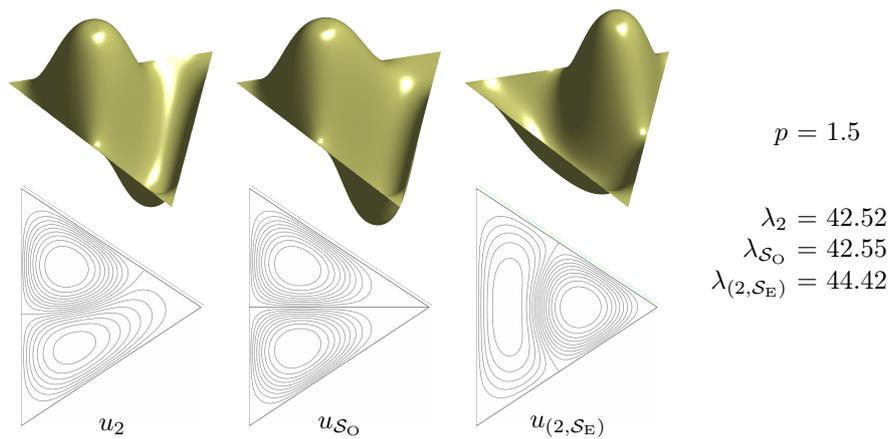

The symmetry properties of the computed $u_2$ change again with
$p$. For this triangle, however, $u_2$ gains more symmetry as $p$
increases (unlike for the triangle with height 1 where $u_2$ lost
symmetry). For $p\leq 1.6$ the nodal line of $u_2$ connects the base
of the triangle with one of its other sides. For $p\geq 1.7$ this
nodal line connects the base with the vertex above the base and $u_2$
is \emph{odd} in $x_2$, i.e., it belongs to
\begin{equation}
  \label{eq:triangleSO}
  \mathcal{S}_\mathrm{O}:=\{u:\Omega\to\R\,|\,u(x_1,x_2)=-u(x_1,-x_2) \}
\end{equation}
as can be seen in Fig.~\ref{fig:triangle3_u2}. It is because of lack
of resolution of the numerical method close to the vertex (where $u_2$
is flat) that the zero contour line in the figure for $p=1.7$ and
$p=8.0$ does not exactly reach the vertex.

For $p\leq 1.6$ there also exist eigenfunctions with symmetry
$\mathcal{S}_\mathrm{O}$. Let $\lambda_{\mathcal{S}_\mathrm{O}}$
denote the smallest eigenvalue with an eigenfunction belonging to
$\mathcal{S}_\mathrm{O}$ (denoted $u_{\mathcal{S}_\mathrm{O}}$). Using
the notation defined in \pref{eq:halftriangle_def} this eigenvalue can
be computed using CDM on $\Omega$ with an additional boundary
condition $u=0$ on $\Gamma_1$ or as the first eigenvalue on the
half-domain $\Omega^\mathrm{half}$. The computed values of
$\lambda_{\mathcal{S}_\mathrm{O}}$ are also listed in
Table~\ref{tab:triangle3_lambda}. For $p\geq 1.7$ the values of
$\lambda_2$ and $\lambda_{\mathcal{S}_\mathrm{O}}$ coincide. For
$p=1.6$ they differ in the sixth digit.

As in the previous computations, various choices of the intermediate
path point $e_\mathrm{M}$ were used to compute $u_2$ on $\Omega$ with
no a priori assumptions on symmetry. In some cases CMPA converged to
different functions depending on this choice (different local mountain
passes). For example, for $p=1.5$ two eigenfunctions were found: one
with a nodal line connecting the base of the triangle with one of its
sides, and another one with a nodal line connecting the two sides and
having an even symmetry in $x_2$. Both eigenfunctions are
(numerically) local mountain pass points of $I$ with respect to the
constraint $S$. The first one is called $u_2$ since it has the
smallest eigenvalue, the second one is called
$u_{(2,\mathcal{S}_\mathrm{E})}$ because of its symmetry (it could
also be understood as a solution of
\pref{eq:evproblem_triangle1_modified} formulated in a similar way for
the triangle with height 3/4). Both eigenfunctions are shown in
Fig.~\ref{fig:triangle3_u2extra} together with
$u_{\mathcal{S}_\mathrm{O}}$ for comparison.

\begin{figure}[t]
  \centering
  \setlength{\unitlength}{1mm}
  \begin{picture}(150,50)(0,2)
    \put(0,0){\includegraphics[width=7.3cm]{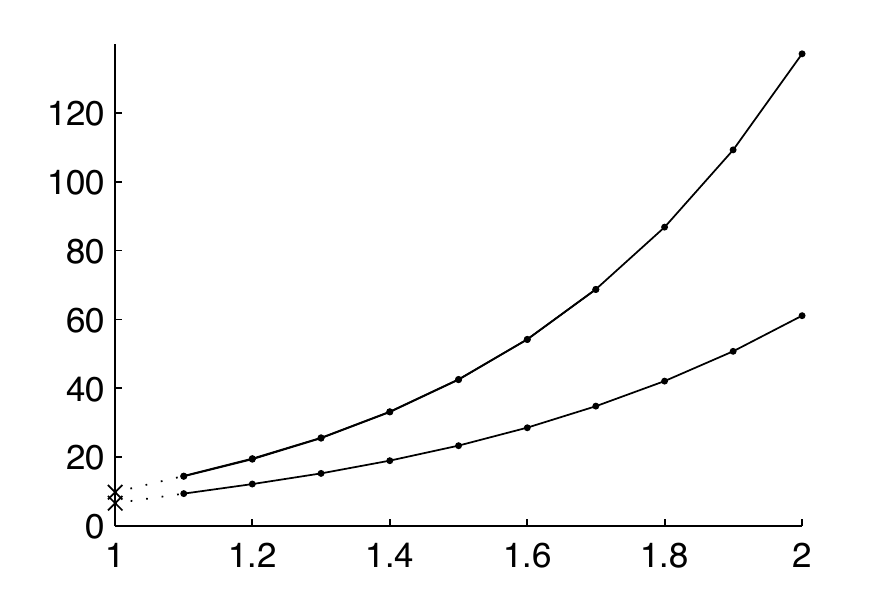}}
    \put(13.4,7.8){\includegraphics[width=1.8cm]{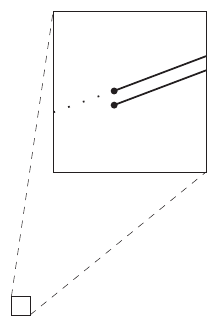}}
    \put(6.5,45.5){$\lambda$}
    \put(69,4.5){$p$}
    \put(56,21){$\lambda_1$}
    \put(56,38){$\lambda_2$}
    \put(26,24.5){\scriptsize $\lambda_2$}
    \put(23,30.5){\scriptsize $\lambda_{\mathcal{S}_\mathrm{O}}$}
    \put(74,0){\includegraphics[width=7.3cm]{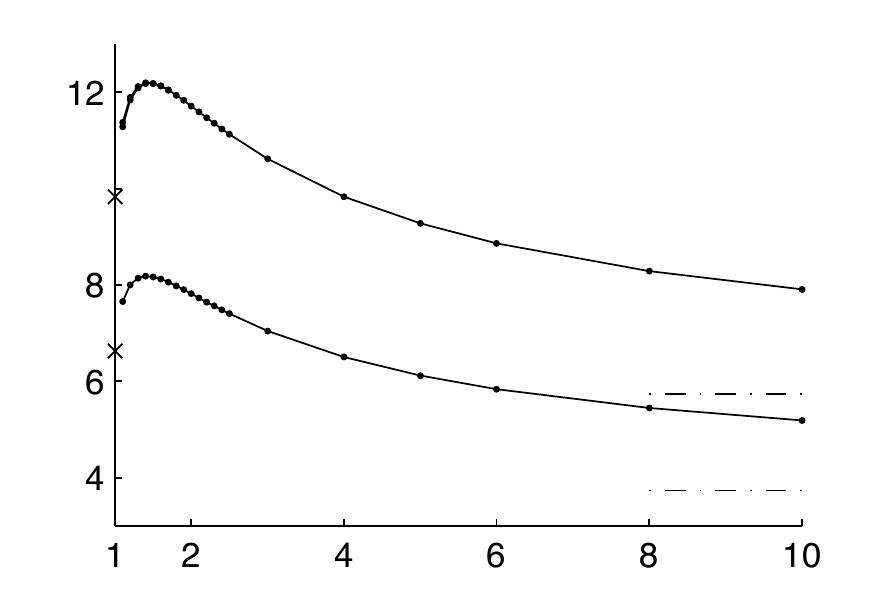}}
    \put(86,29){\includegraphics[height=2cm]{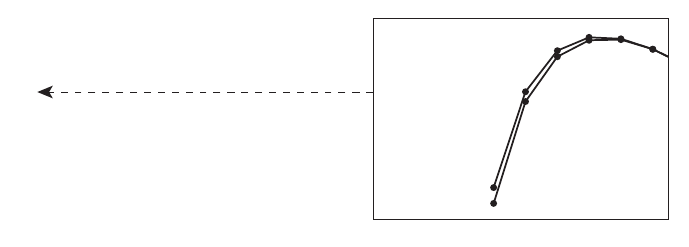}}
    \put(80.5,45.5){$\lambda^{1/p}$}
    \put(143,4){$p$}
    \put(111,12){$(\lambda_1)^{1/p}$}
    \put(111,24){$(\lambda_2)^{1/p}$}
    \put(118,42){\scriptsize  $(\lambda_{\mathcal{S}_\mathrm{O}})^{1/p}$}
    \put(127.5,33){\scriptsize $(\lambda_1)^{1/p}$}
    \put(74,19.5){\small $h_1(\Omega)$}
    \put(70,32){\small $h_1(\Omega^\mathrm{half})$}
    \put(141.5,7.5){\small $\Lambda_1(\Omega)$}
    \put(141.5,15.5){\small $\Lambda_2(\Omega)$}
  \end{picture}  
  \caption{Dependence of the numerically computed eigenvalues for the
    triangle with height 3/4 on $p$.}
  \label{fig:triangle3_lambda_plot}
\end{figure}
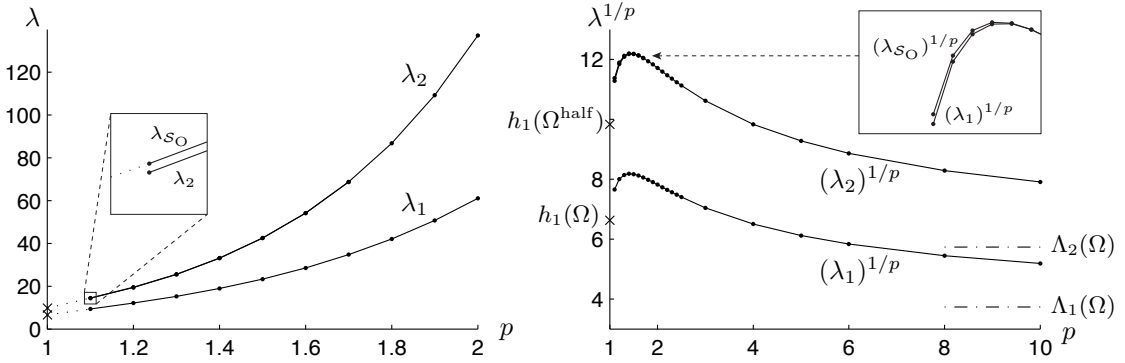

The dependence of the eigenvalues $\lambda_1$, $\lambda_2$, and
$\lambda_{\mathcal{S}_\mathrm{O}}$ on $p$ is further plotted in
Fig.~\ref{fig:triangle3_lambda_plot}. The following limits of
$\lambda_1$ and $\lambda_{\mathcal{S}_\mathrm{O}}$ as $p\to1$ and
those of $\lambda_1$ and $\lambda_2$ as $p\to\infty$ are also marked
in the figure and have these respective values:
\begin{align}
  h_1(\Omega)&=\textstyle\frac{2}{3}(2+\sqrt{13}+\sqrt{6\pi})
  \approx 6.631, &
  \Lambda_1(\Omega)&=\textstyle\frac{2}{3}(2+\sqrt{13})\approx 3.737, \notag\\
  h_1(\Omega^\mathrm{half})&=\textstyle\frac{2}{3}(5+\sqrt{13}+2\sqrt{3\pi})
  \approx 9.830, &
  \Lambda_2(\Omega)&=\textstyle\frac{2}{3}(5+\sqrt{13})\approx 5.737. \notag
\end{align}

\subsection{Equilateral triangle}\label{triangle_equi}
For isosceles triangles close but not equal to an equilateral triangle
a similar observation has been maded as for rectangles close but not
equal to the square: the symmetry properties of the second
eigenfunction $u_2$ change at a certain value $p\neq 2$. According to
the following computations, for an equilateral triangle this change
occurs at $p=2$ (as it does for the square).

Let
\begin{displaymath}\textstyle
  \Omega=\left\{(x_1,x_2)\in\R^2\,\Big|\, x_1\in\left(0,\frac{\sqrt{3}}{2}\right),
    |x_2|<\frac{1}{\sqrt{3}}\left(\frac{\sqrt{3}}{2}-x_1\right)\right\}
\end{displaymath}
be an equilateral triangle with side 1. It was discretized using
32,256 triangles. With the notation introduced in
\pref{eq:halftriangle_def} we can define
$\lambda_{2,\mathcal{S}_\mathrm{E}}$ and
$u_{2,\mathcal{S}_\mathrm{E}}$ as in Sec.~\ref{triangle1} and
$\lambda_{\mathcal{S}_\mathrm{O}}$ and $u_{\mathcal{S}_\mathrm{O}}$ as
in Sec.~\ref{triangle3}.

\begin{table}[h]
  \centering\vspace{1ex}
  \begin{tabular}[t]{|c|c|c|c|c|c|c|c|c|}
    \cline{1-4}\cline{6-9}
    $p$ & $\lambda_1$ & $\lambda_2(=\lambda_{2,\mathcal{S}_\mathrm{E}})$ & $\lambda_{\mathcal{S}_\mathrm{O}}$ & \qquad\qquad &
    $p$ & $\lambda_1$ & $\lambda_2(=\lambda_{\mathcal{S}_\mathrm{O}})$ & $\lambda_{2,\mathcal{S}_\mathrm{E}}$ \\\cline{1-4}\cline{6-9}
    \multicolumn{1}{c}{\vspace{\tableheadsep}} &
    \multicolumn{1}{c}{\rule{\tablecolwidthtreq}{0pt}} &
    \multicolumn{1}{c}{\rule{\tablecolwidthtreq}{0pt}} &
    \multicolumn{1}{c}{\rule{\tablecolwidthtreq}{0pt}} &
    \multicolumn{1}{c}{} & \multicolumn{1}{c}{} &
    \multicolumn{1}{c}{\rule{\tablecolwidthtreq}{0pt}} &
    \multicolumn{1}{c}{\rule{\tablecolwidthtreq}{0pt}} &
    \multicolumn{1}{c}{\rule{\tablecolwidthtreq}{0pt}} \\\cline{1-4}\cline{6-9}
    1.1 & 8.653 & 13.37 & 13.61 && 2.0 & 52.64 & 122.8 & 122.8 \\\cline{1-4}\cline{6-9}
    1.9 & 44.07 & 98.20 & 98.40 && 2.1 & 62.71 & 152.9 & 153.2 \\\cline{1-4}\cline{6-9}
    2.0 & 52.64 & 122.8 & 122.8 && 8.0 & $4.240\cdot10^5$ & $1.483\cdot10^7$ & $1.668\cdot10^7$ \\\cline{1-4}\cline{6-9}
  \end{tabular}\vspace{2ex}
  \caption{Eigenvalues for the equilateral triangle with side 1.}
  \label{tab:triangle_equi_lambda}
\end{table}

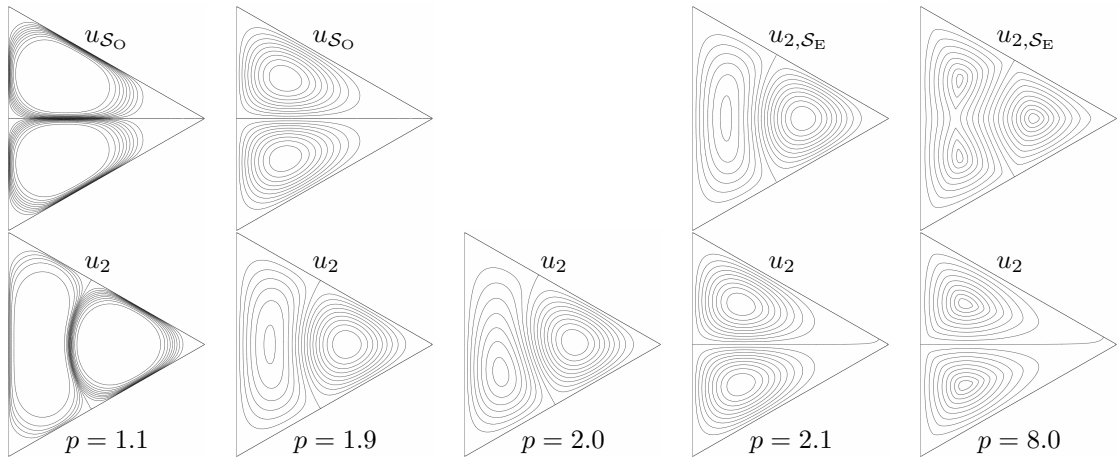
\begin{figure}[h]
  \centering
  \setlength{\unitlength}{1mm}
  \begin{picture}(150,60)
    \put(2,0){\includegraphics[width=26mm]{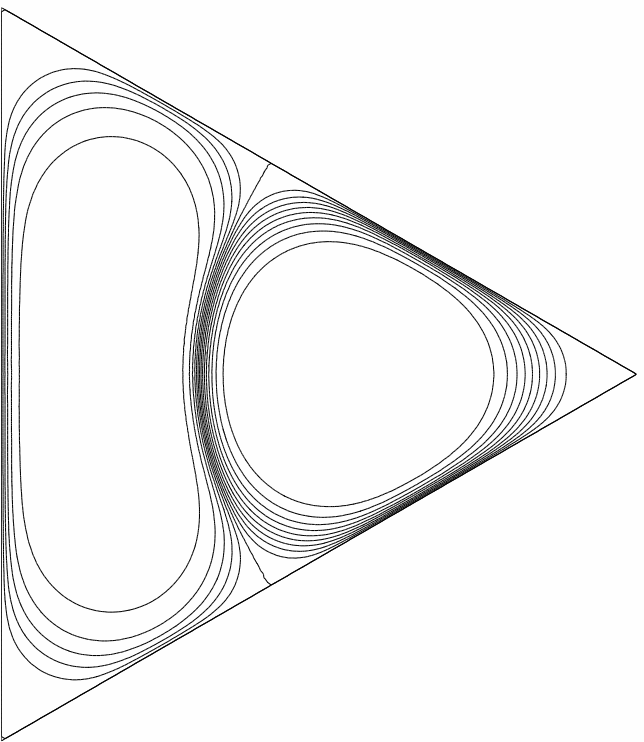}}
    \put(12,25){$u_2$}
    \put(2,30){\includegraphics[width=26mm]{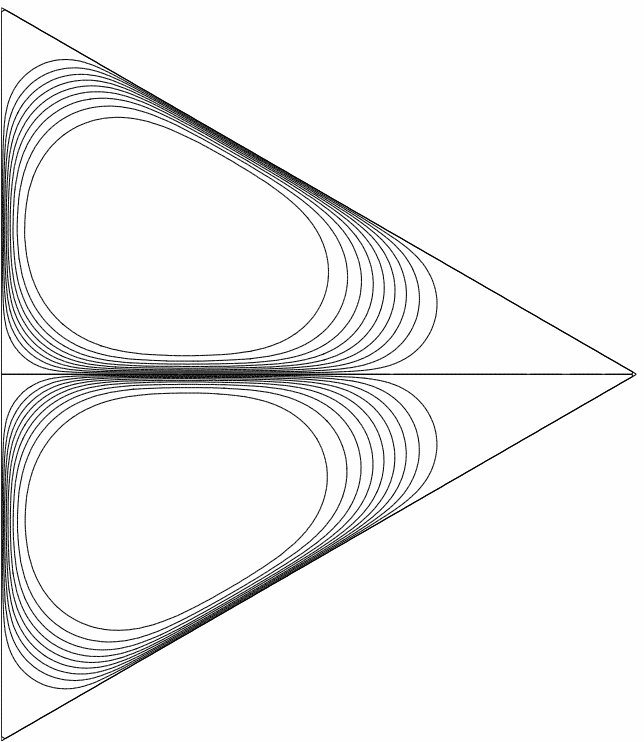}}
    \put(12,55){$u_{\mathcal{S}_\mathrm{O}}$}
    \put(0,0){\makebox(30,4){$p=1.1$}}
    \put(32,0){\includegraphics[width=26mm]{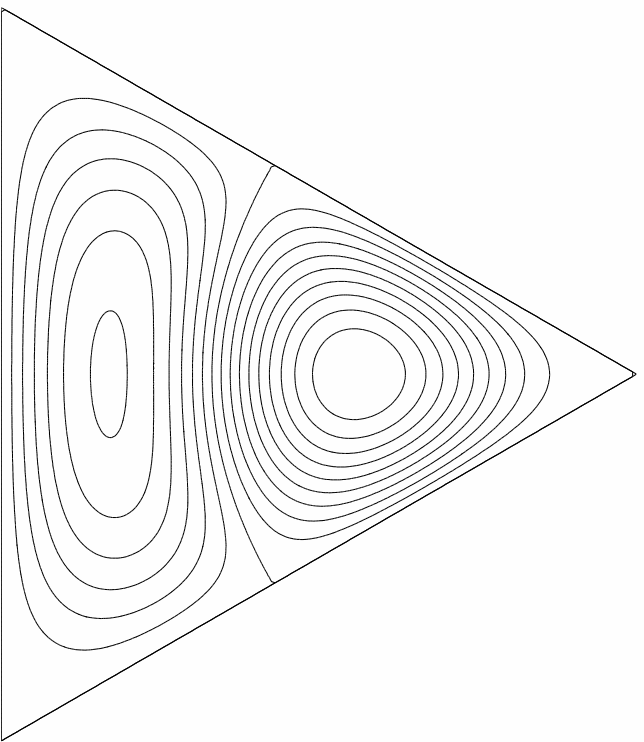}}
    \put(42,25){$u_2$}
    \put(32,30){\includegraphics[width=26mm]{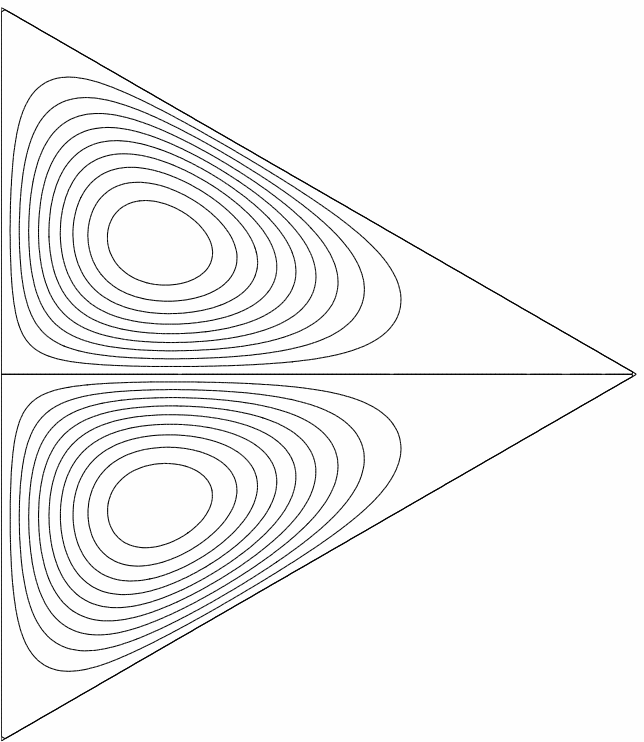}}
    \put(42,55){$u_{\mathcal{S}_\mathrm{O}}$}
    \put(30,0){\makebox(30,4){$p=1.9$}}
    \put(62,0){\includegraphics[width=26mm]{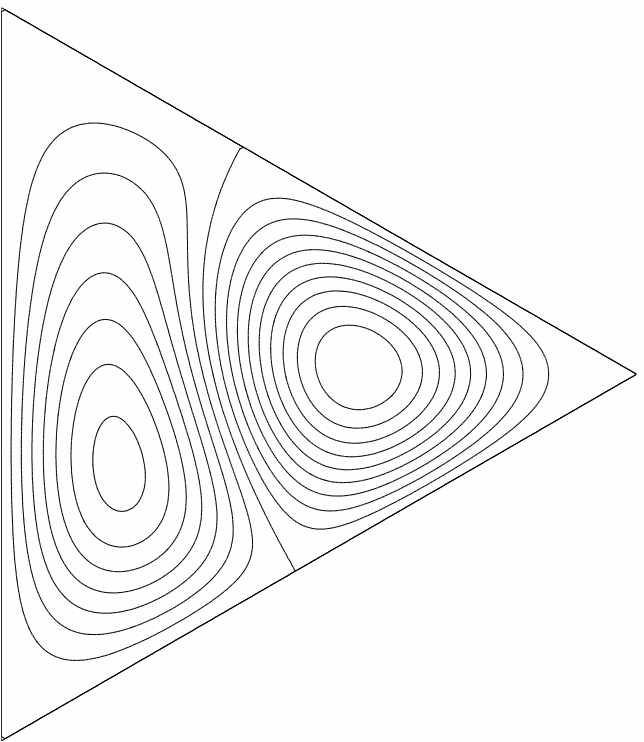}}
    \put(72,25){$u_2$}
    \put(60,0){\makebox(30,4){$p=2.0$}}
    \put(92,0){\includegraphics[width=26mm]{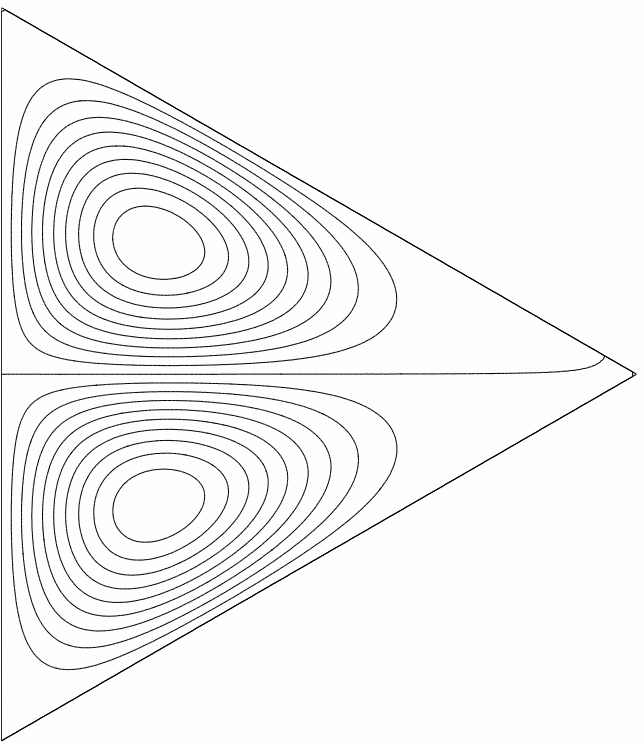}}
    \put(102,25){$u_2$}
    \put(92,30){\includegraphics[width=26mm]{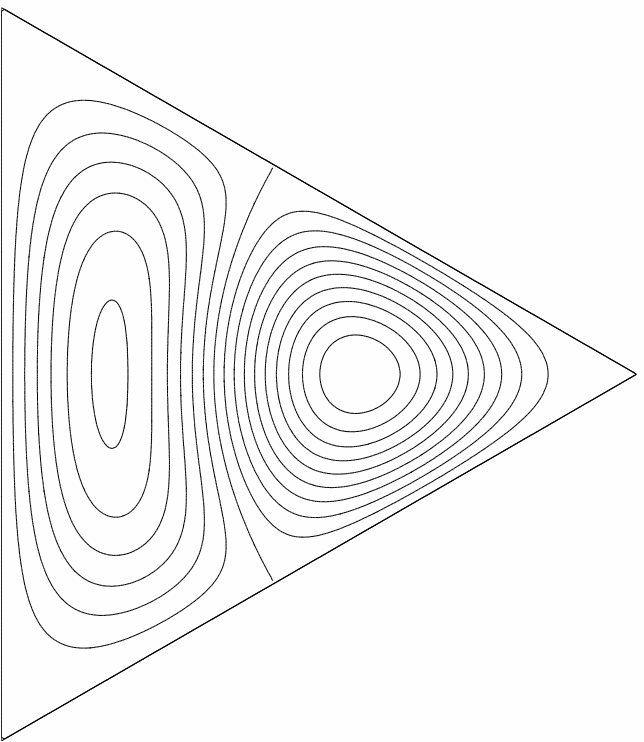}}
    \put(102,55){$u_{2,\mathcal{S}_\mathrm{E}}$}
    \put(90,0){\makebox(30,4){$p=2.1$}}
    \put(122,0){\includegraphics[width=26mm]{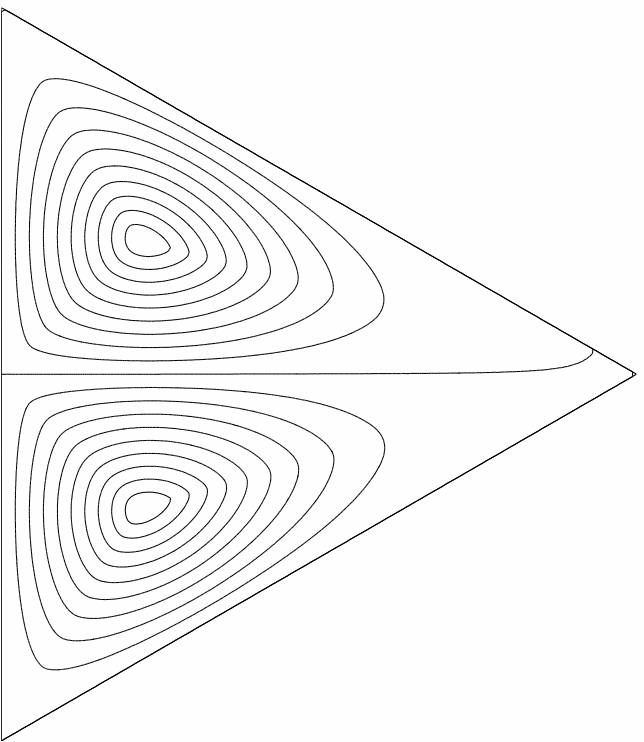}}
    \put(132,25){$u_2$}
    \put(122,30){\includegraphics[width=26mm]{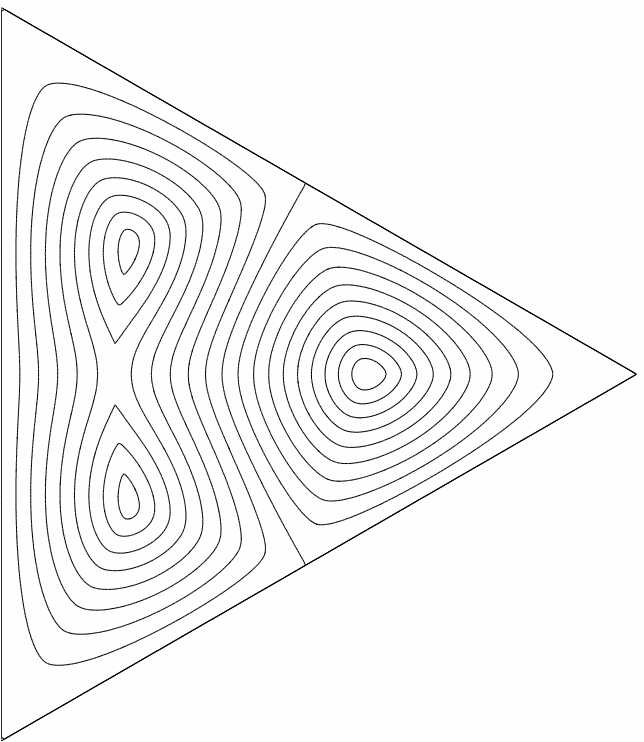}}
    \put(132,55){$u_{2,\mathcal{S}_\mathrm{E}}$}
    \put(120,0){\makebox(30,4){$p=8.0$}}
  \end{picture}
  \caption{The numerically computed eigenfunctions $u_2$,
    $u_{2,\mathcal{S}_\mathrm{E}}$ and $u_{\mathcal{S}_\mathrm{O}}$
    for the equilateral triangle.}
  \label{fig:triangle_equi_u2}
\end{figure}

Our numerical observations are summarized in
Table~\ref{tab:triangle_equi_lambda} and
Fig.~\ref{fig:triangle_equi_u2}: For $p<2$ the second eigenfunction
$u_2$ is even in $x_2$ while for $p>2$ it is odd (up to a rotation of
the triangle by $\pm 2\pi/3$). We note that the values $\lambda_2$
listed in the table were computed with no a priori assumptions on the
symmetry of $u$ and then compared to the computed values
$\lambda_{2,\mathcal{S}_\mathrm{E}}$ and
$\lambda_{\mathcal{S}_\mathrm{O}}$. The corresponding eigenfunctions
$u_2$ are in the bottom row of the figure. For $p=2$ the eigenspace
belonging to $\lambda_2$ is two-dimensional. The member of this
eigenspace $u_2$ to which CMPA converges depends on the initial path,
i.e., on the choice of the intermediate point $e_\mathrm{M}$. The
figure shows one such member. As already mentioned in
Sec.~\ref{triangle3}, it is an artifact of the numerical method that
for $u_2$ and $p=2.1, 8.0$ the zero contour line does not exactly
reach the vertex, where $u_2$ is rather flat.

\section{Remarks on the numerics}\label{sec:numerics_remarks}

\subsection{Dependence on the mesh parameter $h$}
Let $\mathcal{T}^h$ denote the set of all the triangles of a
triangulation of $\Omega^h$. The mesh parameter $h$ was introduced in
Sec.~\ref{sec:fem} as the (smallest) upper bound on the diameter of
the circumscribed circle for triangles of $\mathcal{T}^h$. In this
section the dependence of the computed values of $\lambda_1$ and
$\lambda_2$ on $h$ is investigated. The investigation is conducted for
one particular domain $\Omega$---the rectangular domain used for
computations in Sec.~\ref{rectangle}:
\begin{displaymath}
  \Omega=\left\{(x_1,x_2)\in\R^2\,\big|\, x_1\in(0,2), x_2\in(0,1.75)\right\}.
\end{displaymath}
Four discretizations of this domain are
used. Table~\ref{tab:triangulations} lists details about these
discretizations ordered by the number of triangles. Essentially, a
finer mesh was obtained from a courser one by placing a new vertex in
the middle of each triangle side of the old triangulation, in effect
dividing each triangle in four.
\begin{table}[h]
  \centering
  \begin{tabular}{|c|c|crc|}
    \hline
    $\mathcal{T}^h$ & $h$ & \multicolumn{3}{|c|}{number of triangles} \\\hline
    \multicolumn{3}{c}{\vspace{\tableheadsep}}\\\hline
    course & 0.079 & \rule{5mm}{0pt} & 4,832 & \\\cline{2-5}
    & 0.044 && 19,328 & \\\cline{2-5}
    & 0.022 && 77,312 & \\\cline{2-5}
    fine & 0.011 && 309,248 & \\\hline
  \end{tabular}\vspace{2ex}
  \caption{Triangulations used to discretize the rectangular domain $\Omega$.}
  \label{tab:triangulations}
\end{table}

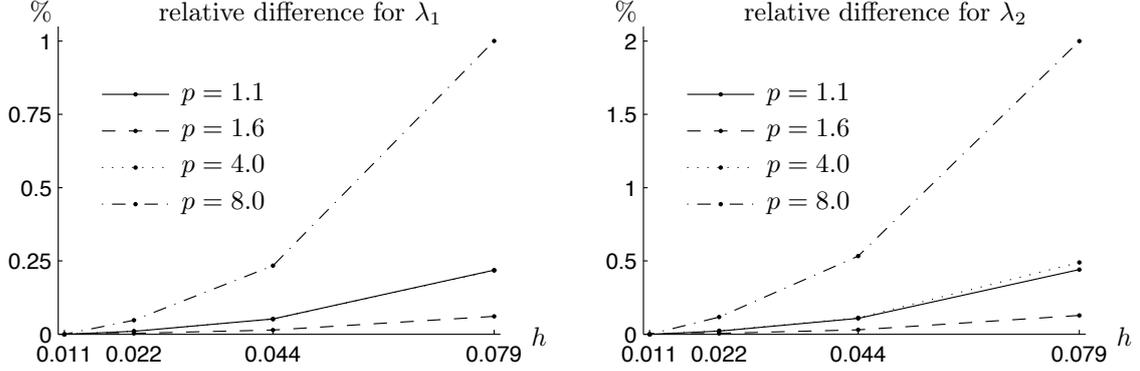
\begin{figure}[h]
  \centering
  \setlength{\unitlength}{1mm}
  \begin{picture}(150,50)(1,0)
    \put(-1,0){\includegraphics[width=7.5cm]{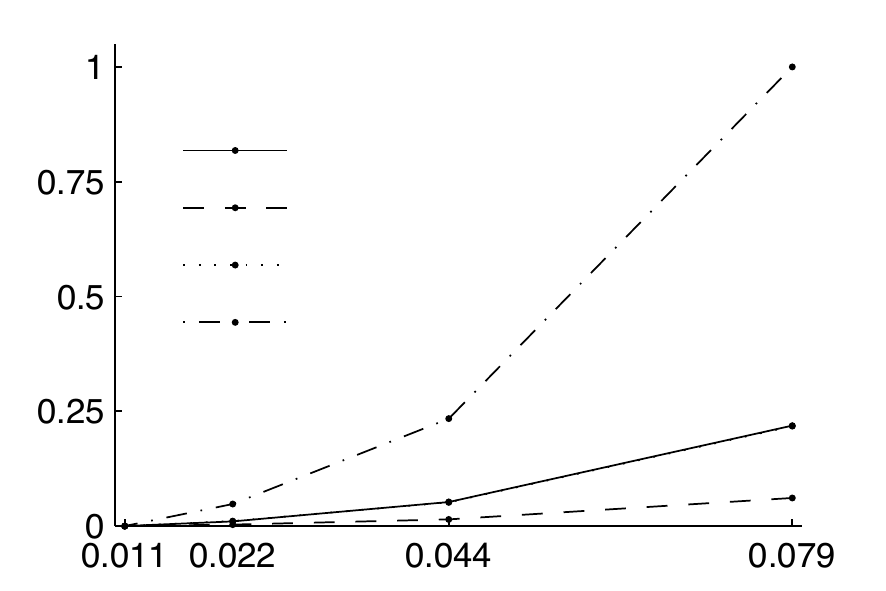}}
    \put(76,0){\includegraphics[width=7.5cm]{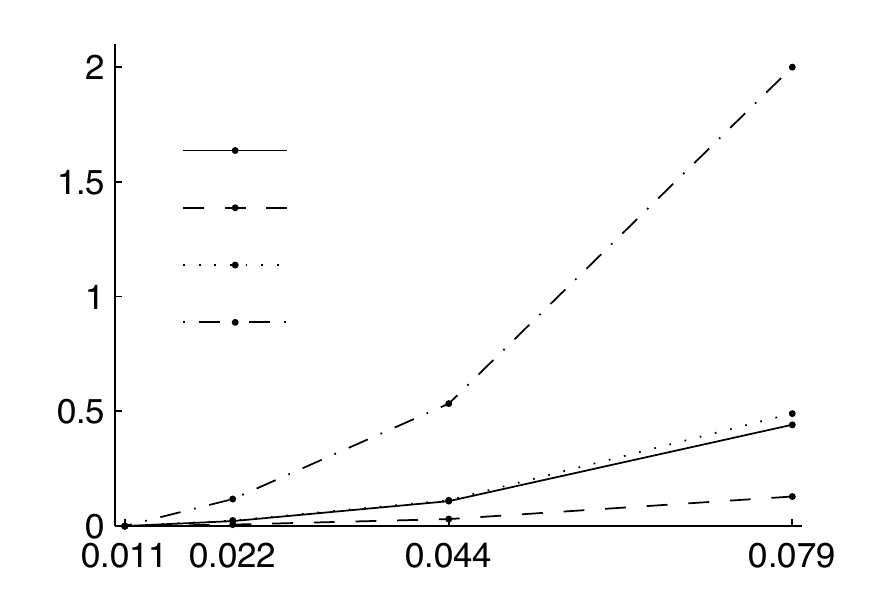}}
    \put(22,47){relative difference for $\lambda_1$}
    \put(5,47){\%}
    \put(71,4){$h$}
    \put(25,36.5){$p=1.1$}
    \put(25,31.7){$p=1.6$}
    \put(25,26.9){$p=4.0$}
    \put(25,22.1){$p=8.0$}
    \put(99,47){relative difference for $\lambda_2$}
    \put(82,47){\%}
    \put(148,4){$h$}
    \put(102,36.5){$p=1.1$}
    \put(102,31.7){$p=1.6$}
    \put(102,26.9){$p=4.0$}
    \put(102,22.1){$p=8.0$}
  \end{picture}  
  \caption{Relative difference
    $\frac{\lambda(\mathcal{T}^h)-\lambda(\mathcal{T}^{.011})}{\lambda(\mathcal{T}^{.011})}\cdot
    100\%$ for $\lambda_1$ (left) and $\lambda_2$ (right) computed on
    a triangulation $\mathcal{T}^h$ with respect to the finest
    triangulation $\mathcal{T}^{.011}$. For $\lambda_1$ the relative
    differences for $p=1.1$ and $p=4.0$ almost coincide and cannot be
    distinguished in the graph.}
  \label{fig:num_convergence_lambda}
\end{figure}

Table~\ref{tab:num_convergence_lambda} shows the values of $\lambda_1$
and $\lambda_2$ computed for the four triangulations characterized by
$h$ and for selected values of
$p$. Figure~\ref{fig:num_convergence_lambda} gives perhaps a more
telling picture: for each $p$ it displays relative differences
$(\lambda(\mathcal{T}^h)-\lambda(\mathcal{T}^{.011}))/\lambda(\mathcal{T}^{.011})$,
where $\lambda(\mathcal{T}^h)$ denotes the eigenvalue computed on the
triangulation $\mathcal{T}^h$. The finest triangulation
$\mathcal{T}^{.011}$ is used as a reference. We observe that the
largest differences occur for large $p$ (here $p=8.0$) and smallest
differences for $p$ close to 2. The differences are about twice as
large for $\lambda_2$ computed by CMPA compared to $\lambda_1$
computed by CDM.

\begin{table}[h]
  \centering
  \begin{tabular}{|c|@{}c@{}|c|c|c|c|c|c|c|c|c|c|c|}
    \cline{3-4}\cline{6-7}\cline{9-10}\cline{12-13}
    \multicolumn{1}{c}{} &\rule{.5ex}{0pt}& \multicolumn{2}{|c|}{$p=1.1$} & & \multicolumn{2}{|c|}{$p=1.6$} & & \multicolumn{2}{|c|}{$p=4.0$} & & \multicolumn{2}{|c|}{$p=8.0$} \\
    \cline{1-1}\cline{3-4}\cline{6-7}\cline{9-10}\cline{12-13}
    $h$ & & $\lambda_1$ & $\lambda_2$ & & $\lambda_1$ & $\lambda_2$ & & $\lambda_1$ & $\lambda_2$ & & $\lambda_1$ & $\lambda_2$ \\\cline{1-1}\cline{3-4}\cline{6-7}\cline{9-10}\cline{12-13}
    \multicolumn{13}{c}{\vspace{\tableheadsep}}\\\cline{1-1}\cline{3-4}\cline{6-7}\cline{9-10}\cline{12-13}
    0.079 && 2.5586 & 3.9507 && 4.2965 & 8.2462 && 14.884 & 91.306 && 50.005 & 2,234.8 \\\cline{1-1}\cline{3-4}\cline{6-7}\cline{9-10}\cline{12-13}
    0.044 && 2.5544 & 3.9376 && 4.2945 & 8.2381 && 14.859 & 90.962 && 49.625 & 2,202.7 \\\cline{1-1}\cline{3-4}\cline{6-7}\cline{9-10}\cline{12-13}
    0.022 && 2.5533 & 3.9342 && 4.2940 & 8.2361 && 14.853 & 90.881 && 49.533 & 2,193.6 \\\cline{1-1}\cline{3-4}\cline{6-7}\cline{9-10}\cline{12-13}
    0.011 && 2.5531 & 3.9334 && 4.2939 & 8.2356 && 14.851 & 90.861 && 49.510 & 2,191.0 \\\cline{1-1}\cline{3-4}\cline{6-7}\cline{9-10}\cline{12-13}
  \end{tabular}\vspace{2ex}
  \caption{Values of $\lambda_1$ and $\lambda_2$ computed for four different triangulations of the rectangular domain $\Omega$ and $p=1.1, 1.6, 4.0, 8.0$.}
  \label{tab:num_convergence_lambda}
\end{table}

\subsection{The Augmented Lagrangian Method}
As described in Sec.~\ref{sec:inverse} this method is used to solve
\pref{eq:constant_rhs} iteratively for a given right-hand side. The
Augmented Lagrangian $\mathcal{L}_r$ defined in
\pref{eq:augmented_lagrangian} depends on a parameter $r>0$. As
observed by the authors of \cite{GlowinskiMarroco} the algorithm is
not very sensitive to the choice of $r$ but the analysis of the
influence of $r$ on the behavior of the algorithm is complicated.

The choice of $r$ has an influence on the speed of convergence of the
algorithm and at the same time on how precise the found numerical
solutions can be. In general, for larger $r$ the algorithm seems to
converge faster but it is able to find only less precise
approximations of the solution. For our computations we tried various
values of $r$ first and then chose the one which seemed to give a
reasonable speed of convergence together with acceptable
residual. This value depended strongly on $p$ and also on the
particular domain $\Omega$. Table~\ref{tab:augmented_r} shows the
dependence of $r$ and of the number of iterations that the algorithm
needed on $p$. For each domain one value of $r$ was chosen from the
given range. Similarly, the number of iterations lay in the given
range. We can observe that for a small $p$ a large $r$ was needed, for
a large $p$ a smaller $r$. For $p$ close to 2 we could choose
$r\approx 1$. The number of iterations needed turned larger for $p$
farther from 2.

\begin{table}[t]
  \centering
  \begin{tabular}{|c|r@{\ --\ }l|r@{\ --\ }l|}
    \hline
    $p$ & \multicolumn{2}{c|}{range of $r$} & \multicolumn{2}{c|}{\# of iterations} \\\hline
    \multicolumn{5}{c}{\vspace{\tableheadsep}}\\\hline
    1.1 & $10^4$ & $10^7$ & 700 & 2,000 \\\hline
    1.2 & $500$ & $2,500$ & 500 & 1,000 \\\hline
    1.8 & $1$ & $1.5$ & 80 & 90 \\\hline
    3.0 & $0.3$ & $0.4$ & 200 & 300 \\\hline
    10.0 & $0.03$ & $0.1$ & 1,200 & 3,000 \\\hline
  \end{tabular}\vspace{2ex}
  \caption{The dependence of the approximate values of $r$ and numbers of iterations
    in the Augmented Lagrangian Method on $p$.}
  \label{tab:augmented_r}
\end{table}

We note that for values of $p$ smaller than 1.1 and larger than 10
(the particular value also depended on the domain) we were not able to
find $r$ giving satisfactory results for our implementation of the
Augmented Lagrangian Method in conjunction with CMPA.

\subsection{CDM and CMPA}
In both the Constrained Descent Method and the Constrained Mountain
Pass Algorithm the measure of convergence is $\|w_u\|$, the
$\sobspace$-norm of the descent direction evaluated at the
approximation $u$ of the eigenfunction which is being computed. The
smallest achieved value depended on the algorithm, on $p$, and in case
of CMPA also on the fact whether there lies another critical point not
far from $u$. For CDM the order of this value was between $10^{-5}$
and $10^{-8}$, for CMPA between $10^{-3}$ and $10^{-7}$. The number of
iteration of CDM was approximately between 10 and 30. The number of
iterations of CMPA varied, it depended on the shape of the initial
path and on $p$, and was anywhere between 100 and 3,000.

\section{Conclusion}\label{sec:conclusion}
In this work a concrete application of the variational numerical
methods of \cite{Ho1} in a Banach space was presented. In particular,
one possible choice of the descent direction required by these methods
was proposed, implemented and tried in computations in the setting of
the Sobolev space $\sobspace$. The computations yielded approximations
of the smallest two Dirichlet eigenvalues and the corresponding
eigenfunctions of the $p$-Laplace operator on several planar domains
for $p$ ranging from 1.1 to 10. This relatively large range made it
possible to study the change of symmetry of the second eigenfunction
with varying $p$ on different domains which was first observed in
\cite{YaoZhou1} for the square and $p$ not far from 2. The computed
eigenvalues seem to agree with the asymptotic behavior known from
theory for $p\to 1$. Our range of $p$ seems to be too small, however,
in order to clearly observe the asymptotic behavior of the eigenpairs
as $p\to\infty$.

Numerical experiments were conducted for the following domains: the
disk, rectangles, and isosceles triangles. We summarize the main
observations about the symmetry of the second eigenfunction $u_2$. For
the disk it was observed that $u_2$ has a straight nodal line dividing
the disk into halves for the whole range of $p$.

For rectangles which are not a square and for small $p$ the second
eigenfunction is odd about its nodal line which is straight and
connects the midpoints of the longer sides. After $p$ crosses some
value $p_0>2$ there are two second eigenfunctions which are mirror
images of each other. Their nodal line is not straight but still
connects the two longer sides.

For the square and $p\neq 2$ there are two second eigenfunctions which
are images of each other under rotations by $\pi/2$ about the center
of the square. For $p<2$ their nodal line is straight and connects the
midpoints of the opposite sides. For $p>2$ the nodal line is a
diagonal of the square. For $p=2$ there are two linearly independent
second eigenfunctions.

The symmetry observations for triangles are based on the family of
isosceles triangles with vertices $(0,-1/2)$,
$(0,1/2)$, and $(\ell,0)$ with base 1 and height $\ell>0$
which are symmetric about the $x_1$-axis. For those shorter than the
equilateral triangle and for small $p$ there are two asymmetric
second eigenfunctions (up to scaling) which are symmetry images of
each other. Their nodal line connects the base with one side of the
triangle. After $p$ crosses some value $p_0<2$ there is only one
eigenfunction $u_2$. It is odd about its nodal line which is straight
and connects the middle of the base with the opposite vertex
(symmetry $\mathcal{S}_\mathrm{O}$).

For triangles longer than the equilateral triangle and for small $p$
there is one eigenfunction $u_2$, it is even about the $x_1$-axis
(symmetry $\mathcal{S}_\mathrm{E}$) and its nodal line connects the
two sides of the triangle. After $p$ crosses some value $p_0>2$ there
are two asymmetric second eigenfunctions which are symmetry images of
each other. Their nodal line still connects the two sides of the
triangle.

For the equilateral triangle and $p\neq 2$ there are three second
eigenfunctions which are images of each other under rotations of
the triangle about its midpoint by $\pm 2\pi/3$. For
$p<2$ their nodal line connects two sides of the triangle and they
have even symmetry about the height coming from the third side. For
$p>2$ the nodal line of the second eigenfunctions follows a height of
the triangle and the eigenfunctions have odd symmetry about this
height. For $p=2$ there are two linearly independent second
eigenfunctions.

Figure~\ref{fig:triangle3_u2extra} indicates that our numerical
methods could be used for finding some higher eigenfunctions and
perhaps for a continuation in $\ell$ to observe the connection between
these eigenfunctions and the second eigenfunctions for the equilateral
triangle. This lies however beyond the scope of this paper.

\appendix
\section{}
Here we give a proof of some claims used in Sections~\ref{sec:descent}
and \ref{sec:CDM}. A subindex notation will be used for general
sequences and does not refer to the enumeration of eigenfunctions and
eigenvalues in this section.
\begin{plapnum_lemma}\label{lemma1}
  Let $(B,\|\cdot\|)$ be a reflexive Banach space with a strictly
  convex norm, $I,J\in C^1(B,\R)$ be two continuously Fr\'echet
  differentiable functionals, and $u$ be a point in $B$ with $J(u)=1$
  which is not a critical point of $I$ with respect to $S:=\{v\in
  B\,|\,J(v)=1\}$. Then the problem
  \begin{displaymath}
    \text{minimize } L(w):=\dualpairing{I'(u)}{w} \quad \text{ subject to } \quad w\in\mathcal{C}:=
    \left\{v\in B\,\left|\, \dualpairing{J'(u)}{v}=0 \text{ and } \|v\|=1\right.\right\}
  \end{displaymath}
  has a unique solution.
\end{plapnum_lemma}

\begin{proof}
  First, we show that $L$ has a negative infimum on $\mathcal{C}$: $L$
  is bounded below on $\mathcal{C}$ by $-\|I'(u)\|_\ast$. It attains
  negative values on $\mathcal{C}$ if there exists $w\in\mathcal{C}$
  such that $L(w)\neq 0$. But if $L\equiv 0$ on $\mathcal{C}$, then we
  would have
  \begin{displaymath}
    \dualpairing{J'(u)}{w}=0 \quad\Rightarrow\quad \dualpairing{I'(u)}{w}=0 \qquad \forall w\in B
  \end{displaymath}
  which would imply existence of $\alpha\in\R$ such that $I'(u)-\alpha
  J'(u)=0$. This is not possible since $u$ is not a critical point of
  $I$ with respect to $S$.

  Let $\{w_n\}\subset\mathcal{C}$ be a minimizing sequence of $L$,
  i.e.,
  \begin{displaymath}
    L(w_n)\to\inf_\mathcal{C} L \in(-\infty,0)\qquad \text{as } n\to\infty.
  \end{displaymath}
  Since this sequence is bounded by 1, the reflexivity of $B$ implies
  existence of a subsequence (still denoted $\{w_n\}$) which converges
  weakly to some $w\in B$ such that $\|w\|\leq 1$. Since $I'(u)$ and $J'(u)$ are
  continuous linear functionals, we obtain
  \begin{displaymath}
    L(w_n)\to L(w)\quad\text{and}\quad \dualpairing{J'(u)}{w_n}\to \dualpairing{J'(u)}{w}
    \qquad \text{as } n\to\infty.
  \end{displaymath}
  This means that $L(w)=\inf_\mathcal{C}L$ and $\dualpairing{J'(u)}{w}=0$.

  To prove that $w$ is a minimizer it remains to show that
  $\|w\|=1$. If $\|w\|<1$, then $\tilde w:=w/\|w\|$ belongs to
  $\mathcal{C}$ and
  \begin{displaymath}
    L(\tilde w)=\frac{L(w)}{\|w\|}<L(w)
  \end{displaymath}
  because $L(w)<0$. But this is a contradiction with the minimality of
  $L(w)$.

  To show uniqueness let $w_1$ and $w_2$ be both minimizers. For
  $\bar w:=\frac{1}{2}(w_1+w_2)\left/\left\|\frac{1}{2}(w_1+w_2)\right\|\right.$
  we obtain
  \begin{displaymath}
    \bar w\in\mathcal{C}\qquad\text{and}\qquad \min_\mathcal{C}L\leq L(\bar w)=
    \frac{\min_\mathcal{C}L}{\left\|\frac{1}{2}(w_1+w_2)\right\|}.
  \end{displaymath}
  Since the minimum is negative, this and the triangle inequality
  imply
  \begin{displaymath}
    \textstyle 1 \leq \left\|\frac{1}{2}w_1+\frac{1}{2}w_2\right\| \leq
    \frac{1}{2}\|w_1\|+\frac{1}{2}\|w_2\| = 1.
  \end{displaymath}
  Hence equality holds in the above inequalities and the strict
  convexity of the norm implies $w_1=w_2$.
\end{proof}

\begin{plapnum_lemma}\label{lemma2}
  Let $I$ and $J$ be defined by \pref{eq:IJ} and $S$ by
  \pref{eq:S}. Further, let $u\in S$ and
\begin{equation}\label{eq:descent_direction1}
  w_u:=-u+\frac{1}{\int_\Omega |u|^{p-2}u \,v_u\,dx}\ v_u,\qquad
  \text{where }v_u:=(-\Delta_p)^{-1}\left(|u|^{p-2}u\right).
\end{equation}
Then $\dualpairing{I'(u)}{w_u}\leq 0$. Equality holds if and only if
$u$ is a critical point of $I$ with respect to $S$ which is the case
if and only if $w_u=0$.
\end{plapnum_lemma}
The proof of this lemma is based on the following inequality which is
a direct consequence of the Cauchy-Schwarz and H\"older
inequalities. Its proof is therefore omitted.
\begin{plapnum_aux_lemma}
  Let $f,g\in\sobspace$, $f\neq 0$. Then
  \begin{displaymath}
    \int_\Omega|\nabla f|^{p-2}\nabla f\nabla g\,dx \leq \|f\|^{p-1}\|g\|.
  \end{displaymath}
  Equality holds if and only if there exists $\nu\geq 0$ such that
  $\nu f=g$.
\end{plapnum_aux_lemma}

\begin{proof}[Proof of Lemma \ref{lemma2}]
  We observe that
  \begin{equation}\label{eq:normvu}
    \int_\Omega |u|^{p-2}u \,v_u\,dx = \int_\Omega (-\Delta_p v_u)v_u\,dx = \|v_u\|^p.
  \end{equation}
  By the definition of $w_u$, \pref{eq:normvu} and the auxiliary lemma
  we obtain
  \begin{equation}
    \label{eq:firstineq}
    \begin{aligned}
      \dualpairing{I'(u)}{w_u} &= -\|u\|^p +
      \frac{1}{\|v_u\|^p}\int_\Omega
      |\nabla u|^{p-2}\nabla u\nabla v_u\,dx \\
      &\stackrel{(\ast)}{\leq} -\|u\|^p +
      \frac{1}{\|v_u\|^p}\|u\|^{p-1}\|v_u\| = \left(-1 +
        \frac{1}{\|v_u\|^{p-1}\|u\|}\right)\|u\|^p.
    \end{aligned}
  \end{equation}
  Using $u\in S$, testing the equation $-\Delta_p v_u=|u|^{p-2}u$ by
  $u$, and applying the auxiliary lemma yields
  \begin{equation}
    \label{eq:secondineq}
    1 = \int_\Omega|u|^p\,dx = \int_\Omega|\nabla v_u|^{p-2}\nabla v_u\nabla u\,dx
    \stackrel{(\ast\ast)}{\leq} \|v_u\|^{p-1}\|u\| .
  \end{equation}
  By combining \pref{eq:firstineq} and \pref{eq:secondineq} we
  conclude that $\dualpairing{I'(u)}{w_u}\leq 0$. Equality holds if
  and only if equality holds in $(\ast)$ and $(\ast\ast)$. According
  to the auxiliary lemma this is the case if and only if $\nu u=v_u$ for
  some $\nu>0$. Finally, we argue that the following are equivalent:
  \begin{enumerate}[label=\rm (\alph*)]
  \item $\nu u=v_u$ for some $\nu>0$,
  \item $u$ is a critical point of $I$ with respect to $S$,
  \item $w_u=0$.
  \end{enumerate}
  Statement (a) is equivalent to $\nu^{p-1}(-\Delta_p u)=|u|^{p-2}u$ and
  hence to (b). If (a) holds, then $\int_\Omega |u|^{p-2}u
  \,v_u\,dx=\nu$ because $u\in S$. Hence $w_u=-u+\frac{1}{\nu}v_u=0$ and
  (c) holds, too. It is obvious that (c) implies (a).
\end{proof}

Before stating the next proposition we recall some known results (let
$p,q\in(0,\infty)$, $\frac{1}{p}+\frac{1}{q}=1$):
\begin{enumerate}[label=\rm (\roman*)]
\item The $p$-Laplace operator $-\plap:\sobspace\to W^{-1,q}(\Omega)$
  is uniformly continuous on bounded sets.
\item The mapping $u\mapsto|u|^{p-2}u:\sobspace\to W^{-1,q}(\Omega)$ is
  compact and uniformly continuous on bounded sets.
\item The inverse $p$-Laplace operator $(-\plap)^{-1}:
  W^{-1,q}(\Omega)\to\sobspace$ is uniformly continuous on bounded
  sets.
\end{enumerate}
Both claims (i) and (ii) follow from standard inequalities found,
e.g., in \cite[Lemmas 5.3 and 5.4]{GlowinskiMarroco}. The compactness
in (ii) follows from the compact embedding of $\sobspace$ in
$L^p(\Omega)$. Claim (iii) follows from standard inequalities found,
i.e., in \cite[Propositions 5.1 and 5.2]{GlowinskiMarroco}.

\begin{plapnum_proposition}\label{proposition1}
Let $I$ and $J$ be defined by \pref{eq:IJ} and $S$ by
  \pref{eq:S}. The initial value problem
  \begin{equation}
  \label{eq:cdm_ivp1}
  \frac{d}{d t}u(t)=w_{u(t)},\qquad u(0)=e_0\in S
\end{equation}
with $w_u$ defined in \pref{eq:descent_direction1} has a unique solution
$u(t)\in S$ defined for $t\in(0,\infty)$. There exists a critical
point $u\in S$ of $I$ with respect to $S$ and a sequence
$\{t_n\}_{n=1}^\infty$ such that $\lim_{n\to\infty}t_n=\infty$ and
$\lim_{n\to\infty}u(t_n)=u$ in $\sobspace$.
\end{plapnum_proposition}
\begin{proof}
  The proof of existence of a solution and its uniqueness follows the
  same lines as the proof of Lemma 5 in \cite{Ho1}. Hence we focus on
  establishing the existence of the sequence $\{t_n\}$.

  Since $0\leq
  I(u(T))=I(e_0)+\int_0^T\dualpairing{I'(u(t))}{w_{u(t)}}dt$ for $T>0$
  and the integrand is non-positive, we obtain
  $\int_0^\infty\left|\dualpairing{I'(u(t))}{w_{u(t)}}\right|dt\leq
  I(e_0)$. Hence there exists a sequence $\{t_n\}_{n=1}^\infty$ with
  $\lim_{n\to\infty}t_n=\infty$ such that for $u_n:=u(t_n)$ and
  $w_n:=w_{u(t_n)}$ it holds:
  \begin{equation}\label{eq:convergence_Iprimewn}
    \dualpairing{I'(u_n)}{w_n}\to 0\quad\text{for } n\to\infty.
  \end{equation}

  We recall that by \pref{eq:descent_direction1} and \pref{eq:normvu} we have
  \begin{equation}
    \label{eq:seq_wn}
    w_n=-u_n+\frac{1}{\|v_n\|^p}v_n,\quad\text{where }v_n:=(-\plap)^{-1}\left(|u_n|^{p-2}u_n\right).
  \end{equation}
  We observe that $\{u_n\}$ is a bounded sequence, hence it converges
  weakly to some $u\in\sobspace$ along a subsequence which we again
  denote $\{u_n\}$. From the compactness of the map
  $u\mapsto|u|^{p-2}u$ and the continuity of the inverse $p$-Laplacian
  it follows that
  \begin{equation}
    \label{eq:convergence_vn}
    v_n\to v:=(-\plap)^{-1}\left(|u|^{p-2}u\right)\qquad \text{strongly in }\sobspace\text{ as }n\to\infty.
  \end{equation}

  Equation \pref{eq:seq_wn} and the fact that $w_n\in T_{u_n}S$ imply
  \begin{equation}\label{eq:wn_in_tangentspace}
    \dualpairing{-\plap(w_n+u_n)}{w_n}=\frac{1}{\|v_n\|^{p(p-1)}}\int_\Omega|u_n|^{p-2}u_n w_n\,dx=0.
  \end{equation}
  Combining \pref{eq:convergence_Iprimewn} and \pref{eq:wn_in_tangentspace} yields
  \begin{equation}
    \label{eq:convergence_to_zero}
    \int_\Omega\left(|\nabla(w_n+u_n)|^{p-2}\nabla(w_n+u_n) -
      |\nabla u_n|^{p-2}\nabla u_n\right)\nabla w_n\,dx \to 0\quad\text{for } n\to\infty.
  \end{equation}
  On the other hand standard estimates \cite[Propositions 5.1 and
  5.2]{GlowinskiMarroco} state that
  \begin{align}
    \int_\Omega\big(|\nabla(w_n+u_n)|^{p-2}\nabla(w_n+u_n) -
      &|\nabla u_n|^{p-2}\nabla u_n\big)\nabla w_n\,dx \notag \\
      &\geq \delta\,\frac{\|w_n\|^2}{(\|w_n+u_n\|+\|u_n\|)^{2-p}} &&\text{for }1<p\leq 2, \\
      &\geq \frac{1}{2^{p-2}}\|w_n\|^p &&\text{for }2\leq p,
  \end{align}
  where $\delta>0$ is a constant which does not depend on $w_n$ and
  $u_n$. These inequalities and \pref{eq:convergence_to_zero} imply
  \begin{equation}
    \label{eq:wn_to_zero}
    w_n\to 0 \qquad \text{strongly in }\sobspace\text{ as }n\to\infty.
  \end{equation}
  This and \pref{eq:seq_wn} in turn imply that $\{u_n\}$ converges
  strongly to $u$ and that
  \begin{equation}
    \label{eq:u_critical_point}
    u=\frac{1}{\|v\|^p}(-\plap)^{-1}\left(|u|^{p-2}u\right),
  \end{equation}
  which means that $u$ is a critical point of $I$ with respect to $S$.
\end{proof}

\begin{plapnum_remark}
  To better understand the implications of the choice of the descent
  direction we remark how the proof of the proposition would change if
  we used the steepest descent direction instead of the descent
  direction given by \pref{eq:descent_direction1}. Up to normalization
  the steepest descent direction $w$ defined by
  \pref{eq:steepestdescentBanach} can be written as the solution of
  \begin{displaymath}
    -\plap w=\plap u+\alpha|u|^{p-2}u
  \end{displaymath}
  for a suitable $\alpha$. Testing this equation by $w$ and using
  $w\in T_u S$ yields
  $\|w\|^p=-\frac{1}{p}\dualpairing{I'(u)}{w}$. Hence equation
  \pref{eq:convergence_Iprimewn} would directly imply $w_n\to 0$. We
  can write
  \begin{displaymath}
    0\leftarrow \|w_n\|^{p-1}=\|-\plap w_n\|_\ast = \left\|-\plap u_n -\alpha_n|u_n|^{p-2}u_n\right\|_\ast=
    \frac{1}{p}\left\|I'(u_n)-\alpha_n J'(u_n)\right\|_\ast,
  \end{displaymath}
  where $\|\cdot\|_\ast$ denotes the norm in the dual space
  $W^{-1,q}(\Omega)$. If we define
  $\|I'|_{S_u}\|:=\inf_{\alpha\in\R}\|I'(u)-\alpha J'(u)\|_\ast$ as in
  \cite{Bon}, \cite{Ho1}, then we would obtain $\|I'|_{S_{u_n}}\|\to
  0$. The Palais-Smale condition under constraints which was
  formulated in \cite{Bon} states in its simplified form that if
  $\{I'(u_n)\}$ is bounded and $\|I'|_{S_{u_n}}\|\to 0$ then $\{u_n\}$
  possesses a convergent subsequence. In
  \cite{CuestadeFigueiredoGossez} it was shown that this condition
  holds in our setting. Hence the choice of the steepest descent
  direction would yield a more ``classical'' proof of the proposition.
\end{plapnum_remark}

\bibliographystyle{amsplain} \bibliography{./plaplacenum_listofpapers}

\end{document}